\documentclass[12pt]{article}

\usepackage{amssymb}

\title{A System of Dependent Types, with an Implementation and a Philosophy}

\author{M. Randall Holmes}

\date{10/22/2016 -- needs extensive updates to embedded Lestrade books.}

\begin{document}

\maketitle

\newpage

\tableofcontents

\newpage

\section{Introduction}

This paper describes the development of some thoughts on philosophy of mathematics, a system of dependent type theory (for the general context of dependent type theories, \cite{barendregt} is a nice reference), and a computer implementation of these ideas, in parallel.  The computer implementation is at least currently nicknamed ``Lestrade" (to go with the author's ``Holmes"):  we hope that we may be forgiven.

The ideas on philosophy of mathematics aim at justification of the current practice of classical mathematics in a manner which may seem more likely to motivate a constructive or even finitist view.  That this can be done is in our view interesting.   The philosophy is certainly adaptable to a frankly constructive view.  All infinities invoked in this scheme are in a suitable sense merely potential:  the system is Aristotelean in its fashion.

The system of dependent types is of a familiar sort motivated by the Curry-Howard isomorphism (see \cite{curryfeys},\cite{howard},\cite{debruijntypes}) and the implementation is recognizably a variant of the system Automath of de Bruijn (for which the omnibus reference is the useful volume \cite{automathbook}; a survey publication by de Bruijn is \cite{thelanguage};  a modern implementation is Freek Wiedijk's \cite{freek}).  There are various modern relations of Automath (source and examples of use of which are more readily available and on a larger scale) which would be more or less distant cousins of Lestrade, such as Coq  (\cite{coq}).  The recent introduction of the ability to introduce rewrite rules by constructing or exhibiting terms of suitable types to justify them, the rewrite rules then being used to justify type equations and rewrite terms, arguably gives us not only a type checking system but a programming language capable of acting on all types of mathematical object and proof:  in particular, it appears to implement the style of programming with rewrite rules described by the author in \cite{holmesrewriting}.   Much more would need to be done to make Lestrade an actual programming environment!

The source code for Lestrade will always be available at \cite{holmessource}:  a recent version is appended to this document as an appendix.

\newpage

\section{A model of mathematical activity}

We present a model of what we might take a mathematician to be doing.

This will be recapitulated with more concrete detail in the description of the implementation of the syntax and semantics of the Lestrade system in section 3.  Most comparisons of Lestrade to its precursor Automath will be deferred to section 3, though some do occur here.

Objects that the ideal mathematician might consider are of various sorts.  Among these sorts are the sort {\tt prop} of propositions, and for each proposition $p$,
a sort ({\tt that} $p$) inhabited by proofs of $p$.

We are of course thinking of the view that mathematical proofs are themselves mathematical objects, which is embodied in the Curry-Howard isomorphism.  In the usual presentation of the Curry-Howard isomorphism, the proof of a conjunction
$A \wedge B$ is to be thought of as a pair consisting of a proof of $A$ and a proof of $B$ and the  proof of an implication is to be thought of as a function taking proofs of $A$ to proofs of $B$.
A proof of $\neg A$ is to be thought of as a function taking proofs of $A$ to proofs of $\perp$, the absurd.  A proof of a universally quantified sentence $(\forall x \in D:A(x))$ is to be thought of as a function taking elements $x$ of the domain $D$ to proofs of $A(x)$ (note that the type of the output will depend on the input).  What we will do will differ from this to some extent:  but this is the general idea of our approach.

We are agnostic as to whether the objects the mathematician talks about in his propositions and proofs are typed or untyped, so we provide support for both approaches
(and we suggest that both parts of this machinery might be useful).  We provide a sort {\tt obj} of untyped mathematical objects.  We provide a sort {\tt type} inhabited
by entities which we call ``type labels", and for each $\tau$ of sort {\tt type}, a corresponding sort ({\tt in} $\tau$):  we refer to objects of this sort as objects of type $\tau$.\footnote{There is a tension here between the word ``sort" and the word ``type":  we call the sorts of Lestrade itself ``sorts", and the special sorts which it uses to represent types of object considered by the ideal mathematician we call ``types":  an object of sort ({\tt in} $\tau$) is an object of type $\tau$.  Of course we may slip and say type when we mean sort, and we generally say ``type" for sorts in Automath itself, as the Automath workers did.}.

It should be noted that though we intend
the terminology to suggest that sorts ({\bf that} $p$) are inhabited by proofs and sorts ({\bf in} $\tau$) are inhabited by typed mathematical objects, the logical framework actually treats  {\bf prop / that} on the one hand and {\bf type / in} on the other in exactly the same way.   In an actual Lestrade book the general axioms for logic on the one hand and manipulation of types on the other are not likely to be exactly parallel.

These are the sorts of {\em entities\/} that we postulate.  Of course, if we postulate additional objects of type {\tt prop} or {\tt type} we simultaneously postulate new sorts of entity.

In addition, we have {\em abstractions\/}, functions taking given objects to new entities.

We introduce abstractions (and in special cases entities) by two processes:  {\em construction} and {\em definition}.

The process of construction is the local version of the introduction of axioms.  This amounts to postulating a function which acts on a concrete finite list of objects of given
sorts to obtain a new entity of a given sort.  We declare a sequence of variables, giving a sort for each variable, and give an entity sort for the output.  A subtlety is that our system of sorts admits dependent sorts:  later variables in the sequence of arguments of the construction may have sorts that depend on earlier variables in the sequence,  and the sort of the output may depend on the input variables.  For example,
a proof that for all $x$ of type $\tau$, $\phi(x)$ (where $\phi$ is a previously given construction taking a variable $x$ of type $\tau$ to output of type {\tt prop}, the natural typing for a predicate) can be thought of as a function with first argument $x$ of sort ({\tt in} $\tau$), and output $xx$ of sort ({\tt that} $\phi(x))$.  Note that the sort of the output depends on the value of the first argument.  Even more subtly, the universal quantifier over type $\tau$ can be taken to be a function with two arguments, the first being a function $\phi$ from
$\tau$ to {\tt prop}, the second being a variable $x$ of sort ({\tt in}  $\tau$), and the output being a proof $xx$ of sort ({\tt that} $\phi(x)$).   Here we see an argument which is an abstraction rather than an entity.   Finally, the universal quantifier over general types can be implemented as a construction taking three arguments, a type $\tau$ of sort {\tt type},
a predicate $\phi$  which is a function from $\tau$ to {\tt prop}, and a variable $x$ of sort ({\tt in} $\tau$), and sending this to output a proof $xx$ of sort ({\tt that} $\phi(x)$).
In this example we see that the sort of an input to an abstraction may depend on the sorts of earlier inputs.

The process of definition allows us to introduce new abstractions (and in special cases entities) whose existence follows from our axioms.  This is done in an entirely familiar way:
introduce a sequence of variables of appropriate types as under the previous heading.   Write out an expression  for an entity in terms of these variables which is well-sorted, using constructions
already given, and compute its sort.  We thereby discover a new abstraction which from any sequence of objects with the appropriate input sorts generates an object of the appropriate output sort (of course an entity sort, and possibly depending on the values of the inputs).

The next ingredient of the philosophical approach we take is a special attitude toward functions.   We do not choose to regard functions as infinite tables determined by applying the function to every possible sequence of inputs of the appropriate sorts.  We take the approach that a function is an actual rule or construction (necessarily speaking informally).  Where $x$ is a natural number variable, we take seriously the idea that the variable expression $x^2+1$ determines a function.  We really have ``variables" in a suitable sense in our metaphysics, as it were.  We are avoiding actual infinite totalities in our metaphysics, so we hardly want a function to be an actual infinite table of the outputs corresponding to given inputs:  we prefer to think of it as a (presumably finite) template into which inputs can be plugged to produce outputs.  The function $(\lambda x.x^2+1)$ is more like an abstraction from the term $x^2+1$ than the infinite table $\{(1,1),(2,4),(3,9),\ldots\}$.  Our attitude that functions are in a sense finite objects is borne out by the fact that abstractions are never introduced except as abstracted from explicit terms with variables in them (which are finite structures):  which is not to say that we identify functions with pieces of text.  We are talking about metaphor here.   The fact that functions are abstracted from expressions, which are certainly finite objects, appears to make it reasonable to suppose that a function is a finite object.  The expression tells us how to compute the value of the function for any input values of appropriate types which are presented to us:  we do not need to know all the values in advance:  the infinity of values of the function is potential, not actual.
 To sumarize, we take the position that a function is a finite object abstracted from a finite expression,
with slots in it (corresponding to variables in the expression) into which objects of the appropriate sorts can be plugged to give an output of the appropriate sort.  The same view applies to the axiomatically given constructions:  it is just that in this case we cannot look under the hood and see how the output is computed (and we do {\em not\/} support any presumption that arbitrary functions introduced in the future will be definable).  We now propose to take seriously the notion of a variable object.  Our suggestion
is that a variable object is an object in a possible world to which we have limited access:  all we know about the object is its sort (which may contain information about previously postulated objects).  One metaphorical way to think of a variable is that this is an object of a given sort which may be given to us in the future.  We use this temporal metaphor explicitly in the latest version of Lestrade by calling the possible worlds ``moves".\footnote{This nomenclature is a recent decision.  In earlier material, what we call here ``moves"  were called ``worlds",
and what we call the ``next move" was called the ``current world" and what we call the ``last move" was called the ``parent world".}

An initial version of this scheme which can be used to present the idea concretely is that we have a sequence of moves at any given point in the process indexed by concrete finite natural numbers 0 to $i+1$ (there are always at least two moves).  Move $i$, which we call the last move, and all lower indexed moves, are inhabited by entities and abstractions which we currently regard as given constant things of their various sorts.  Move $i+1$ we view as inhabited by variable objects, or in terms of our metaphor, things of determinate sort which have not been fixed yet:  we call it the next move.

We can freely declare variables in move $i+1$ of any entity sort we currently have accessible (remember that if we postulate new entities of sort {\tt prop} or {\tt type}, we have thereby postulated new sorts).  The declaration of an object in move $i+1$ may make more sorts available which depend on this object:  these may be used as sorts of further
entities declared in move $i+1$.  We regard all declarations as having been made in an order, with any variable having been declared later than any variable on which its sort depends.
One should remember that ``variables" are always objects at the next move.

We can introduce new primitive constructions:  when a sequence of variables has been postulated at move $i+1$ (some of whose sorts may depend on variables declared earlier)
we may postulate a construction (a new abstraction at move $i$, the last move) which will take any inputs of the given sorts and give an output of a stated entity sort
(which may depend on the input variables).  We may also declare a constant entity in move $i$ of a given fixed entity sort (this can be viewed as a special case where the argument list is empty; we are here introducing an entity constant rather than an entity variable).   The typical format of a construction command is ``construct $f$ which takes arguments $x_1,\ldots,x_n$ (variables declared in move $i+1$ in the order in which they were declared (which ensures that the sort of an $x_i$ can depend on an $x_j$ only if $j<i$) and is such that $f(x_1,\ldots,x_n)$ is of sort $\tau$ (an entity sort which may depend on the $x_i$'s).  Any variable on which the sort of an $x_i$ depends must appear as an $x_j$.

We can introduce new defined abstractions:  any expression using previously given constructions may be taken to determine a new abstraction at the last move (move $i$) taking as arguments a list of variables including all variables actually appearing in the expression (and all variables on which the sorts of variables in the argument list depend, in such a way that sorts of later arguments depend only on earlier arguments).  We can introduce entities by definition as well.  When we introduce an object by definition, we assign it a name:  so we define a new abstraction or entity symbol.  The typical form of a definition is ``define $f(x_1,\ldots,x_n)$ as $t$", where $x_1,\ldots,x_n$ is an argument list of variables just as above
and $t$ is an entity term.  This declares the symbol $f$ as an item at move $i$ whose input sorts are of course the sorts of the $x_i$'s, and whose output sort is the sort of $t$
(which obvously may depend on the $x_i$'s as we expect $t$ itself will so depend).
  
Finally, we observe that the system of moves is dynamic.  At any point, we may fix everything we have postulated so far and open a new move $i+2$ which becomes the next move, whereupon move $i+1$ becomes the last move and we increment the parameter $i$.  Or we can close move $i+1$ (if $i \neq 0$) and return to considering objects
at move $i$ as variable relative to objects at moves of lower index (thus in effect decrementing the parameter $i$).

A subtle point is that it may not be clear that we have given ourselves the ability to declare variables of abstraction sorts.  But we have.  The procedure for doing this is
to open a new move, declare variables of appropriate types desired as input sorts and declare a primitive construction taking this argument list to the of the desired output  sort, then close that move.  The primitive construction declared remains as an object at the next move, a variable abstraction of the desired type.   

We have made the perhaps artificial decision that the output of an abstraction is never an abstraction:  we are in effect reversing the popular operation of ``currying" as far as possible.
This has a striking effect on the language of the computer implementation, at least in the current version.  There is no need for the user ever to write a representation of
an abstraction other than the atomic name under which it was constructed or defined.  [This is not to say that the ability to do this may not be desirable; we think it should be added:  but it is interesting that it is eliminable].  The user's input language contains no variable binding constructions, in fact.   A function is always defined as in calculus when we say $f(x)=x^2+1$, the name $f$ then being provided for the function; it is never introduced as $(\lambda x.x^2+1)$ (anonymously, as it were) though such notations may be generated in sort signatures by
Lestrade when the name under which a function was introduced passes out of scope due to the move at which it was created being closed.

\noindent {\bf 10/10 change allowing more abstraction arguments:} The code of 10/10 while not implementing $\lambda$ terms per se permits easier formation of abstraction arguments.
A term $f(t_1,\ldots,t_m)$ where $m$ is less than the arity $n$ of $f$ will be read as $[(x_{m+1}\ldots x_n) \Rightarrow f(t_1,\ldots,t_m,x_{m+1},\dots,x_n)]$ if the types of the $t_i$'s are appropriate.  This only works if the argument list is explicitly closed with parentheses.  If $f$ is a defined abstraction, it will be definitionally expanded.   Such a term can only appear as an argument.

There is a comment to be made here about the relation between Lestrade and Automath:  Automath has a sublanguage PAL in which $\lambda$-terms are not used (\cite{automathbook}, pp. 79-88), but PAL is strictly weaker than Automath in its logical capabilities:  Lestrade has a more complex context mechanism which allows the logical power added to Automath by adding $\lambda$-terms to be replicated, but without actually having $\lambda$-terms; we expect to add the ability to enter $\lambda$-terms and perhaps other variable binding constructions at some point, but this will add convenience to Lestrade, not logical power.

The system is forced to have an internal representation of anonymous abstractions, however.   Whenever it computes the sort of an abstraction, it produces a dependent type in which the
arguments are bound variables:  a typical sort is $[(x_1:\tau_1),\ldots,(x_n:\tau_n) \Rightarrow \tau]$:  this is an expression with bound variables because the types $\tau_i$ may depend on $x_j$ for $j<i$ and $\tau$ may depend on any $x_i$.  The sort data recorded for a defined abstraction has the body $t$ of the definition as a further piece of information  (the form is $[(x_1:\tau_1),\ldots,(x_n:\tau_n) \Rightarrow (t,\tau)]$, where $\tau$ is the type of $t$):  this ``sort"  is in effect a $\lambda$-term.
Now consider what happens when the last move contains a defined abstraction and the next move is closed.  Primitive constructions declared in the former last move (which is now the next move) become variable abstractions.   Defined abstractions declared in the former last move remain defined, but they are now as it were defined variable
expressions.   When a construction or definition is declared in a way which uses one of these defined abstractions, the name of the defined abstraction must be eliminated from the information recorded at the last move, because if we close the current next move the name of the defined abstraction ceases to be available.  Where the defined abstraction appears in applied position,
it is eliminated by expanding the definition in the obvious way.   Where it appears as an argument, it can be replaced by its sort, which is in effect a lambda term defining it, in which the only terms denoting objects in the same move are bound variables representing its arguments.

One must note that the implementation of this scheme is inevitably recognizably a variant of the Automath type checker.   We are quite happy to acknowledge this.  It differs from the Automath type checker in not having lambda abstraction and function application as operations directly available to the user (in the present version:  for practical reasons we do intend to extend the language with user-defined lambda abstractions\footnote{The recently added implicit argument inference feature seems to make it much less urgent to introduce user-entered $\lambda$ terms.  The still more recently added feature allowing curried abstraction arguments may make this still less urgent} ).

If this is taken to be vague, it is made concrete by the actual implementation.  We support our claim that this scheme implements all mathematical activity by actually presenting declarations implementing standard foundational systems.

A philosophical point to be made about this scheme is that all infinities involved in it are potential.   One is never directly considering an actual infinity of objects at any point.   One manipulates entire sorts not by considering all the objects individually (requiring that the whole sort be given at once) but by the device of generality.  One shows not how to deal with {\em every \/} object of a potentially infinite sort, but how to deal with {\em any\/} object of this sort which may be  presented to one.  We believe that this defuses any objections to impredicative reasoning.  From this standpoint it is clear that {\em definition} is not a logically trivial move:  to define an abstraction (or entity) is to explicitly  make some potential object actual, and since we do not suppose that all potential objects are given to us at once, this is a nontrivial act.

It might be thought that this approach is more suited to constructive reasoning, and indeed it can implement constructive reasoning.  But classical reasoning is also amenable to implementation in this style.

We narrate the introduction of the general universal quantifier operation described in our earlier example.

Initially, we have move 0 and move 1, each containing no declared objects.

Declare a variable $\tau$ of sort {\tt type}  in move 1.

Open move 2:  the last move is now move 1 and move 2 is the next move.

Declare a variable $x_1$ of sort ({\tt in} $\tau$) at move 2.

Construct an abstraction $\phi$ such that $\phi(x_1)$ is of type {\tt prop}.  $\phi$ is introduced as an item
at move 1.

Close move 2.  $\phi$ is now a variable abstraction taking a type $\tau$ argument and returning a proposition (an arbitrary predicate of
type $\tau$ objects).

Declare a variable $x$ of sort ({\tt in} $\tau$).

Construct an abstraction $\forall$ such that $\forall(\tau,\phi,x)$ is a proposition.  $\forall$ is declared at move 0 (an unequivocally constant object).  This is the correct type for
the universal quantifier.  You can see in the extended example in section 6 how we further construct abstractions implementing
the rules for reasoning with the universal quantifier.

\newpage

\section{Formal description of the syntax and sort checking}

There is lots of sample Lestrade input and output to examine in later sections (and some embedded in this section).

\subsection{Lestrade Syntax}

\begin{description}

\item[books made up of lines:]  A Lestrade book (a Lestrade text is called a book, following Automath usage) is a sequence of lines.

\item[resolution of lines into tokens:]  Each line breaks up into tokens.

A token is either

\begin{enumerate}

\item a non-null string which is the concatenation of a single upper case letter or null string followed by a possibly null string of lowercase letters followed by a possibly null string of digits.
One of the three components has to be non-null.

\item a string of special characters taken from

\begin{verbatim}
~!@#$%^&*-+=/<>|?
\end{verbatim}

\item or a single character taken from 

\begin{verbatim}
,:()
\end{verbatim}

\item If a double quote appears in a Lestrade line, the entire remainder of the line is taken to be a token (without the quote).  This is used by some diagnostic commands not listed here:  this item is mainly a warning not to use double quotes in Lestrade lines.\footnote{{\tt parsetest "<text>} will parse the text as a term and display the term and its type.  {\tt parsetest2 "<text>} will parse the text as an entity type.  I have more diagnostic functions which I will probably eventually implement in the Lestrade interface using this device.}

\end{enumerate}

Whitepace is ignored, except that it terminates a token.  {\tt A 1} is two tokens, where {\tt A1} is one.

\item[command names:]  The first token in a line will be a command name, one of

\begin{verbatim}

declare, construct, define, rewritec, rewrited, open, close, clearcurrent  
  [the first line contains all the commands 
  that implement the logical framework]

save, foropen, forclearcurrent

[commands to do with context saving]

showall, showrecent, showdec, showdecs, displayrewrites
  [commands that will cause display of information]

showimplicit, hideimplicit

[turns on or off the display of implicit arguments in sorts]

readfile

[read a first file, output to a second.  The filenames must be
 legal Lestrade identifiers.

load

[load a saved theory]

comment or %, comment1 or %% (logged comments);
>>  for transient comments 

  [comments, logged and  unlogged, respectively; 
  comment or % is followed by a line break]

\end{verbatim}

\item[reserved tokens, identifiers:]  The reserved tokens are the comma, colon and parentheses, and the words {\tt obj}, {\tt prop}, {\tt type}, {\tt that}, and {\tt in}.   All other
tokens are identifiers (there is no reason that the command names cannot be identifiers).  An identifier which appears as the next token after the command {\tt declare}, {\tt construct}, or {\tt define} should not have been previously
used and will be assigned a sort if the line is well-formed and can successfully be executed.   Any other identifier which appears in a line should have been declared already
and assigned a sort.

A piece of information we need is that the sort of an identifier will tell us whether it represents an entity or an abstraction, and if it represents an abstraction the sort will tell us the arity of the abstraction, a positive integer.

\item[logical commands:]  The only lines which contain parsed text after the initial command are {\tt declare}, {\tt construct}, and {\tt define}, and now the new commands
{\tt rewritec} and {\tt rewrited}.

\item[comment commands:]  {\tt comment} or {\tt >>} will be followed by comment
text which will be ignored.   {\tt comment1} is for non-final lines in comments (not followed by a line break when echoed);
{\tt \%} and {\tt \%\%} are alternative versions of {\tt comment}, {\tt comment1}.

\item[display commands:]  The {\tt showdec} command may be followed by a single identifier whose declaration will be shown.  Text following other commands will be ignored.

\item[the readfile command:]  {\tt readfile file1 file2}  will read the Lestrade commands in {\tt file1.lti} and log to {\tt file2.lti}.  The filenames need to be valid Lestrade tokens.

\item[the load command:]  The {\tt load} command will be followed by a token (the name of a Lestrade log file which has been read by the system, it is hoped).\footnote{The double quote token construction can handle file names which are not well formed Lestrade tokens otherwise.}

\item[toggling display of implicit arguments in sorts:]  {\tt showimplicit} turns on display of implicit arguments in sorts; {\tt hideimplicit} turns it off again.

\item[syntactical forms of the logical commands:]  The form of a {\tt declare} line is the keyword {\tt declare}, followed by the undeclared identifier to be declared (which cannot be a reserved token), 
followed by an entity sort term.   These classes of strings will be explained.  Note that there is {\bf not} a colon before the entity sort term as in the following command.

The form of a {\tt construct} line is the keyword {\tt construct}, followed by the undeclared identifier to be declared (which cannot be a reserved token) followed by an argument list,
followed (optionally) by a colon\footnote{The colon in this and following commands is {\bf not} there to set off a following sort term but to terminate an argument list whose length cannot be predicted:  in the construct command its use is optional;  in the define command, where what follows is {\bf not} a sort, its use is mandatory.}, followed by an object sort term.

The form of a {\tt define} line is the keyword {\tt define}, followed by the undeclared identifier to be declared (which cannot be a reserved token), followed by an argument list,
followed by a colon (required), followed by an entity term.

The form of a {\tt rewritec} or {\tt rewrited} command is the respective keyword, followed by the identifier to be declared ({\tt rewritec}) or used as witnessing proof [in this case already declared] ({\tt rewrited})
followed by an argument list, followed by a colon (required), followed by a single identifier.

\item[entity sort terms:]  An entity sort term is either {\tt obj}, or {\tt prop}, or {\tt type}, or the keyword {\tt that} followed by an entity term, or the keyword {\tt in} followed by an entity term.

\item[entity terms:]  An entity term is either an identifier declared to be of entity type, or an identifier declared to be an abstraction of arity $n$, followed by an argument list of length $n$,
or an entity term followed by an identifier declared as an abstraction of entity $n \geq 2$  followed by an argument list of length $n-1$ which cannot be enclosed in parentheses unless it is of length 1
(this is an infix or mixfix term).  An entity term may optionally be enclosed in parentheses.  The precedence order assumed is that of the old computer language APL (every infix or mixfix is of the same precedence, except unary prefix operators, which bind more tightly, and everything groups to the right).

The display functions will always use infix form when an abstraction is of arity 2 and has first argument an entity.  Mixfix forms other than infix are never displayed.   The display functions use lots of parentheses and commas.

\item[argument lists:]  An argument list may be of length 0, in which case it is the null string (null argument lists occur only as the third item in a {\tt construct} or {\tt define} line:  the parser never finds them).  An argument list of positive length $n$ may optionally be enclosed in parentheses unless it is of length greater than one and follows a mixfix operator.  The $n$  items in it are either entity terms or lone identifiers of abstraction type;  individual items may optionally be separated by commas, which may be necessary to avoid an abstraction item from being read as a prefix, infix, or mixfix operator, depending on its arity.  

\noindent {\bf 10/10 change allowing more abstraction arguments:} The code of 10/10 while not implementing $\lambda$ terms per se permits easier formation of abstraction arguments.
A term $f(t_1,\ldots,t_m)$ where $m$ is less than the arity $n$ of $f$ will be read as $[(x_{m+1}\ldots x_n) \Rightarrow f(t_1,\ldots,t_m,x_{m+1},\dots,x_n)]$ if the types of the $t_i$'s are appropriate.  This only works if the argument list is explicitly closed in parentheses.  If $f$ is a defined abstraction, it will be definitionally expanded.   Such a term can only appear as an argument.

{\bf serious warning about the parser:} It is important to note that if the first item in an argument list following a prefix operator is enclosed in parentheses, one must also enclose the entire argument list in parentheses, to avoid reading the first item as the entire argument list.

\end{description}

\subsection{Lestrade Sort Declaration and Checking}

All that we have revealed so far is the syntactical requirements for a line.   There are semantic requirements as well, of course, handled by the sort checking functions of the prover.

\begin{description}

\item[the scheme of moves:]  At any given point the user has a finite list of sort declarations of identifiers which appear in the order in which they were added to moves 0 to $i+1$, where $i$ is a parameter maintained by the program.   Any identifier declared at any move is accessible to the parser and sort-checker.  We call move $i+1$ the next move and we call move $i$ the last move.  Some moves may be assigned names other than the default numerical name derived from their distance from move 0:  moves which do not have such names are assigned their default numerical names.

\item[opening a new move:]  The {\tt open} command causes $i$ to be incremented and a new move $i+1$ to be opened with no definitions in it.  [The command
{\tt open} with an argument will open an existing version of move $i+1$ named by that argument or create a new one.]  The {\tt open} command will fail if one tries to open
a move which does not have the default numerical name from a move which does have the default numerical name.

\item[saving the next move:]  The {\tt save} command will save the current state of moves 1 to $i$ with their current attached names and the next move
with the name supplied to it as an argument, or its current name if no argument is supplied, except that it will refuse to save move $i+1$ with the default numeral name $i+1$.
Moves can no longer be saved with their default numeral names at any level of the tree, nor will it save the next move with a non-default name if the last move has the default name, unless the last move is move 0.  The option of saving the next move with a new name allows Lestrade to fully emulate
the Automath context system (though more verbosely).  The next move will be renamed to the argument of the save command if the save is successful.

\item[closing the next move:]  The {\tt close} command does nothing if $i$ is 0.  If $i$ is positive, the command discards move $i+1$ and all declarations contained in it and decrements the counter $i$.  The close command does not save any declaration information:  saving of a move must be done explicitly.

\item[the clearcurrent command:]  The {\tt clearcurrent} command has the effect of {\tt close} followed by {\tt open}:  it empties the next move of declarations while not changing the index of the next move.  {\tt clearcurrent} is needed as a separate command, however, because move 1 cannot be closed, but can be cleared of declarations.  [{\tt clearcurrent} with an argument will empty the current version of move $i+1$ and replace it with a new one or previously saved one named by the argument].  The command will not allow you to give a move a name other than its default numerical name if its last move has the default numerical name, except in the case where the last move is move 0.  The first named move that is opened in a given session will be opened with {\tt clearcurrent} (or have its name changed to a non-default name by {\tt save}).

\item[discovering saved moves:]   The {\tt foropen} command will show what named moves can currently be opened with {\tt open}
and the {\tt forclearcurrent} command will show what named moves can currently be opened with {\tt clearcurrent}.  The display will just be a list of names.

\item[the display commands:]  {\tt showall} will show all declarations.  {\tt showrecent} will show all declarations at the last and next moves.  {\tt showdecs}
will show the declarations in the last and next moves, one at a time, those in the next move in order of declaration and the others in reverse order.  {\tt q}
will break out of either of the lists in {\tt showdecs}.  {\tt showdec} takes an identifier argument and displays its sort.

\item[the load command:]   This command with argument {\tt <filename1>} (which must be a token) will clear everything and load move 0 as it was when the command
{\tt readline} {\tt <filename1> <filename2>} was last run (along with some internal serial numbers).  Don't put the file name in quotes in the Lestrade interface!  No context information will be saved, just exactly the information in move 0; so this is a rather limited include feature.  Nor is there any way to merge theories.

\item[variables:]  We refer to all identifiers declared at the next move without definitions as {\em variables\/}.  Abstraction variables will always have been introduced with the {\tt construct} command at some stage when the current next move was the last move.

\item[argument list semantic conditions:]  We provide further that all items in the argument lists which follow the keyword and the identifier being declared in {\tt declare}, {\tt construct}, and {\tt define} commands
must be variables.  Moreover, they must appear in the order in which they were declared.  This is an easy way to enforce the constraint that the sort of a variable in such an
argument list may depend on a second variable appearing in the argument list only if this second variable appears earlier in the list.   Conversely, if the sort of a variable in the list depends
on any second variable, this second variable must appear earlier in the argument list.  [This is enforced by temporarily replacing the next move with just the items in the argument list and then declaration/sort checking all items in the argument list:  it is interesting to note that the internal data structure used to represent an argument list in Lestrade is exactly the same data structure used to represent a move.]

\item[semantics of the declare command:]  The command \begin{center}{\tt declare <identifier> <entity sort>}\end{center}  declares the identifier as having the given entity sort (as long as the identifier is undeclared and the entity sort term sort-checks).
The entity sort terms {\tt obj}, {\tt prop}, {\tt type} of course will sort check.  The term {\tt that <entity term>} sort-checks iff the entity term has sort {\tt prop}.   The term {\tt in <entity term>} sort-checks iff the entity term has sort {\tt type}.   

This declaration is added to the next move if it succeeds.

\item[constructing an entity:]  The command \begin{center}{\tt construct <identifier> : <entity sort>}\end{center} (with null argument list) works as the previous line does, except that the identifier is declared (and type checked) at the last move.

\item[general semantics of the construct command:]  The command \begin{center}{\tt construct <identifier> <argument list> : <entity sort>}\end{center} declares an abstraction of arity equal to the (positive) length of the argument list.  If the argument list
is $(x_1,\ldots,x_n)$, with the type of $x_i$ being $\tau_i$,  and the entity sort is $\tau$,  the type recorded is $$((x_1,\tau_1),\ldots,(x_n,\tau_n) \Rightarrow ({\tt ---},\tau),$$ where
the variables $x_i$ are to be understood as bound.   This may be a quite complex dependent type:  observe that each $\tau_i$ may contain occurrences of $x_j$ for $j<i$,
and $\tau$ may contain occurrences of any of the $x_i$'s (and any variable occurring free in $\tau$ or any $\tau_i$ must be one of the $x_i$'s).  Lestrade in fact renames all the bound variables, attaching a fresh numerical tag (the same tag for all variables in this argument list) to each of the $x_i$'s, and this procedure is repeated whenever a substitution is made into an abstraction sort, to avoid bound variable collision problems.
Note that some or all of the $\tau_i$'s may themselves be abstraction sorts, but the output type $\tau$ must be an entity sort.  The symbol {\tt ---} is a marker indicating that this is a primitive construction rather than a defined construction.  This declaration is added to the last move.

\item[defining an entity:]  The command \begin{center}{\tt define <identifier> : <entity term>},\end{center} where the entity term (which we write as $D$)  type checks, is declared with the sort $(D,\tau)$ [as if it were a nullary abstraction],
where $\tau$ is the type of $D$.   The identifier can be expanded to $D$ by prover functions as required. 

\item[general semantics of the define command:]  The command \begin{center}{\tt define <identifier> <argument list> : <entity term>} ,\end{center} where the entity term  (which we write  as $D$) checks with type $\tau$,
will be assigned the sort $((x_1,\tau_1),\ldots,(x_n,\tau_n) \Rightarrow (D,\tau)$:  this will succeed exactly if $\tau$ is the sort of $D$ and $((x_1,\tau_1),\ldots,(x_n,\tau_n) \Rightarrow ({\tt ---},\tau)$ is a well-formed abstraction sort (satisfying expected restrictions on variable dependencies, etc.)  This serves as the recorded type of the identifier
declared (and the declaration is created at the last move).  The identifier can be expanded to its sort (understood as a $\lambda$-term) by prover functions,
and any term $D(t_1,\ldots,t_n)$ which is well-typed can be expanded by substitution in the natural (and rather complex) way.

\item[semantics of the {\tt rewritec} and {\tt rewrited} commands:]  These commands allow construction of a primitive function witnessing validity of a rewrite rule ({\tt rewritec}), or introduction of a rewrite rule by exhibiting a function already constructed or defined of the appropriate type ({\tt rewrited}).  The argument list component consists of an initial segment of the actual argument list of variable arguments for the initial identifier, followed by two complex terms, the pattern {\tt pattern} and the target {\tt target}.  The pattern and the target must be of the same type $\tau$,
and the last variable $P$ in the argument list must be a variable typed as a function from $\tau$ to {\tt prop}.  The identifier after the colon (which must be undeclared)  is declared with type
{\bf that} $P({\tt pattern})$.   The intended type of the initial identifier is then the type sending the argument list before the pattern, with the new variable appended, to the output type
{\bf that} $P({\tt target})$.  All the variables in the argument list must actually appear in {\tt pattern}; $P$ cannot appear in {\tt pattern}; all variables that appear in {\tt target} must appear in {\tt pattern}.  The existence of a function of the indicated type will witness that any proof of $P({\tt pattern})$ can be mapped to a proof of $P({\tt target})$, which certainly looks like evidence that {\tt pattern} must be equal to {\tt target} for all values of the included variables.  The initial identifier in the case of {\tt rewritec} must be undeclared, and is declared with the indicated type if all the conditions hold (and a rewrite rule ${\tt pattern} := {\tt target}$ is recorded, with variables moved to a fresh namespace).  The initial identifier
in the case of {\tt rewrited} must already have been defined or constructed with the indicated type -- if this checks out, a rewrite rule is recorded in the same way.  When a rewrite rule is added to the last move with {\tt rewrited}, any rewrite rule already associated with the same first identifier in the last move  is deleted (but such rules are not deleted from lower indexed moves):  this supports changes in the order in which rewrite rules are applied, as roughly speaking the most recently introduced rewrite rule is attempted first.

\item[expansion of local definitions:]  Whenever a declaration is created in the last move, all defined identifiers in it whose declarations are found in the next move must be expanded (this is one
way that the internal $\lambda$-terms appear; they can also be introduced by the implicit argument inference mechanism):  this must be done because a declaration at the last move needs to continue to make sense if the next move is closed.
Of course, if such expansions reveal dependencies on variables not found in the argument list, an error is reported.

\item[computation of sorts of terms:]  We now discuss assignment of sorts to terms, and computation of identity of sort terms [and, as it turns out, computation of definitional expansions].

A bare identifier (whether an entity or an abstraction) is assigned type by lookup, with the proviso that an identifier typed $(D,\tau)$ (denoting a defined entity) is
simply assigned type $\tau$.   A term $f(t_1,\ldots,t_n)$ (or $t_1 f t_2,\ldots,t_n$) is assigned a sort using a procedure of matching of the declared type
$(x_1,\sigma_1),\ldots,(x_n,\sigma_n) \Rightarrow (D,\sigma)$ of the identifier $f$ (where $D$ can be {\tt ---} if $f$ is constructed or a term if $f$ is defined) with the list of types $\tau_i$ of the $t_i$'s.
If the type $\tau_1$ is identical to the type $\sigma_1$, we continue by replacing $x_1$ with $t_1$ in each $\sigma_i$ and $\sigma$ to obtain $\sigma_i^*$ and $\sigma^*$,
then continuing the matching of $(t_2,\ldots,t_n)$ with $(x_2,\sigma^*_2),\ldots,(x_n,\sigma^*_n)\Rightarrow (D^*,\sigma^*)$.  The lengths of the two argument lists must be the same
for the matching to succeed, and the final type assigned is the final form of $\sigma^*$.  It is worth noting that if $f$ is defined, the term $D^*$ obtained at the end of this procedure is the definitional expansion of $f(t_1,\ldots,t_n)$.

\item[considerations of bound variable naming and definitions]  may allow typographically distinct terms and sorts to be identified:  Abstraction sorts are regarded as identical when one can be converted to the other by renaming of bound variables:  the computational procedure for checking this is similar to the term typing procedure.   Terms are regarded as identical when an expansion of defined terms can convert one to the other:  when a failure to match is encountered, an attempt is made to correct it by expanding definitions [or by rewriting].   Due to the simple form of our terms, this is straightforward to compute.

\item[expansion of definitions:]  Expansion of $F(t_1,\ldots,t_n)$ where the type of $F$ is $((x_1,\tau_1),\ldots,(x_n,\tau_n) \Rightarrow (D,\tau)$ is a straightforward matter of substitution [already described incidentally in the discussion of type matching], and sort safe as long as all terms involved have already been sort-checked.  Have $t_1$ replace $x_1$ in all terms and types, then $t_2$ replace $x_2$, and so forth.
As noted above, a defined abstraction appearing by itself as an argument expands to its sort, which can be understood as a $\lambda$-term.

\item[action of rewrite rules:]  This is a new feature.  When a term is rewritten, the most recently introduced rewrite rule whose pattern matches it is applied to it [the precise order is, most recently introduced rewrite rule in the highest indexed move which contains such a rewrite rule; the use of the most recently introduced rule first makes it easier for the user to reprogram the rewrite system].  Equations between entities and types will succeed if they can be made to work by rewriting, just as they will succeed if they can be made to work by definition expansion.  Rewrites introduced at the last move or moves of lower index are applied, though rewrites at the next move are kept.  The list of rewrites recorded at each move is managed in the same way as the moves are managed by the open, close, clearcurrent commands.

To maintain confluence, there is a subtlety about matching:  when a pattern is matched with a term, proper subterms of the term are rewritten first, and the match succeeds only if it still works after this rewriting.

The entity output of the define command is rewritten.  Types are not rewritten in any command (though rewriting may be used to verify equivalences between types:  the fact that the system does this may indicate that the rewrite system adds something to the logical framework).  

\item[another type inference feature (implicit arguments):]  In the latest update, Lestrade has an extensive ability to accept construction and definition declarations of abstractions with missing arguments.  Lestrade infers and supplies the required missing arguments at declaration time by finding variables in the types of the explicitly given arguments, and it is able to determine where to put them to get a well formed argument list.  When typing an instance of an abstraction with implicit arguments, it deduces the values of the implicit arguments by type matching.  This includes the ability to find abstraction arguments which have not been given explicit names.  For an abstraction implicitly given to be recovered when it appears in a type in applied position, it must have a fully abstract occurrence, that is, one in which its arguments coincide with an initial segment of the bound variables in its context in the correct order [the circumstances under which implicit arguments can be inferred are quite complicated to describe; there are extensive examples in sample files on the Lestrade page, mostly simple ones, but a better essay on this is needed].

Implicit arguments can be recognized in the output because their names have a prepended dot.    The implicit argument mechanism does not affect the core logic of Lestrade at all:
all terms in fact have all usually expected arguments.  It is purely a function of input and output.   Sort displays show all input arguments, implicit and explicit of an abstraction whose type is being displayed:  the default behavior is that only explicit arguments of abstractions appearing in applied position in sort displays are shown:  this behavior can be changed as desired with commands {\tt showimplicit} and {\tt hideimplicit}.  The type matching which supplies implicit arguments may fail if some definitional expansion of the type being matched actually eliminates the information about the implicit argument in the type:  there is no guarantee that the way the implicit argument is presented in the type in the original declaration is deducible semantically from the value of the type:  it is entirely a matter of syntax.   Rewriting cannot be used to assist implicit argument extractions, so far.  

The precise condition under which an implicit abstraction argument appearing in applied position in the sort of an explicit argument
can be matched is that it appears either applied to arguments which do not depend on variables bound in that sort term (in which case it will match an explicit abstraction term in applied position applied to matching terms) or applied to arguments which are an initial segment of the variables bound in the context of that sort term, appearing in the precise order in which they appear bound in the context (in which case it will match a lambda term constructed by Lestrade, if the term matching it does not depend on any other bound variables in the context). [This is obscure and needs supporting examples.]

This feature may mean that it is much less urgent to install user entered lamba terms:  laborious construction of abstractions by the user (notably predicates quantified over) is greatly reduced, as the values of these arguments can often be inferred by higher order matching as anonymous internal lambda terms.  It does also mean that one should look carefully at the actual type of a defined or constructed abstraction, including its implicit arguments:  all computations with an abstraction with implicit arguments will use its full type.

\noindent {\bf 10/10 change allowing more abstraction arguments:} The code of 10/10 while not implementing $\lambda$ terms per se permits easier formation of abstraction arguments.
A term $f(t_1,\ldots,t_m)$ where $m$ is less than the arity $n$ of $f$ will be read as $[(x_{m+1}\ldots x_n) \Rightarrow f(t_1,\ldots,t_m,x_{m+1},\dots,x_n)]$ if the types of the $t_i$'s are appropriate.  This only works if the argument list is explicitly closed with parentheses.  If $f$ is a defined abstraction, it will be definitionally expanded.   Such a term can only appear as an argument.

\end{description}

\newpage 

\subsubsection{Analogies between Lestrade lines and Automath lines; remarks on differences between the Lestrade and Automath type systems}

An Automath line consists of four components:  an indicator of context, an identifier being declared, a definition of that identifier, and a type.  The analogue of the identifier component in a Lestrade line is always evident:  it is the identifier following the command keyword.  We are here considering only {\tt declare}, {\tt construct}, and {\tt define} lines in
a Lestrade text:  other lines do not have analogues in Automath.

In Automath, the context mechanism is fairly simple, consisting of either 0 or a reference to a previous identifier.  One is then defining the new identifier as
indicated by the definition component with a stated type in a context in which the previous identifier and all identifiers one finds by backtracking to its context indicator and iterating are present.

In Lestrade, the context is more complex:  the intersection of the explicit context given by Automath with variables declared in the next move  is present as the argument list to the identifier introduced by a  {\tt construct} or {\tt define} command (all items being listed, with no chained expansion of the context as in Automath), but additional related roles are played by the system of moves.  A variable introduced in Lestrade using the {\tt declare}
command may be thought of as having all previous variables in its (unindicated) context, and certainly at least those variables which explicitly appear in its sort:  note that when a variable appears in an argument list, it must be preceded by all the variables mentioned in its sort.  Argument lists in Lestrade are more flexible than contexts in Automath:  Automath contexts are restricted to forming a tree under the initial segment relation.

The function of the special definition body PN used in Automath for a primitive notion is handled in Lestrade by using the {\tt construct} command rather than the {\tt define} command.  In the {\tt construct} command a type is explicitly given; in the {\tt define} command, Lestrade computes and displays the sort:  the user does not need to supply it.  The function of the definition body {\tt ---}  used for variables in Automath is handled in Lestrade by using the {\tt declare} command, for entities, and by using the {\tt construct} command then closing the next move, to obtain variable abstractions.

The extra trick which allows user-entered lambda terms to be avoided (so far) is that whenever an expression is defined using a given context, what is defined is not a variable
expression specific to that context, in Automath usable in a different context by supplying terms of appropriate types to replace the elements of the original context not shared with the current context (considered as an argument list), but an abstraction (a function) present at the last move, in effect a name for the function implementing the variable expression, which would have to be expressed in Automath as a lambda term.  Further, constructing a primitive function of a given sort at the last move, then closing the next move, gives one a variable
abstraction of that type in the formerly last, now next move.  By the use of these two devices, and the fact that abstractions may be supplied as arguments to other abstractions using their atomic names, one entirely avoids the need to write nonatomic terms of abstraction sorts.  The further feature that abstractions do not have abstraction output ensures that the user never needs to type an abstraction sort.

It does appear that for fluent use of the system it would be advisable to introduce an ability to enter explicit lambda terms as arguments to Lestrade terms.  But it is useful to see that no logical power will be added to our framework if we do this.  We have not so far been convinced that a user ever should have to write an abstraction sort.

In Lestrade, propositions are not themselves sorts, but entities $p$ of a sort {\tt prop}, correlated with sorts {\tt that} $p$ of proofs of $p$.\footnote{Later dialects of Automath have propositions as types; references I have looked at suggest that earlier dialects have a system analogous to the Lestrade system, with $p$ an object of type {\tt prop} and a separate type {\tt Proofs}($p$) inhabited by proofs of $p$; it may always have been the case that the same notation was used in Automath code for the object of type {\tt prop} and the associated type.}  Mathematical types are for us
similarly entities $\tau$ of a sort {\tt type}, correlated with sorts {\tt in} $\tau$ inhabited by specific objects of those types.   Operations on propositions or types can then
be postulated  naturally by declaring abstractions with arguments of type {\tt prop} or {\tt type} (not actual sorts as arguments).

Abstractions are  not entities and cannot be outputs of constructions.\footnote{This is somewhat modified by the 10/10 update allowing abstraction arguments to be constructed by currying.}   But they can be inputs to constructions, and one can create entities correlated with abstractions of a particular sort.  For example, we do not identity proofs of $p \rightarrow q$ with functions from {\tt that} $p$ to {\tt that} $q$ (which are not entities), but instead provide a construction {\tt Ifproof} which sends such a function to an entity of the sort {\tt that} $p \rightarrow q$.  One can see in the extensive Lestrade book given in section 6 that this does not obstruct the usual sorts of reasoning in a system of this kind.  The move system of Lestrade implements reasoning under hypotheses ending with the proof of an implication or reasoning about arbitrarily postulated objects ending in the proof of a universally quantified statement quite naturally.

It is worth noting that Lestrade does not have any specific notion of function application, any more than it has user-written lambda terms.  Lestrade abstractions are always applied to their complete argument lists\footnote{or to partial argument lists in which the missing implicit arguments can be extracted from the sorts of the arguments given explicitly.} in the format in which they were defined (or appear by themselves as arguments\footnote{The 10/10 update allows abstractions with shortened argument lists to represent ``curried" abstraction arguments.}).  The application of Lestrade abstractions to their argument lists corresponds
to the application of atomic defined terms in Automath to argument lists replacing items in the context in which they were originally defined (substitution rather than application), although the atomic Automath terms are as it were variable expressions and the Lestrade abstractions denote functions.

A real difference in strength between Automath and Lestrade can be seen in connection with the odd subtyping supported in later versions of Automath.  Automath regards
a type inhabited by functions from type $\tau$ to type {\tt prop} as a subtype of type {\tt prop}.  Further, it regards an element $p$ of {\tt prop} as being itself a type
(in our terms, {\tt that} $p$ is identified with $p$).   This gives quantification over type $\tau$ for free.  Any predicate of type $\tau$ is itself a proposition.  If a predicate $P$
is inhabited by an object $pp$, then this object is itself a function from elements $n$ of type $\tau$ to proofs of $P(n)$:  so such an inhabitant is a proof of $(\forall n \in \tau:P(n))$.
Further, $\tau$ itself may be a quite complex abstraction type:  we get quantification of all orders for free from the type system.  In Lestrade, it is actually possible to postulate quantification over all entity types and abstraction sorts uniformly, indirectly, in a way suggested at the end of our large example of Lestrade text, but it takes considerable work.  To  provide quantification over any type (or indeed to get any inhabitant of {\tt prop} at all) requires some declarations in Lestrade.

Let's explore why an inhabitant of $P$ is a function from $n \in \tau$ to type $P(n)$.  The underlying idea is that an expression $f(x)$ of type $\sigma(x)$  is an expression
$\lambda x.f(x)$ of type $(\lambda x.\sigma( x))$.  So if $f$ is of type $(\lambda x.\sigma(x))$, $f(x)$ is of type $\sigma(x)$.  So if $P$ is of type $\tau \rightarrow {\tt prop}$,
which is supposed to be an element of type {\tt prop} as well, then an inhabitant of $P$ is a function $f$ such that $f(x)$ is an inhabitant of $P(x)$.  But $P(x)$ is also of type
{\tt prop}, so $f(x)$ is a proof of $P(x)$ for each $x$.   This is very cute.  It also has all the advantages of theft over honest toil, as Russell said in some similar context.

For us, for any type $\tau$ we must  declare a constructor sending any predicate $P$ of type $\tau$ to an element $\forall P$ of {\tt prop}, then further postulate
a function sending any function which takes $x : \tau$ to an element of {\tt that} $P(x)$ to a proof of $\forall P$.  We can declare these uniformly over the entity types
fairly easily, and with a little trickery we can indirectly declare them over all abstraction sorts as well.  But we also have the option of declaring only the quantifiers we want,
and never enabling higher order quantification.  An Automath theory is of necessity a higher order theory.

In his paper \cite{aut} describing his implementation \cite{freek} of Automath, Wiedijk discusses the difference between the $\lambda$-types of dependently typed functions
found in Automath and the $\Pi$-types found in the currently popular type systems.  I am frankly not certain of the place of the Lestrade type system on this dimension.  The difficulties that arise in Automath related to this issue seem to relate to the possibility of the same term appearing as both an object of the system and a sort, and this is ruled out in Lestrade.
The rule that if $F$ types as $G$, then $Fa$ types as $Ga$, which seems to characterize the $\lambda$-approach, holds for Lestrade, but note that $a$ must be a complete argument list,
$Fa$ must be an entity (and so cannot be a function) and $Ga$ must be an entity sort (and the entity sorts have no intersection with the objects of the theory, though entities can be packaged in them via the {\tt that} and {\tt in} constructors).  In fact, where $F:G$ and $F(a):G(a)$ in Lestrade, the four objects mentioned are all of different metasorts:  $F$ is an abstraction,
$G$ is an abstraction sort, $F(a)$ is a entity and $G(a)$ is an entity sort.

\newpage

\subsubsection{A sample Lestrade book with rewriting}

We introduce a brief Lestrade book to illustrate capabilities of the rewrite system.

Watson, referred to in the comments,  is an earlier theorem proving project of mine (see \cite{watson}), an equational prover which incorporated a fairly elaborate scheme of programming using interlocking rewrite rules, described in \cite{holmesrewriting}.  We believe that most features of the Watson programming system are implementable in Lestrade rewrites.  The rule {\tt Assocs} in the second example, which implements regrouping of arbitrarily complex nested sums so that all grouping is to the right, parallels a basic Watson example.

Notice in both examples that the final examples of {\tt define} commands exhibit ``execution behavior".  It would be useful to give an example in which the basic properties
underlying the rewrite rules are proved then the rewrites introduced using the {\tt rewrited} command, rather than having the functions providing evidence for the validity of the rewrite rules introduced by fiat using the {\tt rewritec} command:  though one should also note that {\tt rewritec} is a perfectly respectable way to introduce equational axioms.

\begin{verbatim}


>> Inspector Lestrade says:  
>>   Welcome to the Lestrade Type Inspector,
>>    test matching version of 7/2/2016
>>    1 pm Boise time


declare x obj

>> x:  obj {move 1}

declare y obj

>> y:  obj {move 1}

construct pair x y obj

>> pair:  [(x_1:obj),(y_1:obj) => (---:obj)] {move 0}

construct p1 x obj

>> p1:  [(x_1:obj) => (---:obj)] {move 0}

construct p2 x obj

>> p2:  [(x_1:obj) => (---:obj)] {move 0}

open

     declare x1 obj

>>      x1:  obj {move 2}

     construct P x1 prop

>>      P:  [(x1_1:obj) => (---:prop)] {move 1}

     close

rewritec First x y P, p1 pair x y, x : u

>> u:  that P(p1((x pair y))) {move 1}


>> First:  [(x_1:obj),(y_1:obj),(P_1:[(x1_2:obj) => (---:prop)]),
>>   (u_1:that P_1(p1((x_1 pair y_1)))) => (---:that P_1(x_1))] 
>>   {move 0}

rewritec Second x y P, p2 pair x y, y: v

>> v:  that P(p2((x pair y))) {move 1}


>> Second:  [(x_1:obj),(y_1:obj),(P_1:[(x1_2:obj) => (---:prop)]),
>>   (v_1:that P_1(p2((x_1 pair y_1)))) => (---:that P_1(y_1))] 
>>   {move 0}

open

     declare x1 obj

>>      x1:  obj {move 2}

     declare y1 obj

>>      y1:  obj {move 2}

     define reverse x1 : pair (p2 x1, p1 x1)

>>      reverse:  [(x1_1:obj) => ((p2(x1_1) pair p1(x1_1)):obj)] 
>>        {move 1}

     define reversetest x1 y1 :  reverse (pair x1 y1)

>>      reversetest:  [(x1_1:obj),(y1_1:obj) => (reverse((x1_1 
>>        pair y1_1)):obj)] {move 1}

     close

% notice that Lestrade executes the pair reversal!

define testing x y:  reversetest x y

>> testing:  [(x_1:obj),(y_1:obj) => ((y_1 pair x_1):obj)] 
>>   {move 0}

clearcurrent

% associative law simplication

% I believe I have implemented almost the full power of the Watson
% rewrite rule programming scheme.  The interlock between matching and
% rewriting should make it possible to implement its control structures
% without extra primitives.

construct Nat type

>> Nat:  [(---:type)] {move 0}

declare m in Nat

>> m:  in Nat {move 1}

declare n in Nat

>> n:  in Nat {move 1}

declare p in Nat

>> p:  in Nat {move 1}

construct + m n in Nat

>> +:  [(m_1:in Nat),(n_1:in Nat) => (---:in Nat)] {move 
>>   0}

construct assoc m in Nat

>> assoc:  [(m_1:in Nat) => (---:in Nat)] {move 0}

construct assocs m in Nat

>> assocs:  [(m_1:in Nat) => (---:in Nat)] {move 0}

open

     declare m1 in Nat

>>      m1:  in Nat {move 2}

     construct Pn m1 prop

>>      Pn:  [(m1_1:in Nat) => (---:prop)] {move 1}

     close

rewritec Assocfails m Pn, assoc m, m:u

>> u:  that Pn(assoc(m)) {move 1}


>> Assocfails:  [(m_1:in Nat),(Pn_1:[(m1_2:in Nat) => 
>>   (---:prop)]),(u_1:that Pn_1(assoc(m_1))) => (---:that 
>>   Pn_1(m_1))] {move 0}

rewritec Assocsfails m Pn, assocs m, m:v

>> v:  that Pn(assocs(m)) {move 1}


>> Assocsfails:  [(m_1:in Nat),(Pn_1:[(m1_2:in Nat) => 
>>   (---:prop)]),(v_1:that Pn_1(assocs(m_1))) => (---:that 
>>   Pn_1(m_1))] {move 0}

rewritec Assocrule m n p Pn, (m + n) + p, m + (n + p):w

>> w:  that Pn(((m + n) + p)) {move 1}


>> Assocrule:  [(m_1:in Nat),(n_1:in Nat),(p_1:in Nat),
>>   (Pn_1:[(m1_2:in Nat) => (---:prop)]),(w_1:that Pn_1(((m_1 
>>   + n_1) + p_1))) => (---:that Pn_1((m_1 + (n_1 + p_1))))] 
>>   {move 0}

rewritec Assocsrule m n p Pn, (m + n) + p, assocs(assoc(m + (assocs (n+p)))):x

>> x:  that Pn(((m + n) + p)) {move 1}


>> Assocsrule:  [(m_1:in Nat),(n_1:in Nat),(p_1:in Nat),
>>   (Pn_1:[(m1_2:in Nat) => (---:prop)]),(x_1:that Pn_1(((m_1 
>>   + n_1) + p_1))) => (---:that Pn_1(assocs(assoc((m_1 
>>   + assocs((n_1 + p_1)))))))] {move 0}

declare q in Nat

>> q:  in Nat {move 1}

define test m n p q:(m+n)+(p+q)

>> test:  [(m_1:in Nat),(n_1:in Nat),(p_1:in Nat),(q_1:in 
>>   Nat) => ((m_1 + (n_1 + (p_1 + q_1))):in Nat)] {move 
>>   0}

declare r in Nat

>> r:  in Nat {move 1}

declare s in Nat

>> s:  in Nat {move 1}

define test2 m n p q r s:((m+n)+p)+((q+r)+s)

>> test2:  [(m_1:in Nat),(n_1:in Nat),(p_1:in Nat),(q_1:in 
>>   Nat),(r_1:in Nat),(s_1:in Nat) => ((m_1 + (n_1 + (p_1 
>>   + (q_1 + (r_1 + s_1))))):in Nat)] {move 0}


>> Inspector Lestrade says:  Done reading scratch to rewrites:
>>  type lines or type quit to exit interface

quit

\end{verbatim}

\newpage

\subsection {A sample Lestrade book with implicit arguments}

\begin{verbatim}


>> Inspector Lestrade says:  
>>   Welcome to the Lestrade Type Inspector,
>>    full version of 7/22/2016 (implicit arguments upgrade)
>>   1 pm Boise time


declare p prop

>> p: prop {move 1}


declare q prop

>> q: prop {move 1}


construct & p q prop

>> &: [(p_1:prop),(q_1:prop) => (---:prop)] {move 0}


declare pp that p

>> pp: that p {move 1}


declare qq that q

>> qq: that q {move 1}


construct Andintro pp qq:that p & q

>> Andintro: [(.p_1:prop),(pp_1:that .p_1),(.q_1:prop),
>>   (qq_1:that .q_1) => (---:that (.p_1 & .q_1))] {move 
>>   0}


declare rr2 that p&q

>> rr2: that (p & q) {move 1}


construct Andelim1 rr2:that p

>> Andelim1: [(.p_1:prop),(.q_1:prop),(rr2_1:that (.p_1 
>>   & .q_1)) => (---:that .p_1)] {move 0}


construct Andelim2 rr2:that q

>> Andelim2: [(.p_1:prop),(.q_1:prop),(rr2_1:that (.p_1 
>>   & .q_1)) => (---:that .q_1)] {move 0}


define Ptest pp:Andintro pp pp

>> Ptest: [(.p_1:prop),(pp_1:that .p_1) => (Andintro(.p_1,
>>   pp_1,.p_1,pp_1):that (.p_1 & .p_1))] {move 0}


construct -> p q prop

>> ->: [(p_1:prop),(q_1:prop) => (---:prop)] {move 0}


open

     declare pp2 that p

>>      pp2: that p {move 2}


     construct Ded pp2 that q

>>      Ded: [(pp2_1:that p) => (---:that q)] {move 1}


     close

construct Ifintro Ded:that p -> q

>> Ifintro: [(.p_1:prop),(.q_1:prop),(Ded_1:[(pp2_2:that 
>>   .p_1) => (---:that .q_1)]) => (---:that (.p_1 -> .q_1))] 
>>   {move 0}


open

     declare q2 that q

>>      q2: that q {move 2}


     define qid q2:q2

>>      qid: [(q2_1:that q) => (q2_1:that q)] {move 1}


     close

define Selfimp q: Ifintro qid

>> Selfimp: [(q_1:prop) => (Ifintro(q_1,q_1,[(q2_2:that 
>>   q_1) => (q2_2:that q_1)]):that (q_1 -> q_1))] {move 
>>   0}


declare rr that p-> q

>> rr: that (p -> q) {move 1}


construct Mp pp rr:that q

>> Mp: [(.p_1:prop),(pp_1:that .p_1),(.q_1:prop),(rr_1:
>>   that (.p_1 -> .q_1)) => (---:that .q_1)] {move 0}


open

     declare x obj

>>      x: obj {move 2}


     construct P x prop

>>      P: [(x_1:obj) => (---:prop)] {move 1}


     close

construct Forall P: prop

>> Forall: [(P_1:[(x_2:obj) => (---:prop)]) => (---:prop)] 
>>   {move 0}


declare U that Forall P

>> U: that Forall(P) {move 1}


declare y obj

>> y: obj {move 1}


construct Ug U y that P y

>> Ug: [(.P_1:[(x_2:obj) => (---:prop)]),(U_1:that Forall(.P_1)),
>>   (y_1:obj) => (---:that .P_1(y_1))] {move 0}


open

     declare z obj

>>      z: obj {move 2}


     construct ui z that P z

>>      ui: [(z_1:obj) => (---:that P(z_1))] {move 1}


     close

construct Ui P, ui:that Forall P

>> Ui: [(P_1:[(x_2:obj) => (---:prop)]),(ui_1:[(z_3:obj) 
>>   => (---:that P_1(z_3))]) => (---:that Forall(P_1))] 
>>   {move 0}


open

     declare w obj

>>      w: obj {move 2}


     open

          declare zz that P w

>>           zz: that P(w) {move 3}


          define zzid zz:zz

>>           zzid: [(zz_1:that P(w)) => (zz_1:that P(w))] 
>>             {move 2}


          close

     define Q w:P w -> P w

>>      Q: [(w_1:obj) => ((P(w_1) -> P(w_1)):prop)] {move 
>>        1}


     define zzz w :  Ifintro zzid

>>      zzz: [(w_1:obj) => (Ifintro(P(w_1),P(w_1),[(zz_2:
>>        that P(w_1)) => (zz_2:that P(w_1))]):that (P(w_1) 
>>        -> P(w_1)))] {move 1}


     close

define test P: Ui Q, zzz

>> test: [(P_1:[(x_2:obj) => (---:prop)]) => (Ui([(w_3:
>>   obj) => ((P_1(w_3) -> P_1(w_3)):prop)],[(w_4:obj) => 
>>   (Ifintro(P_1(w_4),P_1(w_4),[(zz_5:that P_1(w_4)) => 
>>   (zz_5:that P_1(w_4))]):that (P_1(w_4) -> P_1(w_4)))]):
>>   that Forall([(w_6:obj) => ((P_1(w_6) -> P_1(w_6)):prop)]))] 
>>   {move 0}


declare r prop

>> r: prop {move 1}


open

     declare outerhyp that (p->q) & (q->r)

>>      outerhyp: that ((p -> q) & (q -> r)) {move 2}


     define firstlink outerhyp :  Andelim1 outerhyp

>>      firstlink: [(outerhyp_1:that ((p -> q) & (q -> 
>>        r))) => (Andelim1((p -> q),(q -> r),outerhyp_1):
>>        that (p -> q))] {move 1}


     define secondlink outerhyp : Andelim2 outerhyp

>>      secondlink: [(outerhyp_1:that ((p -> q) & (q -> 
>>        r))) => (Andelim2((p -> q),(q -> r),outerhyp_1):
>>        that (q -> r))] {move 1}


     open

          declare innerhyp that p

>>           innerhyp: that p {move 3}


          define step1 innerhyp:  Mp innerhyp firstlink outerhyp

>>           step1: [(innerhyp_1:that p) => (Mp(p,innerhyp_1,
>>             q,firstlink(outerhyp)):that q)] {move 2}


          define step2 innerhyp:  Mp (step1 innerhyp,secondlink outerhyp)

>>           step2: [(innerhyp_1:that p) => (Mp(q,step1(innerhyp_1),
>>             r,secondlink(outerhyp)):that r)] {move 2}


          close

     define step3 outerhyp : Ifintro step2

>>      step3: [(outerhyp_1:that ((p -> q) & (q -> r))) 
>>        => (Ifintro(p,r,[(innerhyp_2:that p) => (Mp(q,
>>        Mp(p,innerhyp_2,q,firstlink(outerhyp_1)),r,secondlink(outerhyp_1)):
>>        that r)]):that (p -> r))] {move 1}


     close

define Transimp p q r:  Ifintro step3

>> Transimp: [(p_1:prop),(q_1:prop),(r_1:prop) => (Ifintro(((p_1 
>>   -> q_1) & (q_1 -> r_1)),(p_1 -> r_1),[(outerhyp_2:that 
>>   ((p_1 -> q_1) & (q_1 -> r_1))) => (Ifintro(p_1,r_1,
>>   [(innerhyp_3:that p_1) => (Mp(q_1,Mp(p_1,innerhyp_3,
>>   q_1,Andelim1((p_1 -> q_1),(q_1 -> r_1),outerhyp_2)),
>>   r_1,Andelim2((p_1 -> q_1),(q_1 -> r_1),outerhyp_2)):
>>   that r_1)]):that (p_1 -> r_1))]):that (((p_1 -> q_1) 
>>   & (q_1 -> r_1)) -> (p_1 -> r_1)))] {move 0}


open

     declare x obj

>>      x: obj {move 2}


     construct ev x that P x

>>      ev: [(x_1:obj) => (---:that P(x_1))] {move 1}


     close 

construct Ui2 ev:that Forall P

>> Ui2: [(.P_1:[(x_2:obj) => (---:prop)]),(ev_1:[(x_3:
>>   obj) => (---:that .P_1(x_3))]) => (---:that Forall(.P_1))] 
>>   {move 0}


open

     declare x17 obj

>>      x17: obj {move 2}


     open

          declare ev2 that P x17

>>           ev2: that P(x17) {move 3}


          define evid2 ev2:  ev2

>>           evid2: [(ev2_1:that P(x17)) => (ev2_1:that 
>>             P(x17))] {move 2}


          close

     define theimp x17: Ifintro evid2

>>      theimp: [(x17_1:obj) => (Ifintro(P(x17_1),P(x17_1),
>>        [(ev2_2:that P(x17_1)) => (ev2_2:that P(x17_1))]):
>>        that (P(x17_1) -> P(x17_1)))] {move 1}


     close

define testing P : Ui2 theimp

>> testing: [(P_1:[(x_2:obj) => (---:prop)]) => (Ui2([(x17_3:
>>   obj) => ((P_1(x17_3) -> P_1(x17_3)):prop)],[(x17_4:
>>   obj) => (Ifintro(P_1(x17_4),P_1(x17_4),[(ev2_5:that 
>>   P_1(x17_4)) => (ev2_5:that P_1(x17_4))]):that (P_1(x17_4) 
>>   -> P_1(x17_4)))]):that Forall([(x17_6:obj) => ((P_1(x17_6) 
>>   -> P_1(x17_6)):prop)]))] {move 0}



>> Inspector Lestrade says:  Done reading scratch to Dtest:
>>  type lines or type quit to exit interface

quit

\end{verbatim}

\newpage

\subsection{Formalization of the type system of Lestrade}

We use the notation $T[a/x]$  for substitution of a notation $a$ for a variable $x$ in a notation $T$.

Metasorts {\bf esort} (entity sort), {\bf asort} (abstraction sort), and {\bf arglist} $n$ (argument list sorts for argument lists of length $n$) for each positive $n$
are postulated.  These are not internal sorts of Lestrade at all.  The union of the metasorts {\bf esort} and {\bf asort} is the metasort {\bf sort}.

{\tt prop} and {\tt type} are of metasort {\bf esort}.

If $p$ is of sort {\tt prop},  ({\tt that} $p$) is of metasort {\bf esort}.

If $\tau$ is of sort {\tt type}, ({\tt in} $\tau$) is of metasort {\bf esort}.

All terms of metasort {\bf esort} are built in these ways.

If $x$ is a term of sort $\tau$, then $\tau$ is of metasort {\bf sort}.

A countable supply of variables of each sort is given.

We view a list $[t_1,\ldots,t_n]$ as a function in the sense of the metatheory with domain \{1,\ldots,n\} sending each $i$ to $t_i$:  thus
$[t]$ is different from $t$ and any notation $[t_1,\ldots,t_n]$ is handled.

If $x_i$ is a variable of sort $\tau_i$ for each $i \leq n$, all $x_i$'s are distinct and no $x_i$ occurs in $\tau_j$ for $j\leq i$, then
$[(x_1,\tau_1),\ldots,(x_n,\tau_n)]$ is a term of metasort {\bf arglist} $n$.  All terms of this metasort are built in this way.  Notice
that $\tau_i$'s may be of metasort {\bf esort} or {\bf asort}.

If $t$ is a term of sort $\tau$, then $[t]$ is an argument list of list sort $[(x,\tau)]$ (the list sort being of metasort {\bf arglist} 1) (for any variable $x$ of sort $\tau$, technically speaking not occurring in $\tau$, something which will not naturally happen).  For $n>1$,  $[t_1,\ldots t_n]$ is an argument list of list sort $$[(x_1,\tau_1),\ldots,(x_n,\tau_n)],$$ where no $x_i$ occurs in any $t_j$, iff  the list sort is of metasort {\bf arglist} $n$ and
$t_1$ is of sort $\tau_1$ and $[t_2,\ldots,t_n]$ is of list sort $[(x_2,\tau_2[t_1/x_1]),\ldots,(x_n,\tau_n[t_1/x_1])]$.  For other list sorts it may have, see the definition of equivalence
of objects of metasorts {\bf arglist} $n$ given below [basically a rule for renaming bound variables $x_i$].

If $L$ is of metasort {\tt arglist} $n$ and $\tau$ is of metasort {\bf esort}, then $(L \Rightarrow \tau)$ is of metasort ${\bf asort}$.  All terms of this metasort are built in this way.
Notice that the input sorts found in $L$ may be either abstraction sorts or entity sorts, but the output sort is always an entity sort.   Note that {\bf asort} and {\bf esort} are disjoint.

If $t$ is a term of sort $\tau_1$ and $f$ is a term of sort $[(x_1,\tau_1)] \Rightarrow \tau$ (where $x_1$ does not occur in $t$), then $f[t]$ is a term of type $\tau[t/x_1]$.

If $[t_1,\ldots,t_n]$ is an argument list of list sort $[(x_1,\tau_1),\ldots,(x_n,\tau_n)]$  where no $t_j$ contains any $x_i$, and $f$ is a term of sort $[(x_1,\tau_1),\ldots,(x_n,\tau_n)] \Rightarrow \tau$
then $f[t_1,\ldots,t_n]$ is a term of the same sort as $g[t_2,\ldots,t_n]$, where $g$ is a variable of sort $$[(x_2,\tau_1[t_1/x_1]),\ldots,(x_n,\tau_n[t_1/x_1])] \Rightarrow \tau[t_1/x_1].$$

An equivalence relation on objects of type {\bf arglist} $n$ for each $n$ is defined:  $[(x_1,\tau_1)]$ is equivalent to any $[(x_1^*,\tau_1)]$.  $[(x_1^*,\tau_1^*),\ldots,(x_n^*,\tau_n^*)]$ is equivalent  to $[(x_1,\tau_1),\ldots,(x_n,\tau_n)]$ iff $\tau_1=\tau_1^*$ and $[(x_2^*,\tau_2^*),\ldots,(x_n^*,\tau_n^*)]$ is equivalent  to $$[(x_2,\tau_2[x_1^*/x_1]),\ldots,(x_n,\tau_n[x_1^*/x_1])],$$ where none of the starred variables occur in the unstarred term; the relation is then extended to be transitive.   Equivalent argument list sorts may freely replace one another in all contexts.  This amounts to the observation that the $x_i$'s are bound in this construction and can freely be renamed.  Similarly, $[(x_1,\tau_1),\ldots,(x_n,\tau_n)] \Rightarrow \tau$ is equivalent to $[(x_1^*,\tau_1^*),\ldots,(x_n^*,\tau_n^*)] \Rightarrow \tau^*$ under the same conditions under which any term
$[(x_1,\tau_1),\ldots,(x_n,\tau_n),(x_{n+1},\tau)]$ with $x_{n+1}$ a variable of correct type is equivalent to $$[(x_1^*,\tau_1^*),\ldots,(x_n^*,\tau_n^*),(x_{n+1},\tau^*)].$$  In particular, a term or list of any sort has all equivalent sorts as well.\footnote{An error in this description is most likely to have occurred here:  bound variable renaming is tricky!}

Note that application terms are always of sorts which are of metasort {\bf esort}.  Non-atomic terms of sorts of metasort {\bf asort} are of the shape $$[(x_1,\tau_1),\ldots,(x_n,\tau_n)] \Rightarrow (D,\tau),$$ where $D$ is a term of sort $\tau$.  $([(x_1,\tau_1),\ldots,(x_n,\tau_n)] \Rightarrow (D,\tau))[t_1,\ldots,t_n]$ reduces to $D[t_1/x_1][t_2/x_2]\ldots[t_n/x_n]$, if the argument list is of the correct argument list sort. The Lestrade user does not in the current version write any such terms (they are introduced implicitly by definitions of abstractions) or set up any such reductions, but they do appear in displayed types and such reductions happen in type computations.   Equivalence induced by renaming of bound variables is defined for these terms in essentially the same way as for abstraction sorts and argument list type sorts.

When the Lestrade engine computes a sort for a term or otherwise makes a substitution into a bound variable construction, it renames all bound variables in a way guaranteed to give fresh names.   When attempting to match types, the engine attempts definitional expansion of the types (and also rewriting of the types) if it does not initially see them as equivalent up to renaming of bound variables.

\noindent {\bf 10/10 change allowing more abstraction arguments:} The code of 10/10 while not implementing $\lambda$ terms per se permits easier formation of abstraction arguments.
A term $f(t_1,\ldots,t_m)$ where $m$ is less than the arity $n$ of $f$ will be read as $[(x_{m+1}\ldots x_n) \Rightarrow f(t_1,\ldots,t_m,x_{m+1},\dots,x_n)]$ if the types of the $t_i$'s are appropriate.  This only works if the argument list is explicitly closed with parentheses.  If $f$ is a defined abstraction, it will be definitionally expanded.   Such a term can only appear as an argument.
\newpage

\subsection{A sketch of semantics for a large class of  Lestrade theories}

The referent of a sort is a set.   The referent of a term of a particular sort will be an element of the referent of the sort.  We will refer to the sets which are referents of sorts ({\bf in} $\tau$) as ``types".  [we may further suppose, if we are willing to be boringly classical,  that the referent of {\tt prop} has exactly two elements and that types ({\bf that},$\tau$) are also sets, each such type having one element or none:  these sets may also be viewed as types].

The referent of an abstraction sort $((x_1,\tau_1),\ldots,(x_n,\tau_n) => ({\bf ---},\tau)$ is a space of functions whose range is the set which is the referent of $\tau$
and whose domain is a complex set of $n$-tuples:  the domain of the referent of $((x_1,\tau_1),\ldots,(x_n,\tau_n) => ({\bf ---},\tau)$ is the set of all $n$-element lists whose
whose head $t_1$ belongs to the referent of $\tau_1$ and whose tail belongs to the domain of the referent of $((x_2,\tau_2[t_1/x_1]),\ldots,(x_n,\tau_n[t_1/x_1]) => ({\bf ---},\tau[t_1/x_1])$.

The scope of this type system is better understood by adding dependent product constructions and dependent function space constructions:  where $\tau$ is a type and
$F$ is a function from the type $\tau$  to types, the dependent product $\tau \times F$ of $\tau$ and $F$ is the set of ordered pairs $(t,u)$ where $t \in \tau$ and $u \in F(t)$.   The dependent arrow type $\tau \rightarrow F$ is the set of functions $f$ with domain $\tau$ such that $f(t) \in F(t)$.   Iteration of dependent product will give referents of the domain types of 
all of our abstraction sorts:  the abstraction sorts themselves are dependent function types.

If the cardinality of the domain of types is taken to be inaccessible and each type is taken to be of smaller cardinality than the domain of types,  it is demonstrable that each dependent product type and dependent function type
whose first component is a type and values of whose second component are types is of lower cardinality than the domain of types.    Of course we can produce paradoxes by using {\tt type} as a component of a sort or by asserting existence of constructions which violate cardinality restrictions on types.

Lestrade declarations thus translate into assertion of existence of elements of certain collections of functions in an ambient set theory.

At the end of the sample book in section 6 we exhibit declarations of dependent product and dependent function types:   although the Lestrade user does not directly
type abstraction sorts, the user can develop theories in which these sorts are represented internally by types accessible to the user.  This could be done more faithfully to the
Lestrade sort system:  what is done in the book is a proof of concept.  It is very important to notice that the declarations of dependent product and dependent function types
are new axioms asserted in the logical framework, not consequences of the built-in features of the logical framework:  the logical framework cannot describe its own inner workings unaided.

One can formulate Lestrade theories which do not fit this semantics:  a constructive theory would have quite different semantics, for example.  Some interesting manipulations of
{\tt type} as a type which are used in current type systems would of course not fit with the naively set theoretical scheme of semantics described here.

\newpage

\section{Using Lestrade}

Lestrade is implemented currently in Moscow ML 2.01.  The files can be found on my web page (\cite{holmessource}).  I am planning to produce an implementation in Python 3 (there is a previous implementation in Python, which is deprecated:  it contains errors which I do not intend to correct -- instead I plan to port the ML version back into Python).  We also provide a source {\tt lestrade\_basic.sml} for a version without the rewriting system and other features using type inference which we plan to introduce:  the basic Lestrade features provide full support for our philosophical views and are less likely to be buggy.  I am planning to continue supporting both versions.

I have now installed a function basic() which if called from the ML command line will disable rewriting and implicit arguments, and a command explicit() which will disable implicit arguments but leave rewriting functions active.  These are not internal commands of Lestrade.  fullversion() will restore the usual behavior.

Features which we plan to introduce are user-entered $\lambda$-abstractions (entered as arguments, and only when types of the bound variables are fully deducible from the body of the abstraction).  A blue-sky idea I am thinking about is the implementation of imperative programming constructions:  it does not seem impossible that this can be done in a type-safe manner, and it would be very interesting if it could be done.

Lestrade books are text files with the extension {\tt .lti}.  Edit them with a text editor (fill them with lines in the format described in section 4) and the Lestrade checker should be able to do  something with them.

The basic file handling command of Lestrade is

{\tt readfile "<source file>" "<target file>";}

(to be typed on the ML command line!)

which executes the Lestrade lines in the file {\tt <source file>.lti} and echoes the commands and any responses from the system to {\tt target file>.lti}.  After running the file,
the user can type Lestrade lines in the interface and receive immediate feedback both at the console and echoed into the target file, though our usual approach is to edit the source file
and issue the readfile command again.  To exit the interface, type {\tt quit}.  To avoid poisonous problems with execution of readfile, always end a file of Lestrade commands
run with readfile with the line {\tt quit}.  Execution of files can be chained if subsequent files start with a {\tt load} command:  in this case it is important that file names be Lestrade tokens!

One can also type 

{\tt interface "<target file>";} to enter Lestrade lines at the console and echo the results to the given target file.  If the null string is used, no target file is used (or at least none is intended).

If the source
file contains valid Lestrade commands, the target file will in fact be executable as a source file with the same effects, and will in addition be better formatted and commented
with such things as the sorts of all the terms declared.

If an error is encountered while reading a file, Lestrade will stop and issue an error message.  Hitting return will scroll through any further error messages, and Lestrade will eventually stop the reading of the file and leave the user  in the interface:  it will not continue to execute commands in a file after an error is encountered.

A file should always end with the line {\tt quit}.

\newpage

\section{Distinctive features of our approach, and vague philosophical speculations}

The general fact that mathematical reasoning and construction of mathematical objects can be managed using a type checking system of this kind is already well known, from work with other systems (Automath, Coq and their relations),
and examples of constructions of both kinds under Lestrade appear in section 6.

There are some distinctive features of our approach.

One feature which the reader and the operator of Lestrade will notice is that the user never writes any abstraction sort or $\lambda$-term.   All abstractions which are introduced are assigned names, as if we introduced all functions in the style $f(x) = x^2+1$ instead of the style $f=(\lambda x.x^2+1)$ or $f=(x \mapsto x^2+1)$.  I thought originally
that no $\lambda$-terms were needed at all, but this is not the case:   it is easy to get into a situation where the name of an abstraction passes out of scope when a move is closed,
but the abstraction appears as an argument in a sort, and so is expanded to a $\lambda$-term.  An abstraction can also appear in applied position and pass out of scope,
in which case $\beta$-reduction will occur automatically:  a $\lambda$-term can only appear as an argument, never in applied position.

This could be changed, but so far I have not felt compelled to do so.  I have found it interesting working on a style in which abstractions which might appear anonymously as scopes of variable binding constructions have to be explicitly set up and assigned names in advance.   The need to use a lot of names is limited by the fact that identifiers are regularly freed up when moves are closed.  Fluent handling of abstraction arguments may be improved by the recent introduction of abstractions applied to shortened argument lists, interpreted as complex abstraction terms via currying.

There is a limitation of my type system.   The output sorts of all abstractions are entity sorts (this is superficially modified by the new device of curried abstraction arguments).  This I have no intention of changing:  I take the view that abstractions are not {\em prima facie} first-class objects, and attempts to convert them into entities (first-class objects) must involve postulation of suitable constructions by the user.   For example, the following theory introduces a very powerful ability to code abstractions using entities, which leads to Russell's paradox.   The approach I take to implementing type theory of sets in the previous example has some similarities.

\newpage

\begin{verbatim}
>>Inspector Lestrade says:
  Welcome to the Lestrade Type Inspector,
 version of 6/21/2016 11:15 am Boise time

open
     declare x obj

>>     x:  obj

     construct P x:prop

>>     P:  [(x_1:obj) => (---:prop)]

     close

\end{verbatim}

The constructor {\tt set} is postulated to cast predicates of type {\tt obj} to entities of type {\tt obj}.  The fact that postulating an object to correspond to each predicate is a nontrivial logical move is concealed in the intuitive argument for Russell's ``paradox"  and is in my view the reason why it is a mistake, not a paradox.  The mere existence of {\tt set} is not enough to derive the paradox:  more dubious assumptions must be made, and duly will be!

\begin{verbatim}
construct set P:obj

>>set:  [(P_1:[(x_2:obj) => (---:prop)]) => (---:obj)]

declare x obj

>>x:  obj

declare y obj

>>y:  obj

\end{verbatim}

\newpage

The membership relation is postulated.

\begin{verbatim}

construct E x y:prop

>>E:  [(x_1:obj),(y_1:obj) => (---:prop)]

\end{verbatim}

Here is the crux of the mistake.  We postulate the comprehension axioms, the flat-footed assertion of a one-to-one correspondence, for each predicate $P$, between evidence
for $P(x)$ and evidence for $x \in \{x:P(x)\}$.

\begin{verbatim}

declare x1 that P x

>>x1:  that P(x)

construct comp P, x x1:that E x set P

>>comp:  [(P_1:[(x_2:obj) => (---:prop)]),(x_1:obj),(x1_1:that 
>>   P_1(x_1)) => (---:that (x_1 E set(P_1)))]

declare x2 that E x set P

>>x2:  that (x E set(P))

construct comp2 P, x x2:  that P x

>>comp2:  [(P_1:[(x_2:obj) => (---:prop)]),(x_1:obj),
>>   (x2_1:that (x_1 E set(P_1))) => (---:that P_1(x_1))]

declare p prop

>>p:  prop

declare q prop

>>q:  prop

\end{verbatim}

To complete the execution of our folly, we need to declare the logical operations of implication and negation in a familiar way.  One should note that these declarations  are fine from a constructive standpoint.

\begin{verbatim}

construct Implies p q:prop

>>Implies:  [(p_1:prop),(q_1:prop) => (---:prop)]

construct False:prop

>>False:  [(---:prop)]

declare pp that p

>>pp:  that p

declare rr that Implies p q

>>rr:  that (p Implies q)

construct Mp p q pp rr:that q

>>Mp:  [(p_1:prop),(q_1:prop),(pp_1:that p_1),(rr_1:that 
>>   (p_1 Implies q_1)) => (---:that q_1)]

declare absurd that False

>>absurd:  that False

construct Panic p absurd: that p

>>Panic:  [(p_1:prop),(absurd_1:that False) => (---:that 
>>   p_1)]

define Not p:Implies p False

>>Not:  [(p_1:prop) => ((p_1 Implies False):prop)]

open
     declare pp2 that p

>>     pp2:  that p

     construct Ded pp2:that q

>>     Ded:  [(pp2_1:that p) => (---:that q)]

     close
construct Impliesproof p q Ded:that Implies p q

>>Impliesproof:  [(p_1:prop),(q_1:prop),(Ded_1:[(pp2_2:that 
>>   p_1) => (---:that q_1)]) => (---:that (p_1 Implies 
>>   q_1))]

\end{verbatim}

{\tt Russell} is the Russell predicate $(\lambda x:x \not\in x\}$ while $R$ is the Russell (paradoxical) set.  The argument which follows
is familiar.  {\tt R6} is a proof of the False in move 0, a disaster.

\begin{verbatim}

define Russell x:Not E x x

>>Russell:  [(x_1:obj) => (Not((x_1 E x_1)):prop)]

open
     define R: set Russell

>>     R:  [(set(Russell):obj)]

     declare R1 that E set Russell, set Russell

>>     R1:  that (set(Russell) E set(Russell))

     define R2 R1:comp2 Russell, set Russell, R1

>>     R2:  [(R1_1:that (set(Russell) E set(Russell))) 
>>        => (comp2(Russell,set(Russell),R1_1):that Russell(set(Russell)))]

     define R3 R1:Mp E set Russell, set Russell, False R1 R2 R1

>>     R3:  [(R1_1:that (set(Russell) E set(Russell))) 
>>        => (Mp((set(Russell) E set(Russell)),False,R1_1,
>>        R2(R1_1)):that False)]

     close
define R4:Impliesproof E set Russell, set Russell, False R3

>>R4:  [(Impliesproof((set(Russell) E set(Russell)),False,
>>   [(R1_1:that (set(Russell) E set(Russell))) => (Mp((set(Russell) 
>>   E set(Russell)),False,R1_1,comp2(Russell,set(Russell),
>>   R1_1)):that False)]):that ((set(Russell) E set(Russell)) 
>>   Implies False))]

define R5:comp Russell, set Russell, R4

>>R5:  [(comp(Russell,set(Russell),R4):that (set(Russell) 
>>   E set(Russell)))]

define R6: Mp E set Russell, set Russell, False R5 R4

>>R6:  [(Mp((set(Russell) E set(Russell)),False,R5,R4):that 
>>   False)]

quit

\end{verbatim}

\newpage

An important thing to notice about the crucial last few lines of the argument is that Lestrade is automatically expanding definitions as it type checks these lines.  I should add comments on each line pointing out the definitional expansions.

I hope you enjoyed that!

Usually the proofs of implications are identified with actual functions under the Curry-Howard isomorphism.    My requirement that user defined constructions have entity
output requires that objects {\tt Ifproof p q Ded} are obtained by a user-declared construction from functions  {\tt Ded}, not identified with functions.  This is useful.
There is no need for proofs of implications $p \rightarrow q$ to have the same identity conditions as functions from proofs of $p$ to proofs of $q$, and one can
create situations where there are ``too many" distinct proofs if the proofs are actually functions.  [I seem to recall that universal quantifiers over {\tt prop} lead to some weirdness.]

This goes along with the idea that treating abstractions as entities requires one to declare constructions that effect this reduction.  This casts a kind of light on reasons behind the ``paradoxes of set theory".

Another limitation of the type system is that there are no disjoint union types or existential types ``built in" (though they certainly can be declared).  In this Lestrade is similar to Automath (of which it is arguably a flavor).  Essentially the only built in type constructor of Lestrade is the dependent type construction of functions, of which the usual Curry-Howard implementations of implication and universal quantification are examples.  This gives the implementations of disjunction and the existential quantifier as logical operations an indirect character [much as in Automath].

We have no constructive prejudices, though we observe that Lestrade supports constructive logic perfectly well with some modifications of the basic declarations of logical operations.  We are interested in a more constructive approach for other reasons:  we would like to make an obvious further extension and make Lestrade a programming environment in which programs can be written which will operate effectively on the wide range of mathematical objects which its rich [and easily user-extendible!]  system of sorts makes accessible.  The rewriting features added recently go some way toward implementing this.

Another class of philosophical objections to classical mathematics is addressed by our approach here.   We have no sympathy with predicativist scruples and we are very fond of second-order logic (a logic which supports quantification over universals).   One should note in the Lestrade book that we showed how to quantify over untyped objects, then how to quantify uniformly over any sort ({\bf in} $\tau$), then how to quantify universally over binary relations (abstractions!) from natural numbers to an arbitry type.  This last is an example of second order quantification.

On the other hand, Lestrade does not automatically commit to the ability to quantify over any type, as Automath does, at least in later versions, due to the fact that functions with output in {\tt prop} or {\tt type} are also typed as instances of {\tt prop} or {\tt type} respectively.  It appears to be possible to restrict oneself to the resources of first-order reasoning by suitable choice of primitives.

The induction axiom on the natural number type defined in the book works on predicates of natural numbers defined in any way that any extension of the book might provide for.
We argue that the intended referent of this type is the true type of natural numbers with a second-order axiomatization.  We intend to extend the book to present an axiomatization of Peano arithmetic for contrast, for which a Lestrade specification of first-order logic formulas over the Peano naturals would have to be given, and induction provided only for properties of the Peano naturals expressible by such formulas.  One could then postulate other predicates of the Peano naturals which did not respect this induction axiom, thus supporting reasoning about ``nonstandard" natural numbers which would not be possible for the type of natural numbers currently described.

Similarly, I would like to present first-order versus second-order Zermelo set theory and ZF.

Predicativist scruples seem to us to arise from a disagreement about the nature of generality.  For us, a function is not an infinite table of values.  I do not need to be acquainted with
every natural number $n$ and compute $n^2+1$ for each of them to be acquainted with the function $f(n) = n^2+1$:  I have to know the method of computation.  I do not need
to know every element of the domain $D$ to be able to prove $(\forall x \in D:\phi(x))$:  I need a function which given an $x \in D$ will generate a proof of $\phi(x)$ (note the dependent typing), and this function may be a finite object in the same sense that the expression $n^2+1$ is.

I do not regard the notion ``set of natural number" as vague because I do not have a way of effectively listing all sets of natural numbers.   A set of natural numbers is a gadget which given a natural number input gives a propositional output:  I can recognize such an object without having a clue as to how many such objects there are.  I do not need any familiarity with the full extent of the domain of sets of natural numbers to be able to prove a universal statement about it.   In particular, if I am given zero and a successor operation in a domain which may be supposed to properly include (an implementation of) the natural numbers, I can present a uniform proof that 3 (defined in the obvous way) belongs to every inductive set  as a ``finite gadget".  And I can abstract from this to the notion that 3 is a natural number, defined in the usual impredicative way.  I do not regard universal quantifications as infinitary conjunctions which require for their understanding that
$\phi(x)$ be previously understood for each $x \in D$:  such an understanding does make impredicative definitions circular.  But note that it is obvious that we cannot possibly be thinking that way about such sentences in practice:  our understanding of universal quantification, which makes reasoning about it clearly possible as a finite act, also dispels the apparently problematic character of impredicative reasoning (without removing the signs that impredicative definition  is a very powerful move).

I'm well aware that philosophical speculations can be apparently vague and unsatisfactory.  One of our aims in designing Lestrade is to create an environment in which one has as it were hands-on access to the full range of mathematical abstractions, in the incontrovertible sense that one can actually design the objects and execute proofs of theorems  about them.  An environment in which very general mathematical objects can be manipulated precisely should be an aid to philosophical contemplation of them.  The claim that this environment does so implement the mathematical objects and proofs is itself a definite philosophical claim, and interaction with the software should make it easier to understand and evaluate this claim.

\newpage

\section{A moderately extensive Lestrade book}

To see how this works, we present a moderately extensive development of some basic logical and mathematical concepts in Lestrade.  The considerations in section 2 should give a general idea of what is going on:  details of syntax and command format in Lestrade are given in section 4, and one may look forward to that point to see details of how the book is to be parsed and executed.

This can be compared with text produced under Automath as in \cite{automathtext} and \cite{grundlagen}.  The development of arithmetic by Jutting in \cite{grundlagen} (following Landau's classic \cite{landau}) is available from Freek Wiedijk at the same place as his Automath implementation, \cite{freek}.

This book makes no use of the new implicit argument feature.  I ought to introduce some declarations with implicit arguments and examples of their use (there is a short section of this kind above).

\subsection{Propositional logic of conjunction and implication}

\begin{verbatim}


>>Inspector Lestrade says:  
Welcome to the Lestrade Type Inspector, 
version of 6/21/2016 11:15 am Boise time

declare p prop

>>p:  prop

declare q prop

>>q:  prop

declare pp that p

>>pp:  that p

declare qq that q

>>qq:  that q

\end{verbatim}

We declare the conjunction operation on propositions.

\begin{verbatim}

comment Declare the conjunction operator
construct & p q : prop

>>&:  [(p_1:prop),(q_1:prop) => (---:prop)]

\end{verbatim}

We present the rule of conjunction introduction as a mathematical object, taking the two propositions and proofs of each of them as arguments
and returning a proof of the conjunction.

Notice that while we always declare an abstraction with two arguments using prefix notation, it will be displayed in infix notation (as long as its first argument is not an abstraction), as conjunction is here.
The parser can handle mixfix notation for abstractions with three or more arguments but Lestrade does not choose to display things in this way.

\begin{verbatim}

comment The rule of conjunction
construct Andproof p q pp qq : that p & q

>>Andproof:  [(p_1:prop),(q_1:prop),(pp_1:that p_1),(qq_1:that 
>>   q_1) => (---:that (p_1 & q_1))]

declare rr that p & q

>>rr:  that (p & q)

\end{verbatim}

Similarly, we present the rules of conjunction elimination (simplication) as mathematical objects.

\begin{verbatim}

comment The rules of simplification
construct And1 p q rr :  that p

>>And1:  [(p_1:prop),(q_1:prop),(rr_1:that (p_1 & q_1)) 
>>   => (---:that p_1)]

construct And2 p q rr :  that q

>>And2:  [(p_1:prop),(q_1:prop),(rr_1:that (p_1 & q_1)) 
>>   => (---:that q_1)]

\end{verbatim}

We declare implication just as we declared conjunction.

\begin{verbatim}

comment The implication operator
construct -> p q : prop

>>->:  [(p_1:prop),(q_1:prop) => (---:prop)]

\end{verbatim}

We develop the rule of conditional proof (the deduction theorem) as a mathematical object.  This is more exciting, because one of its inputs is a function.

\begin{verbatim}

comment Development of conditional proof
open
     declare pp2 that p

>>     pp2:  that p

     comment Ded below does not need p or q in its argument list 
     comment because they are not locally variables.
     construct Ded pp2 : that q

>>     Ded:  [(pp2_1:that p) => (---:that q)]

     close
comment Note that once the move at which Ded was constructed closes, 
comment it is a variable of desirable function type
construct Ifproof p q Ded : that p -> q

>>Ifproof:  [(p_1:prop),(q_1:prop),(Ded_1:[(pp2_2:that 
>>   p_1) => (---:that q_1)]) => (---:that (p_1 -> q_1))]

\end{verbatim}

We demonstrate our powers by actually proving a theorem ($P \rightarrow P$).

\begin{verbatim}

comment Now, for fun, we will construct an actual proof
open
     declare pp2 that p

>>     pp2:  that p

     define Ppid pp2 : pp2

>>     Ppid:  [(pp2_1:that p) => (pp2_1:that p)]

     close
define Selfimp p : Ifproof p p Ppid

>>Selfimp:  [(p_1:prop) => (Ifproof(p_1,p_1,[(pp2_2:that 
>>   p_1) => (pp2_2:that p_1)]):that (p_1 -> p_1))]

comment Notice in the sort of Selfimp that Ppid has
comment been expanded as a lambda-term.
comment Develop the rule of modus ponens
declare ss that p -> q

>>ss:  that (p -> q)

\end{verbatim}

We complete the basics of implication by defining the rule of modus ponens as a mathematical object.

\begin{verbatim}

construct Mp p q pp ss : that q

>>Mp:  [(p_1:prop),(q_1:prop),(pp_1:that p_1),(ss_1:that 
>>   (p_1 -> q_1)) => (---:that q_1)]

\end{verbatim}

\subsection{The universal quantifier}

In this section we do the basic development of the universal quantifier.   After this we will return to propositional logic and after that to a little more quantification.

The argument of the universal quantifier is a function, which requires setup here similar to that for the conditional proof rule above.

\begin{verbatim}

comment Opening an environment to set up definition of a predicate variable P
open
     declare xx obj

>>     xx:  obj

     construct P xx : prop

>>     P:  [(xx_1:obj) => (---:prop)]

     close
comment Declaring the universal quantifier
construct Forall P : prop

>>Forall:  [(P_1:[(xx_2:obj) => (---:prop)]) => (---:prop)]

\end{verbatim}

We declare the rule of universal instantiation.  It is worth noting that the parser sometimes requires us to guard abstractions appearing as arguments with commas
so that they will not be applied to what follows them (or what is before them if they are capable of being read as infix or mixifx operators).

\begin{verbatim}

comment Declaring the rule UI of universal instantiation
declare P2 that Forall P

>>P2:  that Forall(P)

declare x obj

>>x:  obj

construct Ui P, P2 x : that P x

>>Ui:  [(P_1:[(xx_2:obj) => (---:prop)]),(P2_1:that Forall(P_1)),
>>   (x_1:obj) => (---:that P_1(x_1))]

comment Note in the previous line that we follow P 
comment with a comma:  an abstraction argument may need to be 
comment guarded with commas so it will not be read as applied.
comment Opening an environment to declare a function 
comment that witnesses provability of a universal statement

\end{verbatim}

We declare the rule of universal quantifier introduction.  Note that the second argument is a function which takes an object $x$ to a proof of $P(x)$, a dependently typed function witnessing the truth of the universal statement.

\begin{verbatim}

open
     declare u obj

>>     u:  obj

     construct Q2 u : that P u

>>     Q2:  [(u_1:obj) => (---:that P(u_1))]

     close
comment The rule of universal generalization
construct Ug P, Q2 : that Forall P

>>Ug:  [(P_1:[(xx_2:obj) => (---:prop)]),(Q2_1:[(u_3:obj) 
>>   => (---:that P_1(u_3))]) => (---:that Forall(P_1))]

\end{verbatim}

\subsection{Negation (made classical) interacts with implication:  proof of the contrapositive theorem}

In this section we introduce $\perp$ ({\tt ??}), a false statement, as a primitive, then introduce logical negation and the biconditional as defined notions.  We then prove the contrapositive theorem and develop derived rules relating implication and negation.  We make our logic classical by declaring the rule of double negation, but a constructive approach is certainly possible.

\begin{verbatim}

comment Develop rules for negation (which will be classical!) 
    and prove the contrapositive theorem.
comment The absurd proposition.
construct ??:prop

>>??:  [(---:prop)]

comment The negation operation.
define ~p:  p -> ??

>>~:  [(p_1:prop) => ((p_1 -> ??):prop)]

\end{verbatim}

Here we introduce the primitive that makes our logic classical.

\begin{verbatim}

comment We make our logic classical:  the rule of double negation
declare maybe that ~ ~p

>>maybe:  that ~(~(p))

construct Dneg p maybe : that p

>>Dneg:  [(p_1:prop),(maybe_1:that ~(~(p_1))) => (---:that 
>>   p_1)]

\end{verbatim}

Here we show that contradictions in the usual sense of the term imply the absurd primitive.  It is worth noting that Lestrade
does recognize that a negative statement has the form of an implication and applies modus ponens without any need for special action to unpack the definition of negation.

\begin{verbatim}

comment Contradictions are absurd.
declare nn that ~p

>>nn:  that ~(p)

define Contra p pp nn :  Mp p ?? pp nn

>>Contra:  [(p_1:prop),(pp_1:that p_1),(nn_1:that ~(p_1)) 
>>   => (Mp(p_1,??,pp_1,nn_1):that ??)]

comment Notice that Lestrade does expand the definition
comment of the negation operation as we expect.

\end{verbatim}

We start the development of the rule of negation introduction.   We have to do a little extra work, because the most direct approach gives us a rule
in which the output is typed in expanded form as an implication.   But this can be fixed.

\begin{verbatim}

open
     declare pp2 that p

>>     pp2:  that p

     construct Negded pp2: that ??

>>     Negded:  [(pp2_1:that p) => (---:that ??)]

     close
define Negintro1 p Negded :  Ifproof p ?? Negded

>>Negintro1:  [(p_1:prop),(Negded_1:[(pp2_2:that p_1) 
>>   => (---:that ??)]) => (Ifproof(p_1,??,Negded_1):that 
>>   (p_1 -> ??))]

comment Negation introduction.  But it is defective in actually 
comment reporting an implication type.  Let's see if we can fix this.
open
     declare proof1 that p -> ??

>>     proof1:  that (p -> ??)

     define Negproofid proof1:proof1

>>     Negproofid:  [(proof1_1:that (p -> ??)) => (proof1_1:that 
>>        (p -> ??))]

     close
define Negfix p :  Ifproof ((p -> ??), ~p , Negproofid)

>>Negfix:  [(p_1:prop) => (Ifproof((p_1 -> ??),~(p_1),
>>   [(proof1_2:that (p_1 -> ??)) => (proof1_2:that (p_1 
>>   -> ??))]):that ((p_1 -> ??) -> ~(p_1)))]

define Negintro p Negded : Mp ((p -> ??), ~p , Negintro1 p Negded , Negfix p)

>>Negintro:  [(p_1:prop),(Negded_1:[(pp2_2:that p_1) 
>>   => (---:that ??)]) => (Mp((p_1 -> ??),~(p_1),(p_1 Negintro1 
>>   Negded_1),Negfix(p_1)):that ~(p_1))]

comment I succeed in defining a proper negation introduction rule
comment using the defined symbol.  This is important because of limitations
comment of circumstances under which Lestrade expands definitions.

\end{verbatim}

We define the biconditional and introduce its rules.  Of course since it is a defined operation we do not need
to declare any new primitives in this connection.

\begin{verbatim}

comment We define the biconditional.
define <-> p q : (p -> q) & (q -> p)

>><->:  [(p_1:prop),(q_1:prop) => (((p_1 -> q_1) & (q_1 
>>   -> p_1)):prop)]

\end{verbatim}

The biconditional elimination rules are variations of modus ponens.

\begin{verbatim}

comment Biconditional elimination rules
declare tt that p <-> q

>>tt:  that (p <-> q)

define Mpb1 p q pp tt : Mp p q pp, And1 ((p -> q), (q -> p), tt)

>>Mpb1:  [(p_1:prop),(q_1:prop),(pp_1:that p_1),(tt_1:that 
>>   (p_1 <-> q_1)) => (Mp(p_1,q_1,pp_1,And1((p_1 -> q_1),
>>   (q_1 -> p_1),tt_1)):that q_1)]

define Mpb2 p q qq tt : Mp q p qq, And2((p->q),(q->p),tt)

>>Mpb2:  [(p_1:prop),(q_1:prop),(qq_1:that q_1),(tt_1:that 
>>   (p_1 <-> q_1)) => (Mp(q_1,p_1,qq_1,And2((p_1 -> q_1),
>>   (q_1 -> p_1),tt_1)):that p_1)]

comment In both of the last two commands, there are subtle parser issues.
comment Before And1, And2, the comma is needed to prevent Andi 
comment from being read as an infix.
comment Because we enclose the argument (p->q) in parentheses
comment we need to enclose the entire argument list in parentheses
comment because a parenthesis after a prefixed abstraction is
comment always interpreted as enclosing an argument list,
comment not a term.
comment the classic Reductio ad Absurdum (which is not the same as neg intro!)

\end{verbatim}

We develop the rule of {\em reductio ad absurdum\/} (which is not the same as negation introduction, though both are carelessly called
proof by contradiction) and the rule that anything can be deduced from a falsehood.

\begin{verbatim}

open
     declare aa that ~p

>>     aa:  that ~(p)

     construct reductioarg aa :  that ??

>>     reductioarg:  [(aa_1:that ~(p)) => (---:that ??)]

     close
define Reductio p reductioarg : Dneg p (Negintro ~p reductioarg)

>>Reductio:  [(p_1:prop),(reductioarg_1:[(aa_2:that ~(p_1)) 
>>   => (---:that ??)]) => ((p_1 Dneg (~(p_1) Negintro reductioarg_1)):that 
>>   p_1)]

comment Everything follows from the False!
declare huh that ??

>>huh:  that ??

open
     declare negp that ~p

>>     negp:  that ~(p)

     define panick negp :  huh

>>     panick:  [(negp_1:that ~(p)) => (huh:that ??)]

     close
define Panic p huh :  Reductio p panick

>>Panic:  [(p_1:prop),(huh_1:that ??) => ((p_1 Reductio 
>>   [(negp_2:that ~(p_1)) => (huh_1:that ??)]):that p_1)]

\end{verbatim}

We develop the biconditional introduction rule.  This is similar to the rule of conditional proof.  We once again have to do a little extra work to get
an output which is actually typed as a biconditional rather than as a conjunction of implications.

\begin{verbatim}

comment We develop the biconditional introduction rule.
comment In this environment we postulate reasoning
comment leading from p to q and q to p
open
     declare pp2 that p

>>     pp2:  that p

     construct Ded1 pp2: that q

>>     Ded1:  [(pp2_1:that p) => (---:that q)]

     declare qq2 that q

>>     qq2:  that q

     construct Ded2 qq2: that p

>>     Ded2:  [(qq2_1:that q) => (---:that p)]

     close
comment Here we prove an initial version,
comment defective in having expanded output
define Biintro1 p q, Ded1, Ded2: 
Andproof ((p->q),(q->p),Ifproof p q Ded1,Ifproof q p Ded2)

>>Biintro1:  [(p_1:prop),(q_1:prop),(Ded1_1:[(pp2_2:that 
>>   p_1) => (---:that q_1)]),(Ded2_1:[(qq2_3:that q_1) 
>>   => (---:that p_1)]) => (Andproof((p_1 -> q_1),(q_1 
>>   -> p_1),Ifproof(p_1,q_1,Ded1_1),Ifproof(q_1,p_1,Ded2_1)):that 
>>   ((p_1 -> q_1) & (q_1 -> p_1)))]

open
     declare bb that p <-> q

>>     bb:  that (p <-> q)

     define biid bb:bb

>>     biid:  [(bb_1:that (p <-> q)) => (bb_1:that (p 
>>        <-> q))]

     close
comment We fix the defective version much as we fixed Negintro above
define Bifix p q: Ifproof (((p->q) & (q->p)),p<->q,biid)

>>Bifix:  [(p_1:prop),(q_1:prop) => (Ifproof(((p_1 -> 
>>   q_1) & (q_1 -> p_1)),(p_1 <-> q_1),[(bb_2:that (p_1 
>>   <-> q_1)) => (bb_2:that (p_1 <-> q_1))]):that (((p_1 
>>   -> q_1) & (q_1 -> p_1)) -> (p_1 <-> q_1)))]

define Biintro p q, Ded1, Ded2: 
Mp (((p->q)&(q->p)),p<->q,Biintro1 (p, q, Ded1, Ded2),Bifix p q)

>>Biintro:  [(p_1:prop),(q_1:prop),(Ded1_1:[(pp2_2:that 
>>   p_1) => (---:that q_1)]),(Ded2_1:[(qq2_3:that q_1) 
>>   => (---:that p_1)]) => (Mp(((p_1 -> q_1) & (q_1 -> 
>>   p_1)),(p_1 <-> q_1),Biintro1(p_1,q_1,Ded1_1,Ded2_1),
>>   (p_1 Bifix q_1)):that (p_1 <-> q_1))]

\end{verbatim}

We prove the contrapositive theorem.  The proof follows the structure of the proof using  my favorite natural deduction strategy for propositional logic exactly.

\begin{verbatim}

comment We prove the contrapositive theorem,
comment (p->q) <-> (~q <-> ~p)
open
     declare aa that p->q

>>     aa:  that (p -> q)

     comment Our goal is to construct a proof of ~q -> ~p
     comment To do this, we need a function from
     comment proofs of ~q to proofs of ~p
     open
          declare bb that ~q

>>          bb:  that ~(q)

          comment Now our goal is to prove ~p.
          comment For this we need a function from 
          comment proofs of p to proofs of ??
          open 
               declare cc that p

>>               cc:  that p

               comment prove q by m.p.
               define dd cc: Mp p q cc aa

>>               dd:  [(cc_1:that p) => (Mp(p,q,cc_1,
>>                  aa):that q)]

               comment and we have a contradiction
               define ee cc: Contra q (dd cc) bb

>>               ee:  [(cc_1:that p) => (Contra(q,dd(cc_1),
>>                  bb):that ??)]

               close
          define ff bb :  Negintro p ee

>>          ff:  [(bb_1:that ~(q)) => ((p Negintro [(cc_2:that 
>>             p) => (Contra(q,Mp(p,q,cc_2,aa),bb_1):that 
>>             ??)]):that ~(p))]

          close
     define gg aa:  Ifproof ((~q),(~p),ff)

>>     gg:  [(aa_1:that (p -> q)) => (Ifproof(~(q),~(p),
>>        [(bb_2:that ~(q)) => ((p Negintro [(cc_3:that 
>>        p) => (Contra(q,Mp(p,q,cc_3,aa_1),bb_2):that ??)]):that 
>>        ~(p))]):that (~(q) -> ~(p)))]

     comment Now we need the function acting in
     comment the other direction
     declare hh that ~q -> ~p

>>     hh:  that (~(q) -> ~(p))

     comment Our goal is p->q so we want to assume p
     open
          declare ii that p

>>          ii:  that p

          comment Now our goal is q, but we will
          comment actually aim for ~~q and so
          comment assume ~q
          open
               declare jj that ~q

>>               jj:  that ~(q)

               comment Now use modus ponens to prove ~p
               define kk jj :  Mp(~q,~p,jj,hh)

>>               kk:  [(jj_1:that ~(q)) => (Mp(~(q),~(p),
>>                  jj_1,hh):that ~(p))]

               comment Now we have a contradiction
               define ll jj : Contra p ii kk jj

>>               ll:  [(jj_1:that ~(q)) => (Contra(p,
>>                  ii,kk(jj_1)):that ??)]

               close
          define mm ii :  Negintro (~q , ll)

>>          mm:  [(ii_1:that p) => ((~(q) Negintro [(jj_2:that 
>>             ~(q)) => (Contra(p,ii_1,Mp(~(q),~(p),jj_2,
>>             hh)):that ??)]):that ~(~(q)))]

          define nn2 ii : Dneg q mm ii

>>          nn2:  [(ii_1:that p) => ((q Dneg mm(ii_1)):that 
>>             q)]

          close
     define oo hh :  Ifproof p q nn2

>>     oo:  [(hh_1:that (~(q) -> ~(p))) => (Ifproof(p,
>>        q,[(ii_2:that p) => ((q Dneg (~(q) Negintro [(jj_3:that 
>>        ~(q)) => (Contra(p,ii_2,Mp(~(q),~(p),jj_3,hh_1)):that 
>>        ??)])):that q)]):that (p -> q))]

     close



define Contrapositive p q:  Biintro ((p->q),(~q -> ~p),gg,oo)

>>Contrapositive:  [(p_1:prop),(q_1:prop) => (Biintro((p_1 
>>   -> q_1),(~(q_1) -> ~(p_1)),[(aa_2:that (p_1 -> q_1)) 
>>   => (Ifproof(~(q_1),~(p_1),[(bb_3:that ~(q_1)) => ((p_1 
>>   Negintro [(cc_4:that p_1) => (Contra(q_1,Mp(p_1,q_1,
>>   cc_4,aa_2),bb_3):that ??)]):that ~(p_1))]):that (~(q_1) 
>>   -> ~(p_1)))],[(hh_5:that (~(q_1) -> ~(p_1))) => (Ifproof(p_1,
>>   q_1,[(ii_6:that p_1) => ((q_1 Dneg (~(q_1) Negintro 
>>   [(jj_7:that ~(q_1)) => (Contra(p_1,ii_6,Mp(~(q_1),~(p_1),
>>   jj_7,hh_5)):that ??)])):that q_1)]):that (p_1 -> q_1))]):that 
>>   ((p_1 -> q_1) <-> (~(q_1) -> ~(p_1))))]

comment Now is a good point to notice that
comment Lestrade definitely saves proof objects in detail.

\end{verbatim}

We develop the derived logical rules which mix implication and negation, modus tollens and proof by contrapositive.

\begin{verbatim}

comment Develop indirect proof strategies for implication.
comment Modus Tollens
declare negc that ~q

>>negc:  that ~(q)

define Mt p q ss negc :  Mp(~q, ~p, negc ,
 Mpb1 ((p -> q),(~q -> ~p),ss,Contrapositive p q))

>>Mt:  [(p_1:prop),(q_1:prop),(ss_1:that (p_1 -> q_1)),
>>   (negc_1:that ~(q_1)) => (Mp(~(q_1),~(p_1),negc_1,Mpb1((p_1 
>>   -> q_1),(~(q_1) -> ~(p_1)),ss_1,(p_1 Contrapositive 
>>   q_1))):that ~(p_1))]

comment Rule of contrapositive or indirect proof
open 
     declare negq that ~q

>>     negq:  that ~(q)

     construct indarg negq :  that ~p

>>     indarg:  [(negq_1:that ~(q)) => (---:that ~(p))]

     close
define Indirect p q indarg :
 Mpb2 ((p->q),(~q -> ~p),Ifproof (~q,~p,indarg),Contrapositive p q)

>>Indirect:  [(p_1:prop),(q_1:prop),(indarg_1:[(negq_2:that 
>>   ~(q_1)) => (---:that ~(p_1))]) => (Mpb2((p_1 -> q_1),
>>   (~(q_1) -> ~(p_1)),Ifproof(~(q_1),~(p_1),indarg_1),
>>   (p_1 Contrapositive q_1)):that (p_1 -> q_1))]

\end{verbatim}

\subsection{The development of disjunction}

We declare the disjunction operation and introduce its constructively valid rules (addition and proof by cases)  then derive the more powerful rules mixing disjunction and negation.

\begin{verbatim}

comment Now start the development of disjunction.
comment disjunction declared
construct v p q:prop

>>v:  [(p_1:prop),(q_1:prop) => (---:prop)]

comment basic disjunction introduction rules (addition)
construct Addition1 p q pp: that p v q

>>Addition1:  [(p_1:prop),(q_1:prop),(pp_1:that p_1) 
>>   => (---:that (p_1 v q_1))]

construct Addition2 p q qq:that p v q

>>Addition2:  [(p_1:prop),(q_1:prop),(qq_1:that q_1) 
>>   => (---:that (p_1 v q_1))]

comment the basic disjunction elimination rule (proof by cases)
declare r prop

>>r:  prop

declare disj that p v q

>>disj:  that (p v q)

open
     declare pp2 that p

>>     pp2:  that p

     construct case1 pp2 : that r

>>     case1:  [(pp2_1:that p) => (---:that r)]

     declare qq2 that q

>>     qq2:  that q

     construct case2 qq2 : that r

>>     case2:  [(qq2_1:that q) => (---:that r)]

     close
construct Cases p q r disj , case1 , case2 : that r

>>Cases:  [(p_1:prop),(q_1:prop),(r_1:prop),(disj_1:that 
>>   (p_1 v q_1)),(case1_1:[(pp2_2:that p_1) => (---:that 
>>   r_1)]),(case2_1:[(qq2_3:that q_1) => (---:that r_1)]) 
>>   => (---:that r_1)]

comment The rule of proof by cases really is quite complicated!

\end{verbatim}

We prove the equivalences whcih support the rules mixing disjunction and negation:  these are $P \vee Q \leftrightarrow \neg P \rightarrow Q$ and $P \vee Q \leftrightarrow \neg Q \rightarrow P$.

\begin{verbatim}

comment Prove the basic equivalence theorem 
   which supports mixed rules for disjunction
comment The theorem is (p v q) <-> (~p -> q)
open
     declare aa that p v q

>>     aa:  that (p v q)

     comment our goal is to prove ~p -> q
     open
          declare bb that ~p

>>          bb:  that ~(p)

          comment prove this by cases
          open
               declare hyp1 that p

>>               hyp1:  that p

               declare hyp2 that q

>>               hyp2:  that q

               define casea2 hyp2 :  hyp2

>>               casea2:  [(hyp2_1:that q) => (hyp2_1:that 
>>                  q)]

               open
                    declare cc that ~q

>>                    cc:  that ~(q)

                    define panic cc : Contra p hyp1 bb

>>                    panic:  [(cc_1:that ~(q)) => (Contra(p,
>>                       hyp1,bb):that ??)]

                    close
               define casea1 hyp1 : Dneg q (Negintro ~q panic)

>>               casea1:  [(hyp1_1:that p) => ((q Dneg 
>>                  (~(q) Negintro [(cc_2:that ~(q)) => 
>>                  (Contra(p,hyp1_1,bb):that ??)])):that 
>>                  q)]

               close
          define gotq bb : Cases p q q aa, casea1, casea2

>>          gotq:  [(bb_1:that ~(p)) => (Cases(p,q,q,
>>             aa,[(hyp1_2:that p) => ((q Dneg (~(q) Negintro 
>>             [(cc_3:that ~(q)) => (Contra(p,hyp1_2,bb_1):that 
>>             ??)])):that q)],[(hyp2_4:that q) => (hyp2_4:that 
>>             q)]):that q)]

          close
     define notpimpq aa :  Ifproof ~p q gotq

>>     notpimpq:  [(aa_1:that (p v q)) => (Ifproof(~(p),
>>        q,[(bb_2:that ~(p)) => (Cases(p,q,q,aa_1,[(hyp1_3:that 
>>        p) => ((q Dneg (~(q) Negintro [(cc_4:that ~(q)) 
>>        => (Contra(p,hyp1_3,bb_2):that ??)])):that q)],
>>        [(hyp2_5:that q) => (hyp2_5:that q)]):that q)]):that 
>>        (~(p) -> q))]

     declare bb that ~p -> q

>>     bb:  that (~(p) -> q)

     open 
          declare cc that ~(p v q)

>>          cc:  that ~((p v q))

          comment this is a hypothesis for reduction ad absurdum
          comment our aim is prove ~p so we can use the hypothesis bb
          open 
               declare pp2 that p

>>               pp2:  that p

               define dd pp2 :  Addition1 p q pp2

>>               dd:  [(pp2_1:that p) => (Addition1(p,
>>                  q,pp2_1):that (p v q))]

               define ee pp2 : Contra(p v q, dd pp2 , cc)

>>               ee:  [(pp2_1:that p) => (Contra((p v 
>>                  q),dd(pp2_1),cc):that ??)]

               close
          define ff cc :  Negintro p ee

>>          ff:  [(cc_1:that ~((p v q))) => ((p Negintro 
>>             [(pp2_2:that p) => (Contra((p v q),Addition1(p,
>>             q,pp2_2),cc_1):that ??)]):that ~(p))]

          define gg2 cc :  Mp (~p,q,ff cc,bb)

>>          gg2:  [(cc_1:that ~((p v q))) => (Mp(~(p),
>>             q,ff(cc_1),bb):that q)]

          define hh cc : Addition2 p q gg2 cc

>>          hh:  [(cc_1:that ~((p v q))) => (Addition2(p,
>>             q,gg2(cc_1)):that (p v q))]

          define ii cc : Contra (p v q,hh cc, cc)

>>          ii:  [(cc_1:that ~((p v q))) => (Contra((p 
>>             v q),hh(cc_1),cc_1):that ??)]

          close
     define jj bb : Reductio (p v q,ii)

>>     jj:  [(bb_1:that (~(p) -> q)) => (((p v q) Reductio 
>>        [(cc_2:that ~((p v q))) => (Contra((p v q),Addition2(p,
>>        q,Mp(~(p),q,(p Negintro [(pp2_3:that p) => (Contra((p 
>>        v q),Addition1(p,q,pp2_3),cc_2):that ??)]),bb_1)),
>>        cc_2):that ??)]):that (p v q))]

     close
define Orthm p q : Biintro (p v q, ~p -> q, notpimpq, jj)

>>Orthm:  [(p_1:prop),(q_1:prop) => (Biintro((p_1 v q_1),
>>   (~(p_1) -> q_1),[(aa_2:that (p_1 v q_1)) => (Ifproof(~(p_1),
>>   q_1,[(bb_3:that ~(p_1)) => (Cases(p_1,q_1,q_1,aa_2,
>>   [(hyp1_4:that p_1) => ((q_1 Dneg (~(q_1) Negintro [(cc_5:that 
>>   ~(q_1)) => (Contra(p_1,hyp1_4,bb_3):that ??)])):that 
>>   q_1)],[(hyp2_6:that q_1) => (hyp2_6:that q_1)]):that 
>>   q_1)]):that (~(p_1) -> q_1))],[(bb_7:that (~(p_1) -> 
>>   q_1)) => (((p_1 v q_1) Reductio [(cc_8:that ~((p_1 
>>   v q_1))) => (Contra((p_1 v q_1),Addition2(p_1,q_1,Mp(~(p_1),
>>   q_1,(p_1 Negintro [(pp2_9:that p_1) => (Contra((p_1 
>>   v q_1),Addition1(p_1,q_1,pp2_9),cc_8):that ??)]),bb_7)),
>>   cc_8):that ??)]):that (p_1 v q_1))]):that ((p_1 v q_1) 
>>   <-> (~(p_1) -> q_1)))]

comment Prove the symmetric version p v q <-> ~q -> p
open
     declare aa that p v q

>>     aa:  that (p v q)

     define bb aa : Mpb1 (p v q,~p -> q,aa,Orthm p q)

>>     bb:  [(aa_1:that (p v q)) => (Mpb1((p v q),(~(p) 
>>        -> q),aa_1,(p Orthm q)):that (~(p) -> q))]

     define cc aa : Mpb1 (~p -> q, ~q -> ~ ~ p,bb aa,Contrapositive ~p q)

>>     cc:  [(aa_1:that (p v q)) => (Mpb1((~(p) -> q),
>>        (~(q) -> ~(~(p))),bb(aa_1),(~(p) Contrapositive 
>>        q)):that (~(q) -> ~(~(p))))]

     open
          declare negq that ~q

>>          negq:  that ~(q)

          define dd negq:  Mp ~q ~ ~ p negq cc aa

>>          dd:  [(negq_1:that ~(q)) => (Mp(~(q),~(~(p)),
>>             negq_1,cc(aa)):that ~(~(p)))]

          define yesp negq :  Dneg p dd negq

>>          yesp:  [(negq_1:that ~(q)) => ((p Dneg dd(negq_1)):that 
>>             p)]

          close
     define ee aa :  Ifproof ~q p yesp

>>     ee:  [(aa_1:that (p v q)) => (Ifproof(~(q),p,[(negq_2:that 
>>        ~(q)) => ((p Dneg Mp(~(q),~(~(p)),negq_2,cc(aa_1))):that 
>>        p)]):that (~(q) -> p))]

     declare ff that ~q -> p

>>     ff:  that (~(q) -> p)

     comment Prove that ~p implies q then use Orthm
     open
          declare negp that ~p

>>          negp:  that ~(p)

          comment prove q by reductio
          open 
               declare negq that ~q

>>               negq:  that ~(q)

               define pfollows negq :  Mp ~q p negq ff

>>               pfollows:  [(negq_1:that ~(q)) => (Mp(~(q),
>>                  p,negq_1,ff):that p)]

               define disaster negq :  Contra p, pfollows negq negp

>>               disaster:  [(negq_1:that ~(q)) => (Contra(p,
>>                  pfollows(negq_1),negp):that ??)]

               close
          define kk negp :  Reductio q disaster

>>          kk:  [(negp_1:that ~(p)) => ((q Reductio 
>>             [(negq_2:that ~(q)) => (Contra(p,Mp(~(q),
>>             p,negq_2,ff),negp_1):that ??)]):that q)]

          close
     define ll ff :  Ifproof ~p q kk

>>     ll:  [(ff_1:that (~(q) -> p)) => (Ifproof(~(p),
>>        q,[(negp_2:that ~(p)) => ((q Reductio [(negq_3:that 
>>        ~(q)) => (Contra(p,Mp(~(q),p,negq_3,ff_1),negp_2):that 
>>        ??)]):that q)]):that (~(p) -> q))]

     define mm ff :  Mpb2 (p v q,~p -> q,ll ff,Orthm p q)

>>     mm:  [(ff_1:that (~(q) -> p)) => (Mpb2((p v q),
>>        (~(p) -> q),ll(ff_1),(p Orthm q)):that (p v q))]

     close
define Orthm2 p q :  Biintro (p v q, ~q -> p, ee, mm)

>>Orthm2:  [(p_1:prop),(q_1:prop) => (Biintro((p_1 v 
>>   q_1),(~(q_1) -> p_1),[(aa_2:that (p_1 v q_1)) => (Ifproof(~(q_1),
>>   p_1,[(negq_3:that ~(q_1)) => ((p_1 Dneg Mp(~(q_1),~(~(p_1)),
>>   negq_3,Mpb1((~(p_1) -> q_1),(~(q_1) -> ~(~(p_1))),Mpb1((p_1 
>>   v q_1),(~(p_1) -> q_1),aa_2,(p_1 Orthm q_1)),(~(p_1) 
>>   Contrapositive q_1)))):that p_1)]):that (~(q_1) -> 
>>   p_1))],[(ff_4:that (~(q_1) -> p_1)) => (Mpb2((p_1 v 
>>   q_1),(~(p_1) -> q_1),Ifproof(~(p_1),q_1,[(negp_5:that 
>>   ~(p_1)) => ((q_1 Reductio [(negq_6:that ~(q_1)) => 
>>   (Contra(p_1,Mp(~(q_1),p_1,negq_6,ff_4),negp_5):that 
>>   ??)]):that q_1)]),(p_1 Orthm q_1)):that (p_1 v q_1))]):that 
>>   ((p_1 v q_1) <-> (~(q_1) -> p_1)))]

\end{verbatim}

We derive stronger disjunction introduction rules and the rules of disjunctive syllogism.

\begin{verbatim}

comment Develop the full dress disjunction introduction rule
comment reversal of numbering is due to proving the less preferred
open
     declare negq that ~q

>>     negq:  that ~(q)

     construct thusp negq : that p

>>     thusp:  [(negq_1:that ~(q)) => (---:that p)]

     close

define Disjintro p q thusp: 
Mpb2 (p v q, ~q -> p, Ifproof ~q p thusp, Orthm2 p q)

>>Disjintro:  [(p_1:prop),(q_1:prop),(thusp_1:[(negq_2:that 
>>   ~(q_1)) => (---:that p_1)]) => (Mpb2((p_1 v q_1),(~(q_1) 
>>   -> p_1),Ifproof(~(q_1),p_1,thusp_1),(p_1 Orthm2 q_1)):that 
>>   (p_1 v q_1))]

open
     declare negp that ~p

>>     negp:  that ~(p)

     construct thusq negp : that q

>>     thusq:  [(negp_1:that ~(p)) => (---:that q)]

     close
define Disjintro2 p q thusq: 
Mpb2 (p v q, ~p -> q, Ifproof ~p q thusq, Orthm p q)

>>Disjintro2:  [(p_1:prop),(q_1:prop),(thusq_1:[(negp_2:that 
>>   ~(p_1)) => (---:that q_1)]) => (Mpb2((p_1 v q_1),(~(p_1) 
>>   -> q_1),Ifproof(~(p_1),q_1,thusq_1),(p_1 Orthm q_1)):that 
>>   (p_1 v q_1))]

comment Rules of disjunctive syllogism
declare line1 that p v q

>>line1:  that (p v q)

declare line2 that ~q

>>line2:  that ~(q)

define Ds1 p q line1 line2 :  
Mp (~q, p, line2, Mpb1 (p v q, ~q -> p, line1, Orthm2 p q)) 

>>Ds1:  [(p_1:prop),(q_1:prop),(line1_1:that (p_1 v q_1)),
>>   (line2_1:that ~(q_1)) => (Mp(~(q_1),p_1,line2_1,Mpb1((p_1 
>>   v q_1),(~(q_1) -> p_1),line1_1,(p_1 Orthm2 q_1))):that 
>>   p_1)]

declare line3 that p v q

>>line3:  that (p v q)

declare line4 that ~p

>>line4:  that ~(p)

define Ds2 p q line3 line4 :  
Mp (~p, q, line4, Mpb1 (p v q, ~p -> q, line3, Orthm p q)) 

>>Ds2:  [(p_1:prop),(q_1:prop),(line3_1:that (p_1 v q_1)),
>>   (line4_1:that ~(p_1)) => (Mp(~(p_1),q_1,line4_1,Mpb1((p_1 
>>   v q_1),(~(p_1) -> q_1),line3_1,(p_1 Orthm q_1))):that 
>>   q_1)]

\end{verbatim}

\subsection{The existential quantifier and a quantifier proof}

In this section we introduce the existential quantifier and its primitive rules, then prove the quantifier theorem $(\forall x:P(x) \rightarrow Q(x)) \wedge (\forall x:Q(x) \rightarrow R(x)) \rightarrow (\forall x:P(x) \rightarrow R(x))$.

\begin{verbatim}

comment The existential quantifier
construct Exists P : prop

>>Exists:  [(P_1:[(xx_2:obj) => (---:prop)]) => (---:prop)]

\end{verbatim}

The existential quantifier introduction rule (EG).

\begin{verbatim}

comment the rule EG (existential introduction)
declare ev that P x

>>ev:  that P(x)

construct Eg P, x ev :  that Exists P

>>Eg:  [(P_1:[(xx_2:obj) => (---:prop)]),(x_1:obj),(ev_1:that 
>>   P_1(x_1)) => (---:that Exists(P_1))]

\end{verbatim}

The existential quantifier elimination rule (EI).  This is rather complicated!

\begin{verbatim}

comment the rule EI (existential elimination)
declare g prop

>>g:  prop

declare ex that Exists P

>>ex:  that Exists(P)

open
     declare w obj

>>     w:  obj

     declare ev2 that P w

>>     ev2:  that P(w)

     construct wi w ev2 :  that g

>>     wi:  [(w_1:obj),(ev2_1:that P(w_1)) => (---:that 
>>        g)]

     close
construct Ei P, g, ex, wi :  that g

>>Ei:  [(P_1:[(xx_2:obj) => (---:prop)]),(g_1:prop),(ex_1:that 
>>   Exists(P_1)),(wi_1:[(w_3:obj),(ev2_3:that P_1(w_3)) 
>>   => (---:that g_1)]) => (---:that g_1)]

\end{verbatim}

The proof of $(\forall x:P(x) \rightarrow Q(x)) \wedge (\forall x:Q(x) \rightarrow R(x)) \rightarrow (\forall x:P(x) \rightarrow R(x))$.  Notice that while we never write
quantified statements with  variable binding constructions explicit, they do actually appear in the display of the final result, because the identifiers used to specify the component abstractions  pass out of scope.

\begin{verbatim}

comment A quantifier proof
open
     declare xx obj

>>     xx:  obj

     construct Pp xx :prop

>>     Pp:  [(xx_1:obj) => (---:prop)]

     construct Qq xx : prop

>>     Qq:  [(xx_1:obj) => (---:prop)]

     construct Rr xx:prop

>>     Rr:  [(xx_1:obj) => (---:prop)]

     define Ss xx: (Pp xx) -> (Qq xx)

>>     Ss:  [(xx_1:obj) => ((Pp(xx_1) -> Qq(xx_1)):prop)]

     define Tt xx: (Qq xx) -> (Rr xx)

>>     Tt:  [(xx_1:obj) => ((Qq(xx_1) -> Rr(xx_1)):prop)]

     define Uu xx:  (Pp xx) -> (Rr xx)

>>     Uu:  [(xx_1:obj) => ((Pp(xx_1) -> Rr(xx_1)):prop)]

     declare ss2   that Forall Ss

>>     ss2:  that Forall(Ss)

     declare tt2   that Forall Tt

>>     tt2:  that Forall(Tt)

     comment Our goal is to prove Forall Uu
     open
          declare yy obj

>>          yy:  obj

          comment Our goal is to show (Pp yy) -> (Rr yy)
          open
               declare ppyy that Pp yy

>>               ppyy:  that Pp(yy)

               define imp1 :  Ui Ss, ss2 yy

>>               imp1:  [(Ui(Ss,ss2,yy):that Ss(yy))]

               define line5 ppyy: Mp (Pp yy, Qq yy, ppyy, imp1)

>>               line5:  [(ppyy_1:that Pp(yy)) => (Mp(Pp(yy),
>>                  Qq(yy),ppyy_1,imp1):that Qq(yy))]

               define imp2 : Ui Tt, tt2 yy

>>               imp2:  [(Ui(Tt,tt2,yy):that Tt(yy))]

               define line6 ppyy: Mp (Qq yy, Rr yy,line5 ppyy,imp2)

>>               line6:  [(ppyy_1:that Pp(yy)) => (Mp(Qq(yy),
>>                  Rr(yy),line5(ppyy_1),imp2):that Rr(yy))]

               close
          define line7 yy: Ifproof (Pp yy, Rr yy,line6)

>>          line7:  [(yy_1:obj) => (Ifproof(Pp(yy_1),
>>             Rr(yy_1),[(ppyy_2:that Pp(yy_1)) => (Mp(Qq(yy_1),
>>             Rr(yy_1),Mp(Pp(yy_1),Qq(yy_1),ppyy_2,Ui(Ss,
>>             ss2,yy_1)),Ui(Tt,tt2,yy_1)):that Rr(yy_1))]):that 
>>             (Pp(yy_1) -> Rr(yy_1)))]

          close
     define Univimp1 ss2 tt2: Ug Uu, line7

>>     Univimp1:  [(ss2_1:that Forall(Ss)),(tt2_1:that 
>>        Forall(Tt)) => (Ug(Uu,[(yy_2:obj) => (Ifproof(Pp(yy_2),
>>        Rr(yy_2),[(ppyy_3:that Pp(yy_2)) => (Mp(Qq(yy_2),
>>        Rr(yy_2),Mp(Pp(yy_2),Qq(yy_2),ppyy_3,Ui(Ss,ss2_1,
>>        yy_2)),Ui(Tt,tt2_1,yy_2)):that Rr(yy_2))]):that 
>>        (Pp(yy_2) -> Rr(yy_2)))]):that Forall(Uu))]

     declare conj1 that (Forall Ss) & (Forall Tt)

>>     conj1:  that (Forall(Ss) & Forall(Tt))

     define Univimp2 conj1 : 
Univimp1 (And1(Forall Ss, Forall Tt,conj1),And2(Forall Ss, Forall Tt,conj1))

>>     Univimp2:  [(conj1_1:that (Forall(Ss) & Forall(Tt))) 
>>        => ((And1(Forall(Ss),Forall(Tt),conj1_1) Univimp1 
>>        And2(Forall(Ss),Forall(Tt),conj1_1)):that Forall(Uu))]

     close
define Univimp Pp, Qq, Rr : 
Ifproof ((Forall Ss)&(Forall Tt),Forall Uu,Univimp2)

>>Univimp:  [(Pp_1:[(xx_2:obj) => (---:prop)]),(Qq_1:[(xx_3:obj) 
>>   => (---:prop)]),(Rr_1:[(xx_4:obj) => (---:prop)]) => 
>>   (Ifproof((Forall([(xx_5:obj) => ((Pp_1(xx_5) -> Qq_1(xx_5)):prop)]) 
>>   & Forall([(xx_6:obj) => ((Qq_1(xx_6) -> Rr_1(xx_6)):prop)])),
>>   Forall([(xx_7:obj) => ((Pp_1(xx_7) -> Rr_1(xx_7)):prop)]),
>>   [(conj1_8:that (Forall([(xx_9:obj) => ((Pp_1(xx_9) 
>>   -> Qq_1(xx_9)):prop)]) & Forall([(xx_10:obj) => ((Qq_1(xx_10) 
>>   -> Rr_1(xx_10)):prop)]))) => (Ug([(xx_11:obj) => ((Pp_1(xx_11) 
>>   -> Rr_1(xx_11)):prop)],[(yy_12:obj) => (Ifproof(Pp_1(yy_12),
>>   Rr_1(yy_12),[(ppyy_13:that Pp_1(yy_12)) => (Mp(Qq_1(yy_12),
>>   Rr_1(yy_12),Mp(Pp_1(yy_12),Qq_1(yy_12),ppyy_13,Ui([(xx_14:obj) 
>>   => ((Pp_1(xx_14) -> Qq_1(xx_14)):prop)],And1(Forall([(xx_15:obj) 
>>   => ((Pp_1(xx_15) -> Qq_1(xx_15)):prop)]),Forall([(xx_16:obj) 
>>   => ((Qq_1(xx_16) -> Rr_1(xx_16)):prop)]),conj1_8),yy_12)),
>>   Ui([(xx_17:obj) => ((Qq_1(xx_17) -> Rr_1(xx_17)):prop)],
>>   And2(Forall([(xx_18:obj) => ((Pp_1(xx_18) -> Qq_1(xx_18)):prop)]),
>>   Forall([(xx_19:obj) => ((Qq_1(xx_19) -> Rr_1(xx_19)):prop)]),
>>   conj1_8),yy_12)):that Rr_1(yy_12))]):that (Pp_1(yy_12) 
>>   -> Rr_1(yy_12)))]):that Forall([(xx_20:obj) => ((Pp_1(xx_20) 
>>   -> Rr_1(xx_20)):prop)]))]):that ((Forall([(xx_21:obj) 
>>   => ((Pp_1(xx_21) -> Qq_1(xx_21)):prop)]) & Forall([(xx_22:obj) 
>>   => ((Qq_1(xx_22) -> Rr_1(xx_22)):prop)])) -> Forall([(xx_23:obj) 
>>   => ((Pp_1(xx_23) -> Rr_1(xx_23)):prop)])))]

\end{verbatim}

\subsection{Sample declarations of theories of typed objects:  natural numbers and the simple theory of types}

Lestrade is not exclusively deveoted to constructing proof objects.  We give some very compact definitions for typed objects -- natural numbers and the sets of a model of the simple typed theory of sets.

\begin{verbatim}

comment  Declarations of typed objects
comment The type of (true) natural numbers.  The theory of these
comment objects will be second order arithmetic.  Peano arithmetic
comment will be defined:  it will be instructive how hard it is to do this.

\end{verbatim}

The natural numbers as a type, the successor operation, and 1 are declared.

\begin{verbatim}

construct Nat : type

>>Nat:  [(---:type)]

construct 1 : in Nat

>>1:  [(---:in Nat)]

declare n in Nat

>>n:  in Nat

construct Succ n : in Nat

>>Succ:  [(n_1:in Nat) => (---:in Nat)]

\end{verbatim}

The declaration of mathematical induction.   A serious development would include quantifiers over the natural numbers.

\begin{verbatim}

open
     declare n2 in Nat

>>     n2:  in Nat

     construct Pn n2 : prop

>>     Pn:  [(n2_1:in Nat) => (---:prop)]

     close
declare basis that Pn 1

>>basis:  that Pn(1)

open
     declare k in Nat

>>     k:  in Nat

     declare indhyp that Pn k

>>     indhyp:  that Pn(k)

     construct indstep k indhyp :  that Pn Succ k

>>     indstep:  [(k_1:in Nat),(indhyp_1:that Pn(k_1)) 
>>        => (---:that Pn(Succ(k_1)))]

     close
construct Induction n Pn, basis, indstep :  that Pn n

>>Induction:  [(n_1:in Nat),(Pn_1:[(n2_2:in Nat) => (---:prop)]),
>>   (basis_1:that Pn_1(1)),(indstep_1:[(k_3:in Nat),(indhyp_3:that 
>>   Pn_1(k_3)) => (---:that Pn_1(Succ(k_3)))]) => (---:that 
>>   Pn_1(n_1))]

\end{verbatim}

Equality and its rules are declared.

\begin{verbatim}

comment We introduce the declarations for the properties
comment of equality of natural numbers.
declare m in Nat

>>m:  in Nat

declare m2 in Nat

>>m2:  in Nat

construct Eqn n m : prop

>>Eqn:  [(n_1:in Nat),(m_1:in Nat) => (---:prop)]


\end{verbatim}

Equality elimination (the rule of substitution)

\begin{verbatim}

comment We develop the substitution rule (equality elimination)
declare eqev that Eqn m m2

>>eqev:  that (m Eqn m2)

declare pnpf that Pn m

>>pnpf:  that Pn(m)

construct Subs Pn, m m2 eqev pnpf: that Pn m2

>>Subs:  [(Pn_1:[(n2_2:in Nat) => (---:prop)]),(m_1:in 
>>   Nat),(m2_1:in Nat),(eqev_1:that (m_1 Eqn m2_1)),(pnpf_1:that 
>>   Pn_1(m_1)) => (---:that Pn_1(m2_1))]

\end{verbatim}

Equality introduction (indiscerniblity).

\begin{verbatim}

comment We develop the equality introduction rule (Leibniz)
open
     open
          declare n3 in Nat

>>          n3:  in Nat

          construct Pn2 n3:  prop

>>          Pn2:  [(n3_1:in Nat) => (---:prop)]

          close
     declare pnn that Pn2 n

>>     pnn:  that Pn2(n)

     construct eqpf Pn2, pnn:  that Pn2 m

>>     eqpf:  [(Pn2_1:[(n3_2:in Nat) => (---:prop)]),
>>        (pnn_1:that Pn2_1(n)) => (---:that Pn2_1(m))]

     close
construct Eqnproof n m, eqpf :  that n Eqn m

>>Eqnproof:  [(n_1:in Nat),(m_1:in Nat),(eqpf_1:[(Pn2_2:[(n3_3:in 
>>   Nat) => (---:prop)]),(pnn_2:that Pn2_2(n_1)) => (---:that 
>>   Pn2_2(m_1))]) => (---:that (n_1 Eqn m_1))]

\end{verbatim}

We prove the trivial theorem of reflexivity of equality.

\begin{verbatim}

comment We test the equality introduction rule
comment by proving reflexivity of equality.
open
     open
          declare n3 in Nat

>>          n3:  in Nat

          construct Pn2 n3:prop

>>          Pn2:  [(n3_1:in Nat) => (---:prop)]

          close
     declare pnn that Pn2 n

>>     pnn:  that Pn2(n)

     define eqpftest Pn2, pnn: pnn

>>     eqpftest:  [(Pn2_1:[(n3_2:in Nat) => (---:prop)]),
>>        (pnn_1:that Pn2_1(n)) => (pnn_1:that Pn2_1(n))]

     close
define Refln n : Eqnproof n n, eqpftest

>>Refln:  [(n_1:in Nat) => (Eqnproof(n_1,n_1,[(Pn2_2:[(n3_3:in 
>>   Nat) => (---:prop)]),(pnn_2:that Pn2_2(n_1)) => (pnn_2:that 
>>   Pn2_2(n_1))]):that (n_1 Eqn n_1))]

\end{verbatim}

We had to declare equality in order to declare the other two Peano axioms:  here they are.

\begin{verbatim}

construct Pa3 n :  that ~(Succ n Eqn 1)

>>Pa3:  [(n_1:in Nat) => (---:that ~((Succ(n_1) Eqn 1)))]

construct Pa4 n m :  that (Succ n Eqn Succ m) -> n Eqn m

>>Pa4:  [(n_1:in Nat),(m_1:in Nat) => (---:that ((Succ(n_1) 
>>   Eqn Succ(m_1)) -> (n_1 Eqn m_1)))]

comment These definitions are by no means exhaustive.  One wants
comment to declare quantifiers over natural numbers for example.

\end{verbatim}

Here are a set of very economical declarations for the simple type theory of sets.  Once again, I have not declare the quantifiers that are surely wanted.

\begin{verbatim}


comment Declarations for second order type theory.

\end{verbatim}

I could declare the types without using natural numbers at all, but since I have them I will use them.

\begin{verbatim}

construct level n : type

>>level:  [(n_1:in Nat) => (---:type)]

comment level n is what we usually call type n.  
   The bottom type will be type 1.
declare n3 in Nat

>>n3:  in Nat

declare x10 in level n3

>>x10:  in level(n3)

declare y10 in level Succ n3

>>y10:  in level(Succ(n3))

\end{verbatim}

Here is the membership relation.  It is ternary:  we must unavoidably put a natural number first to indicate the type of the element.

\begin{verbatim}

comment Declare the membership relation (with a type argument)
construct E n3 x10 y10 : prop

>>E:  [(n3_1:in Nat),(x10_1:in level(n3_1)),(y10_1:in 
>>   level(Succ(n3_1))) => (---:prop)]

\end{verbatim}

Here is the set abstract primitive, taking predicates of type $n$ objects to sets of type $n+1$.

\begin{verbatim}

comment Declare the set abstract constructor
open
     declare x11 in level n3

>>     x11:  in level(n3)

     construct Pt x11 : prop

>>     Pt:  [(x11_1:in level(n3)) => (---:prop)]

     close

\end{verbatim}

Here are the comprehension axioms, declared in a most economical way.

\begin{verbatim}

comment Declare the comprehension axioms
construct setof n3 Pt :  in level Succ n3

>>setof:  [(n3_1:in Nat),(Pt_1:[(x11_2:in level(n3_1)) 
>>   => (---:prop)]) => (---:in level(Succ(n3_1)))]

declare compev1 that E(n3,x10,setof n3 Pt)

>>compev1:  that E(n3,x10,(n3 setof Pt))

construct Comp1 n3 x10, Pt :  that Pt x10

>>Comp1:  [(n3_1:in Nat),(x10_1:in level(n3_1)),(Pt_1:[(x11_2:in 
>>   level(n3_1)) => (---:prop)]) => (---:that Pt_1(x10_1))]

declare compev2 that Pt x10

>>compev2:  that Pt(x10)

construct Comp2 n3 x10, Pt :  that E(n3,x10,setof n3 Pt)

>>Comp2:  [(n3_1:in Nat),(x10_1:in level(n3_1)),(Pt_1:[(x11_2:in 
>>   level(n3_1)) => (---:prop)]) => (---:that E(n3_1,x10_1,
>>   (n3_1 setof Pt_1)))]

\end{verbatim}

Here is the extensionality axiom, whose force is that things having the same elements (and belonging to a successor type) themselves belong to the same sets.
In the presence of comprehension, this is enough:  objects with the same elements will thus be indiscernible.  Of course a definition of equality (and definitions of quantifiers)
would appear in a full treatment.

\begin{verbatim}

comment Declare the extensionality axiom
declare xx10 in level Succ n3

>>xx10:  in level(Succ(n3))

declare yy10 in level Succ n3

>>yy10:  in level(Succ(n3))

declare ww10 in level Succ(Succ n3)

>>ww10:  in level(Succ(Succ(n3)))

declare xinw that (Succ n3) E xx10 ww10

>>xinw:  that E(Succ(n3),xx10,ww10)

open
     declare z11 in level n3

>>     z11:  in level(n3)

     declare zinx that n3 E z11 xx10

>>     zinx:  that E(n3,z11,xx10)

     declare ziny that n3 E z11 yy10

>>     ziny:  that E(n3,z11,yy10)

     construct xincy z11 zinx :  that n3 E z11 yy10

>>     xincy:  [(z11_1:in level(n3)),(zinx_1:that E(n3,
>>        z11_1,xx10)) => (---:that E(n3,z11_1,yy10))]

     construct yincx z11 ziny :  that n3 E z11 xx10

>>     yincx:  [(z11_1:in level(n3)),(ziny_1:that E(n3,
>>        z11_1,yy10)) => (---:that E(n3,z11_1,xx10))]

     close

construct Extensionality n3 xx10 yy10 ww10,xinw, xincy, yincx :  
that (Succ n3) E yy10 ww10

>>Extensionality:  [(n3_1:in Nat),(xx10_1:in level(Succ(n3_1))),
>>   (yy10_1:in level(Succ(n3_1))),(ww10_1:in level(Succ(Succ(n3_1)))),
>>   (xinw_1:that E(Succ(n3_1),xx10_1,ww10_1)),(xincy_1:[(z11_2:in 
>>   level(n3_1)),(zinx_2:that E(n3_1,z11_2,xx10_1)) => 
>>   (---:that E(n3_1,z11_2,yy10_1))]),(yincx_1:[(z11_3:in 
>>   level(n3_1)),(ziny_3:that E(n3_1,z11_3,yy10_1)) => 
>>   (---:that E(n3_1,z11_3,xx10_1))]) => (---:that E(Succ(n3_1),
>>   yy10_1,ww10_1))]

\end{verbatim}

\subsection{Quantifiers over general types}

In this section we present the primitives supporting quantification over all types of sort {\tt type}.  Quantifiers over the sort {\tt type} itself could be introduced independently
but would of course be more dangerous (it would be easier to utter paradoxes).  These tools could be used to define quantifiers over the natural numbers and over each of the levels
of type theory without any need for new primitives.

This is an illustration of the fact that we have a general ability to construct and define operations on a domain of types.

The text is cloned from the text for the quantifiers over the sort {\tt obj} above, and the comments have mostly not been revised to reflect the changes
in identifiers.

\begin{verbatim}

comment Declaring the universal quantifier for general types.
declare tau type

>>tau:  type

open
     declare uu in tau

>>     uu:  in tau

     construct Ptt uu : prop

>>     Ptt:  [(uu_1:in tau) => (---:prop)]

     close
construct Forallt tau Ptt: prop

>>Forallt:  [(tau_1:type),(Ptt_1:[(uu_2:in tau_1) => 
>>   (---:prop)]) => (---:prop)]

comment Declaring the rule UI of universal instantiation (for general types)
declare Ptt2 that Forallt tau Ptt

>>Ptt2:  that (tau Forallt Ptt)

declare xt in tau

>>xt:  in tau

construct Uit tau Ptt, Ptt2 xt : that Ptt xt

>>Uit:  [(tau_1:type),(Ptt_1:[(uu_2:in tau_1) => (---:prop)]),
>>   (Ptt2_1:that (tau_1 Forallt Ptt_1)),(xt_1:in tau_1) 
>>   => (---:that Ptt_1(xt_1))]

comment Note in the previous line that we follow P 
comment with a comma:  an abstraction argument may need to be 
comment guarded with commas so it will not be read as applied.
comment Opening an environment to declare a function 
comment that witnesses provability of a universal statement
open
     declare ut in tau

>>     ut:  in tau

     construct Qt2 ut : that Ptt ut

>>     Qt2:  [(ut_1:in tau) => (---:that Ptt(ut_1))]

     close
comment The rule of universal generalization (for general types)
construct Ugt tau Ptt, Qt2 : that Forallt tau Ptt

>>Ugt:  [(tau_1:type),(Ptt_1:[(uu_2:in tau_1) => (---:prop)]),
>>   (Qt2_1:[(ut_3:in tau_1) => (---:that Ptt_1(ut_3))]) 
>>   => (---:that (tau_1 Forallt Ptt_1))]

comment The existential quantifier (for general types)
construct Existst tau Ptt : prop

>>Existst:  [(tau_1:type),(Ptt_1:[(uu_2:in tau_1) => 
>>   (---:prop)]) => (---:prop)]

comment the rule EG (existential introduction) (for general types)
declare evt that Ptt xt

>>evt:  that Ptt(xt)

construct Egt tau Ptt, xt evt :  that Existst tau Ptt

>>Egt:  [(tau_1:type),(Ptt_1:[(uu_2:in tau_1) => (---:prop)]),
>>   (xt_1:in tau_1),(evt_1:that Ptt_1(xt_1)) => (---:that 
>>   (tau_1 Existst Ptt_1))]

comment the rule EI (existential elimination) (for general types)
declare gt prop

>>gt:  prop

declare ext that Existst tau Ptt

>>ext:  that (tau Existst Ptt)

open
     declare wt in tau

>>     wt:  in tau

     declare evt2 that Ptt wt

>>     evt2:  that Ptt(wt)

     construct wit wt evt2 :  that gt

>>     wit:  [(wt_1:in tau),(evt2_1:that Ptt(wt_1)) => 
>>        (---:that gt)]

     close
construct Eit tau Ptt, gt, ext, wit :  that gt

>>Eit:  [(tau_1:type),(Ptt_1:[(uu_2:in tau_1) => (---:prop)]),
>>   (gt_1:prop),(ext_1:that (tau_1 Existst Ptt_1)),(wit_1:[(wt_3:in 
>>   tau_1),(evt2_3:that Ptt_1(wt_3)) => (---:that gt_1)]) 
>>   => (---:that gt_1)]

\end{verbatim}

\subsection{Equality and the definite description operator for untyped and typed objects}

\begin{verbatim}

comment  Equality uniqueness and definite description
declare y obj

>>y:  obj

comment Equality of untyped objects
construct = x y : prop

>>=:  [(x_1:obj),(y_1:obj) => (---:prop)]

comment Develop equality introduction rule (indiscernibility)
open
     open 
          declare x2 obj

>>          x2:  obj

          construct Peq2 x2: prop

>>          Peq2:  [(x2_1:obj) => (---:prop)]

          close
     declare pxev that Peq2 x

>>     pxev:  that Peq2(x)

     construct pyev Peq2, pxev : that Peq2 y

>>     pyev:  [(Peq2_1:[(x2_2:obj) => (---:prop)]),(pxev_1:that 
>>        Peq2_1(x)) => (---:that Peq2_1(y))]

     close
construct Eqintro x y pyev :that x = y

>>Eqintro:  [(x_1:obj),(y_1:obj),(pyev_1:[(Peq2_2:[(x2_3:obj) 
>>   => (---:prop)]),(pxev_2:that Peq2_2(x_1)) => (---:that 
>>   Peq2_2(y_1))]) => (---:that (x_1 = y_1))]

comment Construct equality elimination rule (substitution)
declare xyeqev that x = y

>>xyeqev:  that (x = y)

declare pxev that P x

>>pxev:  that P(x)

construct Eqelim P, x y xyeqev pxev :  that P y

>>Eqelim:  [(P_1:[(xx_2:obj) => (---:prop)]),(x_1:obj),
>>   (y_1:obj),(xyeqev_1:that (x_1 = y_1)),(pxev_1:that 
>>   P_1(x_1)) => (---:that P_1(y_1))]

comment The same rules for equality, adapted to general types
declare yt in tau

>>yt:  in tau

construct eqt tau xt yt : prop

>>eqt:  [(tau_1:type),(xt_1:in tau_1),(yt_1:in tau_1) 
>>   => (---:prop)]

comment Develop equality introduction rule (indiscernibility)
open
     open 
          declare x2 in tau

>>          x2:  in tau

          construct Peqt2 x2: prop

>>          Peqt2:  [(x2_1:in tau) => (---:prop)]

          close
     declare pxevt that Peqt2 xt

>>     pxevt:  that Peqt2(xt)

     construct pyevt Peqt2, pxevt : that Peqt2 yt

>>     pyevt:  [(Peqt2_1:[(x2_2:in tau) => (---:prop)]),
>>        (pxevt_1:that Peqt2_1(xt)) => (---:that Peqt2_1(yt))]

     close
construct Eqintrot tau xt yt pyevt :that tau eqt xt yt

>>Eqintrot:  [(tau_1:type),(xt_1:in tau_1),(yt_1:in tau_1),
>>   (pyevt_1:[(Peqt2_2:[(x2_3:in tau_1) => (---:prop)]),
>>   (pxevt_2:that Peqt2_2(xt_1)) => (---:that Peqt2_2(yt_1))]) 
>>   => (---:that eqt(tau_1,xt_1,yt_1))]

comment Construct equality elimination rule (substitution)
declare xyeqevt that tau eqt xt yt

>>xyeqevt:  that eqt(tau,xt,yt)

declare pxevt that Ptt xt

>>pxevt:  that Ptt(xt)

construct Eqelimt tau Ptt, xt yt xyeqevt pxevt :  that Ptt yt

>>Eqelimt:  [(tau_1:type),(Ptt_1:[(uu_2:in tau_1) => 
>>   (---:prop)]),(xt_1:in tau_1),(yt_1:in tau_1),(xyeqevt_1:that 
>>   eqt(tau_1,xt_1,yt_1)),(pxevt_1:that Ptt_1(xt_1)) => 
>>   (---:that Ptt_1(yt_1))]

comment The definite description operator
declare atleast1 that Exists P

>>atleast1:  that Exists(P)

open
     declare x1 obj

>>     x1:  obj

     declare x2 obj

>>     x2:  obj

     declare thatpx1 that P x1

>>     thatpx1:  that P(x1)

     declare thatpx2 that P x2

>>     thatpx2:  that P(x2)

     construct atmost1 x1 x2 thatpx1 thatpx2 : that x1 = x2

>>     atmost1:  [(x1_1:obj),(x2_1:obj),(thatpx1_1:that 
>>        P(x1_1)),(thatpx2_1:that P(x2_1)) => (---:that 
>>        (x1_1 = x2_1))]

     close
construct The P, atleast1 atmost1 : obj

>>The:  [(P_1:[(xx_2:obj) => (---:prop)]),(atleast1_1:that 
>>   Exists(P_1)),(atmost1_1:[(x1_3:obj),(x2_3:obj),(thatpx1_3:that 
>>   P_1(x1_3)),(thatpx2_3:that P_1(x2_3)) => (---:that 
>>   (x1_3 = x2_3))]) => (---:obj)]

construct Theprop P, atleast1 atmost1 :  that P The P, atleast1 atmost1

>>Theprop:  [(P_1:[(xx_2:obj) => (---:prop)]),(atleast1_1:that 
>>   Exists(P_1)),(atmost1_1:[(x1_3:obj),(x2_3:obj),(thatpx1_3:that 
>>   P_1(x1_3)),(thatpx2_3:that P_1(x2_3)) => (---:that 
>>   (x1_3 = x2_3))]) => (---:that P_1(The(P_1,atleast1_1,
>>   atmost1_1)))]

comment The definite description operator (for general types)
declare atleastt1 that Existst tau Ptt

>>atleastt1:  that (tau Existst Ptt)

open
     declare x1 in tau

>>     x1:  in tau

     declare x2 in tau

>>     x2:  in tau

     declare thatpx1 that Ptt x1

>>     thatpx1:  that Ptt(x1)

     declare thatpx2 that Ptt x2

>>     thatpx2:  that Ptt(x2)

     construct atmostt1 x1 x2 thatpx1 thatpx2 : that tau eqt x1 x2

>>     atmostt1:  [(x1_1:in tau),(x2_1:in tau),(thatpx1_1:that 
>>        Ptt(x1_1)),(thatpx2_1:that Ptt(x2_1)) => (---:that 
>>        eqt(tau,x1_1,x2_1))]

     close
construct Thet tau Ptt, atleastt1 atmostt1 : in tau

>>Thet:  [(tau_1:type),(Ptt_1:[(uu_2:in tau_1) => (---:prop)]),
>>   (atleastt1_1:that (tau_1 Existst Ptt_1)),(atmostt1_1:[(x1_3:in 
>>   tau_1),(x2_3:in tau_1),(thatpx1_3:that Ptt_1(x1_3)),
>>   (thatpx2_3:that Ptt_1(x2_3)) => (---:that eqt(tau_1,
>>   x1_3,x2_3))]) => (---:in tau_1)]

construct Thepropt tau Ptt, atleastt1 atmostt1 : 
 that Ptt Thet tau Ptt, atleastt1 atmostt1

>>Thepropt:  [(tau_1:type),(Ptt_1:[(uu_2:in tau_1) => 
>>   (---:prop)]),(atleastt1_1:that (tau_1 Existst Ptt_1)),
>>   (atmostt1_1:[(x1_3:in tau_1),(x2_3:in tau_1),(thatpx1_3:that 
>>   Ptt_1(x1_3)),(thatpx2_3:that Ptt_1(x2_3)) => (---:that 
>>   eqt(tau_1,x1_3,x2_3))]) => (---:that Ptt_1(Thet(tau_1,
>>   Ptt_1,atleastt1_1,atmostt1_1)))]

\end{verbatim}

\subsection{Declarations for complex type theories, up to the level of bootstrapping Lestrade's own abstraction sorts}

The declarations here will support reasoning in Church's simple type theory  of \cite{churchtypes} (or modern variations); they are further augmented with dependent product and function types of a sort which would be required to emulate Lestrade's own system of abstraction sorts.

\begin{verbatim}

% Church's type theory

% one point type

construct One type

>> One:  [(---:type)] {move 0}

construct Unique : in One

>> Unique:  [(---:in One)] {move 0}

declare xx1 in One

>> xx1:  in One {move 1}

construct Oneproof xx1 :  that One eqt xx1 Unique

>> Oneproof:  [(xx1_1:in One) => (---:that eqt(One,xx1_1,
>>   Unique))] {move 0}

% cartesian product construction

declare sigma type

>> sigma:  type {move 1}

construct X tau sigma : type

>> X:  [(tau_1:type),(sigma_1:type) => (---:type)] {move 
>>   0}

declare xt in tau

>> xt:  in tau {move 1}

declare ys in sigma

>> ys:  in sigma {move 1}

construct pair tau sigma xt ys : in tau X sigma

>> pair:  [(tau_1:type),(sigma_1:type),(xt_1:in tau_1),
>>   (ys_1:in sigma_1) => (---:in (tau_1 X sigma_1))] {move 
>>   0}

declare zp in tau X sigma

>> zp:  in (tau X sigma) {move 1}

construct pi1 tau sigma zp :  in tau

>> pi1:  [(tau_1:type),(sigma_1:type),(zp_1:in (tau_1 
>>   X sigma_1)) => (---:in tau_1)] {move 0}

construct pi2 tau sigma zp : in sigma

>> pi2:  [(tau_1:type),(sigma_1:type),(zp_1:in (tau_1 
>>   X sigma_1)) => (---:in sigma_1)] {move 0}

construct Xexact tau sigma zp : 
   that (tau X sigma) eqt zp, 
   pair tau sigma (pi1 tau sigma zp) (pi2 tau sigma zp)

>> Xexact:  [(tau_1:type),(sigma_1:type),(zp_1:in (tau_1 
>>   X sigma_1)) => (---:that eqt((tau_1 X sigma_1),zp_1,
>>   pair(tau_1,sigma_1,pi1(tau_1,sigma_1,zp_1),pi2(tau_1,
>>   sigma_1,zp_1))))] {move 0}

% power set type constructor (use this to build bool from one point type)

construct Pow tau type :  type

>> Pow:  [(tau_1:type) => (---:type)] {move 0}

open

     declare xt2 in tau

>>      xt2:  in tau {move 2}

     construct tausub xt2 :  prop

>>      tausub:  [(xt2_1:in tau) => (---:prop)] {move 
>>        1}

     close

construct Setc tau tausub : in Pow tau

>> Setc:  [(tau_1:type),(tausub_1:[(xt2_2:in tau_1) => 
>>   (---:prop)]) => (---:in Pow(tau_1))] {move 0}

declare Ac  in Pow tau

>> Ac:  in Pow(tau) {move 1}

construct Ec tau xt Ac :prop

>> Ec:  [(tau_1:type),(xt_1:in tau_1),(Ac_1:in Pow(tau_1)) 
>>   => (---:prop)] {move 0}

declare ev1 that tausub xt

>> ev1:  that tausub(xt) {move 1}

declare ev2 that tau Ec xt tau Setc tausub

>> ev2:  that Ec(tau,xt,(tau Setc tausub)) {move 1}

construct Compc1 tau xt ,tausub, ev1 :  that tau Ec xt tau Setc tausub

>> Compc1:  [(tau_1:type),(xt_1:in tau_1),(tausub_1:[(xt2_2:in 
>>   tau_1) => (---:prop)]),(ev1_1:that tausub_1(xt_1)) 
>>   => (---:that Ec(tau_1,xt_1,(tau_1 Setc tausub_1)))] 
>>   {move 0}

construct Compc2 tau xt ,tausub, ev2 :  that tausub xt

>> Compc2:  [(tau_1:type),(xt_1:in tau_1),(tausub_1:[(xt2_2:in 
>>   tau_1) => (---:prop)]),(ev2_1:that Ec(tau_1,xt_1,(tau_1 
>>   Setc tausub_1))) => (---:that tausub_1(xt_1))] {move 
>>   0}

declare Bc in Pow tau

>> Bc:  in Pow(tau) {move 1}

open

     declare xt1 in tau

>>      xt1:  in tau {move 2}

     declare xtina1 that tau Ec xt1 Ac

>>      xtina1:  that Ec(tau,xt1,Ac) {move 2}

     construct aincb xt1 xtina1 :  that tau Ec xt1 Bc

>>      aincb:  [(xt1_1:in tau),(xtina1_1:that Ec(tau,
>>        xt1_1,Ac)) => (---:that Ec(tau,xt1_1,Bc))] {move 
>>        1}

     declare xtinb1 that tau Ec xt1 Bc

>>      xtinb1:  that Ec(tau,xt1,Bc) {move 2}

     construct binca xt1 xtinb1 :  that tau Ec xt1 Ac

>>      binca:  [(xt1_1:in tau),(xtinb1_1:that Ec(tau,
>>        xt1_1,Bc)) => (---:that Ec(tau,xt1_1,Ac))] {move 
>>        1}

     close

construct Extc tau Ac Bc , aincb, binca : that (Pow tau) eqt Ac Bc

>> Extc:  [(tau_1:type),(Ac_1:in Pow(tau_1)),(Bc_1:in 
>>   Pow(tau_1)),(aincb_1:[(xt1_2:in tau_1),(xtina1_2:that 
>>   Ec(tau_1,xt1_2,Ac_1)) => (---:that Ec(tau_1,xt1_2,Bc_1))]),
>>   (binca_1:[(xt1_3:in tau_1),(xtinb1_3:that Ec(tau_1,
>>   xt1_3,Bc_1)) => (---:that Ec(tau_1,xt1_3,Ac_1))]) => 
>>   (---:that eqt(Pow(tau_1),Ac_1,Bc_1))] {move 0}

% arrow type constructor

construct => tau sigma : type

>> =>:  [(tau_1:type),(sigma_1:type) => (---:type)] {move 
>>   0}

open

     declare var in tau

>>      var:  in tau {move 2}

     construct lambdabody var : in sigma

>>      lambdabody:  [(var_1:in tau) => (---:in sigma)] 
>>        {move 1}

     close

construct Lambda tau sigma lambdabody : in tau => sigma

>> Lambda:  [(tau_1:type),(sigma_1:type),(lambdabody_1:[(var_2:in 
>>   tau_1) => (---:in sigma_1)]) => (---:in (tau_1 => sigma_1))] 
>>   {move 0}

declare Fc in tau => sigma

>> Fc:  in (tau => sigma) {move 1}

declare Gc in tau => sigma

>> Gc:  in (tau => sigma) {move 1}

declare xt2 in tau

>> xt2:  in tau {move 1}

construct Applyc tau sigma Fc, xt2 :  in sigma

>> Applyc:  [(tau_1:type),(sigma_1:type),(Fc_1:in (tau_1 
>>   => sigma_1)),(xt2_1:in tau_1) => (---:in sigma_1)] 
>>   {move 0}

construct Beta tau sigma lambdabody, xt2 :
    that sigma eqt Applyc tau sigma (Lambda tau sigma lambdabody) xt2
    lambdabody xt2

>> Beta:  [(tau_1:type),(sigma_1:type),(lambdabody_1:[(var_2:in 
>>   tau_1) => (---:in sigma_1)]),(xt2_1:in tau_1) => (---:that 
>>   eqt(sigma_1,Applyc(tau_1,sigma_1,Lambda(tau_1,sigma_1,
>>   lambdabody_1),xt2_1),lambdabody_1(xt2_1)))] {move 
>>   0}

%  There remains extensionality for arrow types

open

     declare xt3 in tau

>>      xt3:  in tau {move 2}

     construct sameval xt3 : 
     that sigma eqt (Applyc tau sigma Fc xt3) (Applyc tau sigma Gc xt3)

>>      sameval:  [(xt3_1:in tau) => (---:that eqt(sigma,
>>        Applyc(tau,sigma,Fc,xt3_1),Applyc(tau,sigma,Gc,
>>        xt3_1)))] {move 1}

     close

construct Extfnc tau sigma Fc Gc sameval : that (tau => sigma) eqt Fc Gc

>> Extfnc:  [(tau_1:type),(sigma_1:type),(Fc_1:in (tau_1 
>>   => sigma_1)),(Gc_1:in (tau_1 => sigma_1)),(sameval_1:[(xt3_2:in 
>>   tau_1) => (---:that eqt(sigma_1,Applyc(tau_1,sigma_1,
>>   Fc_1,xt3_2),Applyc(tau_1,sigma_1,Gc_1,xt3_2)))]) => 
>>   (---:that eqt((tau_1 => sigma_1),Fc_1,Gc_1))] {move 
>>   0}

% add dependent product and dependent function types, which 

% allow internalization of abstraction sorts of the Lestrade framework.

% declarations for dependent types

open

     declare ys5 in tau

>>      ys5:  in tau {move 2}

     construct Rhofun ys5 : type

>>      Rhofun:  [(ys5_1:in tau) => (---:type)] {move 
>>        1}

     close

% dependent product construction

construct Xx tau Rhofun : type

>> Xx:  [(tau_1:type),(Rhofun_1:[(ys5_2:in tau_1) => (---:type)]) 
>>   => (---:type)] {move 0}

declare xt5 in tau

>> xt5:  in tau {move 1}

declare ys5 in Rhofun xt5

>> ys5:  in Rhofun(xt5) {move 1}

construct paird tau Rhofun, xt5 ys5 : in tau Xx Rhofun

>> paird:  [(tau_1:type),(Rhofun_1:[(ys5_2:in tau_1) => 
>>   (---:type)]),(xt5_1:in tau_1),(ys5_1:in Rhofun_1(xt5_1)) 
>>   => (---:in (tau_1 Xx Rhofun_1))] {move 0}

declare zp5 in tau Xx Rhofun

>> zp5:  in (tau Xx Rhofun) {move 1}

construct Pi1 tau Rhofun, zp5 :  in tau

>> Pi1:  [(tau_1:type),(Rhofun_1:[(ys5_2:in tau_1) => 
>>   (---:type)]),(zp5_1:in (tau_1 Xx Rhofun_1)) => (---:in 
>>   tau_1)] {move 0}

construct Pi2 tau Rhofun, zp5 : in Rhofun (Pi1 tau Rhofun, zp5)

>> Pi2:  [(tau_1:type),(Rhofun_1:[(ys5_2:in tau_1) => 
>>   (---:type)]),(zp5_1:in (tau_1 Xx Rhofun_1)) => (---:in 
>>   Rhofun_1(Pi1(tau_1,Rhofun_1,zp5_1)))] {move 0}

construct Xxexact tau Rhofun, zp5 : 
   that (tau Xx Rhofun) eqt zp5, paird tau Rhofun,
    (Pi1 tau Rhofun, zp5) (Pi2 tau Rhofun, zp5)

>> Xxexact:  [(tau_1:type),(Rhofun_1:[(ys5_2:in tau_1) 
>>   => (---:type)]),(zp5_1:in (tau_1 Xx Rhofun_1)) => (---:that 
>>   eqt((tau_1 Xx Rhofun_1),zp5_1,paird(tau_1,Rhofun_1,
>>   Pi1(tau_1,Rhofun_1,zp5_1),Pi2(tau_1,Rhofun_1,zp5_1))))] 
>>   {move 0}

% dependent function type constructor

construct =>> tau Rhofun : type

>> =>>:  [(tau_1:type),(Rhofun_1:[(ys5_2:in tau_1) => 
>>   (---:type)]) => (---:type)] {move 0}

open

     declare var in tau

>>      var:  in tau {move 2}

     construct lambdabodyd var : in Rhofun var

>>      lambdabodyd:  [(var_1:in tau) => (---:in Rhofun(var_1))] 
>>        {move 1}

     close

construct Lambdad tau Rhofun, lambdabodyd : in tau =>> Rhofun

>> Lambdad:  [(tau_1:type),(Rhofun_1:[(ys5_2:in tau_1) 
>>   => (---:type)]),(lambdabodyd_1:[(var_3:in tau_1) => 
>>   (---:in Rhofun_1(var_3))]) => (---:in (tau_1 =>> Rhofun_1))] 
>>   {move 0}

declare Fd in tau =>> Rhofun

>> Fd:  in (tau =>> Rhofun) {move 1}

declare Gd in tau =>> Rhofun

>> Gd:  in (tau =>> Rhofun) {move 1}

declare xt6 in tau

>> xt6:  in tau {move 1}

construct Applyd tau Rhofun, Fd, xt6 :  in Rhofun xt6

>> Applyd:  [(tau_1:type),(Rhofun_1:[(ys5_2:in tau_1) 
>>   => (---:type)]),(Fd_1:in (tau_1 =>> Rhofun_1)),(xt6_1:in 
>>   tau_1) => (---:in Rhofun_1(xt6_1))] {move 0}

construct Betad tau Rhofun, lambdabodyd, xt6 : 
   that (Rhofun xt6) eqt Applyd tau Rhofun, 
   (Lambdad tau Rhofun, lambdabodyd) xt6 lambdabodyd xt6

>> Betad:  [(tau_1:type),(Rhofun_1:[(ys5_2:in tau_1) => 
>>   (---:type)]),(lambdabodyd_1:[(var_3:in tau_1) => (---:in 
>>   Rhofun_1(var_3))]),(xt6_1:in tau_1) => (---:that eqt(Rhofun_1(xt6_1),
>>   Applyd(tau_1,Rhofun_1,Lambdad(tau_1,Rhofun_1,lambdabodyd_1),
>>   xt6_1),lambdabodyd_1(xt6_1)))] {move 0}

%  There remains extensionality for arrow types

open

     declare xt7 in tau

>>      xt7:  in tau {move 2}

     construct samevald xt7 : 
       that (Rhofun xt7) eqt (Applyd tau Rhofun, Fd xt7) 
      (Applyd tau Rhofun, Gd xt7)

>>      samevald:  [(xt7_1:in tau) => (---:that eqt(Rhofun(xt7_1),
>>        Applyd(tau,Rhofun,Fd,xt7_1),Applyd(tau,Rhofun,
>>        Gd,xt7_1)))] {move 1}

     close

construct Extfnd tau Rhofun,  Fd Gd samevald : 
   that (tau =>> Rhofun) eqt Fd Gd

>> Extfnd:  [(tau_1:type),(Rhofun_1:[(ys5_2:in tau_1) 
>>   => (---:type)]),(Fd_1:in (tau_1 =>> Rhofun_1)),(Gd_1:in 
>>   (tau_1 =>> Rhofun_1)),(samevald_1:[(xt7_3:in tau_1) 
>>   => (---:that eqt(Rhofun_1(xt7_3),Applyd(tau_1,Rhofun_1,
>>   Fd_1,xt7_3),Applyd(tau_1,Rhofun_1,Gd_1,xt7_3)))]) => 
>>   (---:that eqt((tau_1 =>> Rhofun_1),Fd_1,Gd_1))] {move 
>>   0}

%% further remarks about internalization:  Pow One
%% implements prop.  Then all the propositional operations
%% correspond to type constructors just given, with all types
% that p actually being identified with either One or Empty.

construct Empty : type

>> Empty:  [(---:type)] {move 0}

declare xnot in Empty

>> xnot:  in Empty {move 1}

construct notthere xnot : that ??

>> notthere:  [(xnot_1:in Empty) => (---:that ??)] {move 
>>   0}

%% this means that the entire logical framework can 
%% be internalized, at least in its classical version:
%% the full type system of abstraction sorts
% can be studied internally to Lestrade.


quit

\end{verbatim}

\newpage

\section{Appendix:  the source {\tt lestrade.sml} as of October 22 2016, 7  pm  Boise time}

{\tiny

\begin{verbatim}



(*  implement the description of my Automath logic in 
implementation.tex.  *)

(* full version with type inference features:  rewriting is installed, and I intend to install
user-entered lambda terms and deduced arguments.  These features are more dangerous (the implementation
is more likely to be buggy), so I am maintaining the old version in lestrade_basic.fsml.  I have just
installed the capability of deducing implicit arguments; I did preserve a snapshot of the software
before this major upgrade. *)

(* there are now internal functions which will turn off new features and give the behavior of old
versions:  this goes with my conviction that the main type checking engine is now working quite reliably,
and any new features should touch it in as few places as possible, if at all, and should be possible to
turn off by disabling a few functions *)

(*
dated notes

10/21 The readfile command can be issued in the interface. This seems to work reasonably well; testing
still needed.

10/20  Parser refined so that the colon is always optional in the construct command:
the parser knows that an argument list has ended when it encounters a reserved identifier.

10/11 another modification to consider:  allow more general deduction of an abstraction applied to
a concretely given list of arguments in the implicit arguments feature.

10/10 Construction of function constants by in effect currying
(shortening argument lists to defined or constructed operators)
is now implemented, which allows much more compact construction
of functions with internal parameters.

10/9  another projected change:  allow other methods of expressing abstraction arguments,
 both application terms with missing arguments (interpreted as if curried) and 
non-polymorphic abstractions (require that all types be deducible without explicit types 
on input variables, as in the rewrite system.  No code change yet.

10/2  projected possible changes, no modification in code yet.
A dependency system for type 0 postulates and definitions:  the
idea would be to be able to implement constructed objects and
functions by giving definitions with the right types.  Defined objects
depend on the primitive used in their declarations; constructed objects
seem to depend in a sense on other constructed objects which mention
them in their sorts:  at least, it seems that an implement command
would require implementations of all such constructions to be given.

The other possibility I am considering is overloading.

User-entered lambda terms remain a desideratum if I can figure out
to do them neatly.

Testing of rewrited is needed.

Work on executing all commands within the interface.  Being able to
run readfile in some way in the interface would save lots of headaches.
Setting margins in the interface is just common sense.  A device
for allowing carriage returns in command lines and enforcing the margins
on command lines would be useful.

9/27 major upgrades:  changed display so that 
implicit arguments of functions applied in sort displays
are not shown (unless the command showimplicit is run).  
Fixed a bug in the implicit arguments mechanism, so that
yet more arguments can be deduced successfully.

9/7/2016 change world n to move n.  Slight change in metaphor.

8/20/2016  Attempt to debug possible problems with 
interaction of implicit argument and rewriting features.

This needs testing:  examples of rewritec and rewrited 
commands with implicit arguments present are needed.
My belief is that rewritec was already set up to work correctly;  
rewrited needed an argument list fix inserted,
which is now there but needs testing.

8/12/2016 Installed the version toggles as user commands.  They do make sense.

Major (and possibly dangerous) upgrade:  abstraction types 
and lambda terms are now pretty printed to some extent:
breaks are indented more if they are inside abstraction types
 and lambda terms and the display functions always
break at the end of such a term (possibly picking up a parenthesis
 or two or a parenthesis and a comma).

Guarded rewrite patterns from possible bound variable conflicts 
by using negative indexes on variables in their namespace.

More cleanup of program text and comments.

8/11/2016  Continuing reading the document with an eye to intelligibility.
  I did experimentally apply the reindexing of bound variables to actual 
declarations which was being applied to displays.  I believe this is sound, 
and the files actually run;
this fixes an esthetic problem.  I also fixed a problem with displays of
 showrecent and presumably showall.  The program should also break out of
 scripts with errors faster.

The reindexing change is perilous:  I need to be sure that all
situations where substitutions are made into abstraction types
 are preceded by bound variable renaming.

8/10/2016  Planning to read through to document the code. 
 Killing the old comments; new ones will accumulate, no doubt 

--END dated notes *)


(* utility functions *)


(* controls the greeting message when you enter the interface *)

val GREETED=ref false;

(* the file to which all system activity is logged *)

val LOGFILE = ref (TextIO.openOut("default"));

(* close the log file *)

fun closelog() = (TextIO.flushOut(!LOGFILE);
TextIO.closeOut(!LOGFILE);LOGFILE:=TextIO.openOut("default"));

(* say = system messages.   These go to standard output and also
 as temporary comments (ignored and not persisting when logs are 
executed) in the log. *)

fun Flush() = (TextIO.flushOut(TextIO.stdOut);TextIO.flushOut(!LOGFILE));
fun say s = (TextIO.output(TextIO.stdOut,"\nInspector Lestrade says:  "^s^"\n\n");
Flush();TextIO.output(!LOGFILE,"\n>> Inspector Lestrade says:  "^s^"\n\n");Flush());

(* saypause = system messages which are errors; these will terminate scripts being run *)

val BREAKOUT = ref false;

fun saypause s = (say (s^"\n>> Hit return to continue\n\n");if (!BREAKOUT) = false 
then TextIO.inputLine(TextIO.stdIn) else "";BREAKOUT:=true)

(* the current version.  This is also the greeting the system gives 
(and puts at the head of the log) when you enter the interface *)

fun versiondate() = say 
("\n>>   Welcome to the Lestrade Type Inspector,\n>>    "^
"full version of 10/21/2016 (readfile is interface command)\n>>   4 pm Boise time\n");

(* Lestrade version toggles.  These are now internal Lestrade commands. *)

val REWRITEVER = ref true;

val IMPLICITVER = ref true;

fun basic() = (REWRITEVER := false;  IMPLICITVER:=false);

fun explicit() = (REWRITEVER:=true;IMPLICITVER:=false);

fun fullversion() = (REWRITEVER := true;IMPLICITVER:=true);

(* find the type associated with an identifier in a world or argument list *)

(* an item (n,t,u) in one of these lists has n the numerical
age of the item, t a term (usually a suitably packaged identifier, but there is
an exception) and u a type. *)

fun find s nil = nil |

find s ((n,t,u)::L) = if s=t then [u] else find s L;

(* more general find function for, e.g., matching lists *)

fun abstractfind s nil = nil |

abstractfind s ((t,u)::L) = if s=t then [u] else abstractfind s L;

(* drop function for general lists *)

fun abstractdrop s nil = nil |

abstractdrop s ((t,u)::L) = 

   if s = t then abstractdrop s L
   else (t,u)::(abstractdrop s L);

(* version which overwrites earlier matches *)

fun abstractmerge nil L = L |

   abstractmerge ((s,t)::M) L =

      (s,t)::(abstractmerge M (abstractdrop s L));

(* intersection of arguments in two type dec lists -- no compatibility check *)

fun intersection nil L = nil |

intersection ((n,t,u)::L) M = if find t M <> nil then ((n,t,u)::(intersection L M))
else intersection L M;

fun union nil L = L |

union ((n,t,u)::L) M = if find t M = nil then ((n,t,u)::(union L M))
else union L M;

fun unionoflist nil = nil |

unionoflist (L::M) = union L (unionoflist M);

(* find the age of the declaration of a term in a world.  
This is used to check the requirement that (explicit) parameters 
appear in a function definition or declaration in the order in which they are declared.   
The age of a defined abstraction declaration is always zero. *)

fun age s nil = nil |

age s ((n,t,u)::L) = if s=t then [n] else age s L;

(* check a condition for all elements of a list *)

fun testall test nil = true |

testall test (s::L) = test s andalso testall test L;

(* namespace reindexing utilities *)

(* find a number in a list of non negative integers:  used in namespace
reindexing for display in showdec *)

(* These functions are now being applied to
types entered in declarations, not only to their displays. *)

fun Find n nil = ~1 |

Find n ((p,q)::L) = if n=p then q else Find n L;


val Freshindex = ref 0;

val Indexlist = ref [(~1,~1)];

val _ = Indexlist := nil;

fun Renumber n = if n=0 then 0 else let val N = Find n (!Indexlist) in

if N= ~1

then (Indexlist:=(n,!Freshindex)::(!Indexlist);Freshindex:=1+(!Freshindex);(!Freshindex)-1)

else N end

fun Reset() = (Freshindex := 1;Indexlist := nil);

(* the internal representation of basic types and objects *)

(* this type declaration seems to be marvellous, 
as it seems to capture all sorts of things that we talk about, all the way from
mathematical objects to worlds *)

datatype EntType (* basic sorts of entities *) =

obj (* mathematical objects *) |

prop (* propositions *) |

TYPE (* types of mathematical object *) |

that of Entity (* sort of proofs of a given proposition *) |

IN of Entity | (* the sort associated with a particular object of sort TYPE *)

error

and AbstType =

World of (int*Argument*Type) list (* this is the metasort of worlds, and    
                      also of
                      types assigned to abstractions *)

(* the integer is ``age", for ensuring correct order in
argument lists. Age is 0 for defined notions *)

and Type (* general sorts of entities and abstractions *) =

EType of EntType (* type of an entity *)|

AType of AbstType (* type of an abstraction *)

and Entity (* first order entities: 
typed and untyped objects, propositions, types, and proofs *) = 

Unknown (* postulated, unknown--this is the --- 
appearing in the type in an abstraction declaration *) |

Error |

Ent of string*int (* the numeral indicates namespace when nonzero, 
used for renaming bound variables in dependent types and lambda terms *) |

App of string*int*(Argument list)  (* the numeral again indicates namespace *)

and Argument (* first and second order objects, 
entities and abstractions, as they appear in argument lists *) =

EntArg of Entity |

AbstArg of string*int | (* numeral is again namespace *)

Lambda of AbstType;  (* lambda terms appearing in types *)

(* the user never enters a lambda-term as an abstraction argument, only an identifier, 
but the system generates lambda-terms when an identifier passes out of scope, 
and also in the implicit arguments mechanism *)

(* change all namespace indices to their additive inverses;  used for rewriting patterns *)

fun Negindex1 obj = obj |(* mathematical objects *) 

Negindex1 prop = prop | (* propositions *) 

Negindex1 TYPE = TYPE (* types *) |

Negindex1 (that E) = that (Negindex4 E)  (* type of proofs of a proposition *) |

Negindex1 (IN E) = IN (Negindex4 E) |

Negindex1 error = error

and Negindex2 (World L) = World(map (fn (i,A,T) => (i,Negindex5 A, Negindex3 T)) L)

and Negindex3 (EType E) = EType (Negindex1 E) (* type of an entity *)|

Negindex3 (AType A) = AType (Negindex2 A) (* type of an abstraction *)

and Negindex4 Unknown = Unknown (* postulated, unknown--this is the --- 
appearing in the type in an abstraction declaration *) |

Negindex4 Error = Error |

Negindex4 (Ent(s,n)) = (Ent(s,~n)) (* the numeral indicates namespace, 
used for renaming bound variables in dependent types and lambda terms *) |

Negindex4 (App(s,n,L)) = (App(s,~n, map Negindex5 L)) 
      (* the numeral again indicates namespace *)

and Negindex5 (EntArg E) = EntArg(Negindex4 E) |

Negindex5 (AbstArg(s,n)) = AbstArg(s,~n)| (* numeral is again namespace *)

Negindex5 (Lambda A) = Lambda (Negindex2 A);  (* lambda terms appearing in types *)

(* namespace reindexing function:  this is a utility which avoids runaway indices
on new bound variables in declaration displays *)

(* This function is now applied
to the actual recorded types; in earlier versions it was only applied to displays *)

(* this required care that renaming of bound variables (the
function renamespace) be applied before *any* substitution into
an abstraction type or lambda term *)

fun Reindex1 obj = obj |(* mathematical objects *) 

Reindex1 prop = prop | (* propositions *) 

Reindex1 TYPE = TYPE (* types *) |

Reindex1 (that E) = that (Reindex4 E)  (* type of proofs of a proposition *) |

Reindex1 (IN E) = IN (Reindex4 E) |

Reindex1 error = error

and Reindex2 (World L) = World(map (fn (i,A,T) => (i,Reindex5 A, Reindex3 T)) L)

and Reindex3 (EType E) = EType (Reindex1 E) (* type of an entity *)|

Reindex3 (AType A) = AType (Reindex2 A) (* type of an abstraction *)

and Reindex4 Unknown = Unknown (* postulated, unknown--this is the ---
 appearing in the type in an abstraction declaration *) |

Reindex4 Error = Error |

Reindex4 (Ent(s,n)) = (Ent(s,Renumber n)) 
(* the numeral indicates namespace, 
used for renaming bound variables in dependent types and lambda terms *) |

Reindex4 (App(s,n,L)) = (App(s,Renumber n, map Reindex5 L)) 
(* the numeral again indicates namespace *)

and Reindex5 (EntArg E) = EntArg(Reindex4 E) |

Reindex5 (AbstArg(s,n)) = AbstArg(s,Renumber n)| (* numeral is again namespace *)

Reindex5 (Lambda A) = Lambda (Reindex2 A);  (* lambda terms appearing in types *)

(* CONTEXT is the list of worlds in use *)

val CONTEXT = ref  [World nil, World nil];

(* the list of rewrite rules declared or justified *)

val REWRITES = ref [[("bogus",(Unknown,Unknown))],[("bogus",(Unknown,Unknown))]];

val _ = REWRITES := [nil,nil];

(* the list of names attached to worlds in the current sequence *)

val WORLDNAMES = ref ["1","0"];

(* the list of worlds which have been saved *)

val SAVEDWORLDS = ref [(!WORLDNAMES,hd(!CONTEXT))];

val _ = SAVEDWORLDS:=nil;

(* saved rewrite rules, really a component of the previous though implemented separately *)

val SAVEDREWRITES = ref [(!WORLDNAMES,hd(!REWRITES))];

val _ = SAVEDREWRITES:=nil;

(* full theories, just the declarations in world 0, saved when files are run
and recovered by the load command *)

val SAVEDTHEORIES = ref [("bogus",(0,0,hd(rev(!CONTEXT))))];

val _ = SAVEDTHEORIES := nil;

(* determine the number of the current world and attach 
its name if it has a non-default name (not simply its numeral index)  *)

fun worldname0 n nil = "bogus" |

worldname0 0 L = if hd L = makestring(length L-1) then "" else ":"^(hd L) |

worldname0 n L = worldname0 (n-1) (tl L);

fun worldname n = worldname0 n (!WORLDNAMES);

(* does a world (other than world 0) have a trivial name (i.e., simply its numerical index) *)

fun defaultworld L = (length L >=2 andalso hd L = makestring((length L) -1));

(* recover lists of worlds which can be opened using
open or clearcurrent commands *)

fun savedfor M nil = "" |

   savedfor M ((s,t)::L) = 

   if tl s = M then (hd s)^"\n"^(savedfor M L)
   else (savedfor M L);

fun savedforopen() = savedfor (!WORLDNAMES) (!SAVEDWORLDS);

fun savedforclearcurrent() = savedfor (tl(!WORLDNAMES)) (!SAVEDWORLDS);

(* pretty printing *)

val EXTRAINDENTS = ref 0;  (* keep track of extra indentation needed in display of abstraction sorts *)

(* break after the first comma, colon, or space after n characters *)

(* also, indent displays more deeply as one moves to deeper worlds *)

(* also, indent more deeply at breaks inside abstraction sorts or lambda terms
and always break after a closing bracket (possibly picking up some following punctuation) *)

fun tolinebreak n nil = nil |

tolinebreak n (#"]":: #")":: #",":: L) = (EXTRAINDENTS:=(!EXTRAINDENTS)-1;[#"]", #")", #","]) |

tolinebreak n (#"]":: #")":: #")":: L) = (EXTRAINDENTS:=(!EXTRAINDENTS)-1;[#"]", #")", #")"]) |

tolinebreak n (#"]":: #")":: L) = (EXTRAINDENTS:=(!EXTRAINDENTS)-1;[#"]", #")"]) |

tolinebreak n (#"\n":: #"\n":: L) = [#"\n", #"\n"]@(tolinebreak 0 L) |

tolinebreak 0 (c::L) = (
if c = #"[" then EXTRAINDENTS := 1+(!EXTRAINDENTS) else ();
if c = #"]" then EXTRAINDENTS := (!EXTRAINDENTS)-1 else ();
if c = #"," orelse c = #":" orelse c= #"]" orelse (c = #" " 
andalso (L = nil orelse hd L <> #" ")) then [c] else c::(tolinebreak 0 L)) |

tolinebreak n (#" ":: #" ":: L) = #" ":: #" ":: tolinebreak (n-2) L |


tolinebreak n (x:: #" ":: #" ":: L) = 

(
if x = #"[" then EXTRAINDENTS := 1+(!EXTRAINDENTS) else ();
if x = #"]" then EXTRAINDENTS := (!EXTRAINDENTS)-1 else ();tolinebreak n (x:: #" " :: L)) |

tolinebreak n (#"\n" :: #" " :: L) = tolinebreak n (#"\n" :: L) |

tolinebreak n (#"\n" :: L) = tolinebreak n (#" "::L)  |

tolinebreak n (c::L) = (
if c = #"[" then EXTRAINDENTS := 1+(!EXTRAINDENTS) else ();
if c = #"]" then EXTRAINDENTS := (!EXTRAINDENTS)-1 else ();
if c= #"]" then [#"]"] else
c::(tolinebreak (n-1) L));

fun restlinebreak n nil = nil |

restlinebreak n (#"]":: #")":: #",":: L) = L |

restlinebreak n (#"]":: #")":: #")":: L) = L |

restlinebreak n (#"]":: #")":: L) = L |

restlinebreak n (#"\n":: #"\n":: L) = restlinebreak 0 L |

restlinebreak 0 (c::L) = if c = #"," orelse c= #":" orelse c = #"]" 
orelse (c = #" " andalso (L = nil orelse hd L <> #" "))  then L else (restlinebreak 0 L) |

restlinebreak n (#" ":: #" ":: L) = restlinebreak (n-2) L |

restlinebreak n (x:: #" ":: #" ":: L) = restlinebreak n (x:: #" " :: L) |

restlinebreak n (#"\n" :: #" " :: L) = restlinebreak n (#"\n" :: L) |

restlinebreak n (#"\n" :: L) = restlinebreak n (#" " :: L) |

restlinebreak n (c::L) = if c = #"]" then L else (restlinebreak (n-1) L);

val MARGIN = ref 40;

fun setmargin n = MARGIN:= n;

fun indents n = if n<= 0 then "" else (indents (n-1))^"     ";

fun INDENTS() = indents (((length(!CONTEXT)-2))+(!EXTRAINDENTS));

fun initial s = implode(tolinebreak(!MARGIN)(explode s));

fun final s = implode(restlinebreak(!MARGIN)(explode s));

fun despace0 (#" "::L) = despace0 L |

despace0 L = L

and despace s = implode(despace0(explode s)) 

fun prettyprint s = let val I = INDENTS() in 
(initial ((I)^(despace s)))^(if final (I^(despace s)) = "" then ""
 else "\n   "^(I)^(prettyprint(final ((I)^(despace s))))) end;

(* this produces the same output as prettyprint with each output
formatted for the log file as a temporary comment *)

fun prettyprint2 s = let val I = INDENTS() in (initial 
((I^(despace s))))^(if final ((I)^(despace s))  = "" then "" else
"\n>>   "^(prettyprint2(final ((I^(despace s)))))) end;

fun prettyprint3 s = "\n>> "^(prettyprint2 s);

(* post a pretty printed piece of Lestrade notation to standard output *)

fun say0 s= (EXTRAINDENTS:=0;TextIO.output(TextIO.stdOut,(prettyprint s)^"\n");Flush());

(* post a pretty printed piece of Lestrade notation to standard output and the log file *)

fun say1 s= (say0 s;EXTRAINDENTS:=0;TextIO.output(!LOGFILE,(prettyprint3 s)^"\n");Flush());

(* find the Argument reference of the string s in
the list of worlds L (this returns the type of the 
entity or argument named by s if there is one, which does
contain enough information to tell whether the object
is declared and whether it is an entity or an abstraction *)

fun Find s (World L) = find s L;

(* similarly this returns the age of any declaration in any world of s -- 
the main use of this is to identify defined abstractions, which have age 0 *)

fun Age s (World L) = age s L;

(* stringdef returns the singleton list of 
the type of an identifier using a given world list argument,
paired with the numerical index of the world it is found in,
 or nil for error *)

fun pi1(x,y)=x;


fun pi2(x,y)=y;

fun stringdef s nil = nil |

stringdef s L =

let val A = Find (EntArg(Ent (s,0))) (hd L) in

if A <> nil then  [(hd A,0)] else

let val B = Find (AbstArg (s,0)) (hd L) in

if B <> nil then [(hd B,0)] else

let val C = stringdef s (tl L) in

if C = nil then nil else [(pi1(hd C),pi2(hd C)+1)] end

end

end;

(* this returns the singleton list of the type of an identifier in the Lestrade context,
paired with the numerical index of the world it is found in,
or nil for error *)

fun stringtype s =  stringdef s (!CONTEXT);

(* build function to extract explicit argument list *)

val DISPLAYIMPLICIT = ref false

fun showimplicit() = DISPLAYIMPLICIT := true

fun hideimplicit() = DISPLAYIMPLICIT := false

fun purgeimplicit (World L) nil = nil |

purgeimplicit (World ((n,EntArg(Ent(s,m)),t)::L)) (u::M) =

if s<> "" andalso hd(explode s) = #"." then purgeimplicit (World L) M
else (u::(purgeimplicit (World L) M)) |

purgeimplicit (World ((n,AbstArg(s,m),t)::L)) (u::M) =

if s<> "" andalso hd(explode s) = #"." then purgeimplicit (World L) M
else (u::(purgeimplicit (World L) M)) |

purgeimplicit (World(x::L)) (u::M) = u::(purgeimplicit (World L) M);

(* get abstraction type (world) from a general type *)

fun getabstype (AType x) = x|

getabstype x = (saypause "getabstype error";World nil);


fun explicitlist s n L = if (!DISPLAYIMPLICIT) orelse n<>0 then L
else purgeimplicit (getabstype((pi1(hd(stringtype s))))) L;

(* display functions *)

fun isentarg (EntArg x) = true |

isentarg x = false;

(* display1 displays entity types *)

fun display1 obj = "obj" |

display1 prop = "prop" |

display1 TYPE = "type" |

display1 (that P) = "that "^(display2 P) |

display1 (IN P) = "in "^(display2 P) |

display1 error = "error"

(* display2 displays entity terms.  Note that --- 
is reserved for the pseudo-entity Unknown used
 as output for primitive abstractions, and ??? for Error. *)

and display2 (Ent(s,n)) = if n=0 then s else s^"_"^(makestring n)|

display2(App(s,n,nil)) = display2(Ent(s,n)) |

display2 (App(s,n,LL)) = 

let val L = explicitlist s n LL in

if length L = 2 andalso isentarg (hd L) then 
"("^(display4 (hd L))^" "^s^(if n=0 then ""
 else "_"^(makestring n))^" "^(display4(hd(tl L)))^")"

else s^(if n=0 then "" else "_"^(makestring n))^"("^(display3 L) end |

display2 Unknown = "---" |

display2 Error = "???"

(* display3 displays argument lists *)

and display3 [a] = (display4 a)^")" |

display3 (a::L) = (display4 a)^","^(display3 L) |

display3 nil = "*?*?)"

and display4 (EntArg x) = display2 x |

display4 (AbstArg(s,n)) = display2(Ent(s,n)) |

display4 (Lambda x) = "["^(display5 x)

(* display5 displays dependent types and abstraction terms.
  This is the part of the output language about which the parser knows nothing. *)

and display5 (World [(n,a,t)]) = "("^(display4 a)^":"^(display6 t)^")]" |

display5 (World((n,a,t)::L)) = "("^(display4 a)^":"^(display6 t)^
(if length L = 1 then ") => " else "),")^(display5 (World L)) |

display5 (World nil) = "(?*?*?*?)"

(* display6 displays general types *)

and display6 (EType x) = display1 x |

display6 (AType(World x)) =  "["^(display5 (World x)); 

(* display a world.  This is the same type displayed
by display5, with different intent *)

fun displayworld (World nil) = "" |

displayworld (World((n,a,t)::L)) = (displayworld (World L))
^(INDENTS())^(display4 a)^":"^(display6 t)^"\n\n";

(* display the indexed list of worlds.  
displacement corrects world indices in the showrecent command below *)

fun displayworlds displacement nil = "\n\n" |
displayworlds displacement L = "\n\nMove "^(makestring(displacement+(length L)-1))
^(worldname (length (!WORLDNAMES)-(displacement+(length L))) )
^":\n\n"^(displayworld (hd L))^(displayworlds displacement (tl L));

(* display a rewrite list *)

fun displayarewritelist nil = "\n\n" |

displayarewritelist ((s,(t,u))::L) = (s^":  "
^(display2 t)^" := "^(display2 u)^"\n")^(displayarewritelist L);

fun displayrewrites0 nil = "\n\n" |

displayrewrites0 (L::M) = (displayarewritelist L)^(displayrewrites0 M);

fun displayrewrites() = say (displayrewrites0 (tl(!REWRITES)));

(* a utility -- coerce an argument to an entity. *)

fun deent (EntArg x) = x|
deent x = Error;

(* display all worlds -- this shows all declarations in detail *)

fun showall () = say0 (displayworlds 0 (!CONTEXT));

(* display the current world and the parent world 
-- these are the worlds in which you can actually make declarations,
 and the showall display is likely to be huge.  
So is the showrecent display if world 0 is displayed! *)

fun showrecent() = say0 
(displayworlds (length(!CONTEXT)-2) [hd(!CONTEXT),hd(tl(!CONTEXT))]);

(* user command:  open a new world (or a previously saved world) *)

fun Open s = (if defaultworld (!WORLDNAMES) 
andalso not(defaultworld(s::(!WORLDNAMES)))
 then saypause "Cannot follow default move with named move" else
(let val W = abstractfind (s::(!WORLDNAMES)) (!SAVEDWORLDS)
 and R= abstractfind (s::(!WORLDNAMES)) (!SAVEDREWRITES) in
if W = nil then
(CONTEXT := (World nil)::(!CONTEXT); REWRITES:= (nil::(!REWRITES)))
else (CONTEXT := (hd W)::(!CONTEXT);REWRITES := (hd R)::(!REWRITES))  end; 
WORLDNAMES := s::(!WORLDNAMES)));



(* display the declaration information of an identifier on standard output.  Send
this to the log as well. Bound variable indices are rationalized.  *)

fun showdec s = let val S = stringtype s in 

if S = nil then saypause (s^" is not declared")

else (Reset(); say1(s^":  "^(display6(Reindex3(pi1 (hd S))))^
" {move "^(makestring(length(!CONTEXT)-1-(pi2 (hd S))))
^(worldname(pi2(hd S)))^"}\n\n"))end;

(* utility for converting a world to a list of type declarations *)

fun deworld (World L) = L;

(* showdecs will display declarations one by one.  Useful if world 0
is being displayed.  It goes through the current world in order, then the
parent world in reverse order *)

fun showdecs0 nil = () |

showdecs0 ((i,EntArg(Ent(s,0)),t)::L) = 
(showdec s;if TextIO.input(TextIO.stdIn)="q\n" then ()
 else showdecs0 L)|

showdecs0 ((i,AbstArg(s,0),t)::L) = 
(showdec s;if TextIO.input(TextIO.stdIn)="q\n" then ()
 else showdecs0 L)|

showdecs0 L = ();

fun showdecs() = (say "Hit return after each declaration or q to quit"; 
say "Next move declarations"; showdecs0(deworld(hd(!CONTEXT)));
say "Present move declarations:"; showdecs0(rev(deworld(hd(tl(!CONTEXT))))));

(* age of an identifier in a given context packaged in a singleton list
 -- useful mainly for identifying defined abstractions *)

fun stringage s nil = nil |

stringage s L =

let val A = Age (EntArg(Ent (s,0))) (hd L) in

if A <> nil then A else

let val B = Age (AbstArg (s,0)) (hd L) in

if B <> nil then B else

stringage s (tl L)

end

end;

(* as the previous command, in the Lestrade context *)

fun stringAge s = stringage s (!CONTEXT);

(* utility identifies an entity type as opposed to an abstraction type *)

fun isenttype (EType x) = true |
isenttype (AType x) = false;

(* declaration checking for all types against a context given as an argument *)

fun (* check entity types *) deccheck1 L obj = true

| deccheck1 L prop = true

| deccheck1 L TYPE = true

| deccheck1 L (that P) = deccheck2 L P 

| deccheck1 L (IN P) = deccheck2 L P

| deccheck1 L error = false

and deccheck2 L Unknown = true |

deccheck2 L Error = false |

(* check entities *) deccheck2 L (Ent (s,n)) = n<>0 orelse let val S = stringdef s L in 
(if S = nil then saypause ("Did not find entity "^s^" (deccheck2)")else (); Flush();     

S <> nil andalso isenttype(pi1(hd S))) end |

deccheck2 L (App(s,n,M)) =
      (n<>0 orelse
      let val S = stringdef s  L in

(if S = nil then saypause ("Did not find abstraction "^s^" (deccheck 2)") else();  Flush();

      S <> nil andalso not (isenttype(pi1(hd S)))) end)
      andalso testall (deccheck3 L) M

      

and (* check arguments *) deccheck3 L (EntArg s) = deccheck2 L s

|  deccheck3 L (AbstArg (s,n)) = n<>0 orelse let val S = stringdef s  L 
in  (if S = nil then saypause ("Did not find abstraction "^s^" (deccheck3)")else();Flush();

      S <> nil andalso not (isenttype(pi1(hd S)))) end

| deccheck3 L (Lambda T) = deccheck4 L T

and (* check worlds *) deccheck4 L (World M) =

testall (fn (n,a,t) => deccheck3 L a andalso deccheck5 L t) M

and (* check general types *) deccheck5 L (EType x) = deccheck1 L x |

deccheck5 L (AType x) = deccheck4 L x;

(* the index for generating fresh names in new namespaces,
 for bound variable renaming in dependent types and abstractions *)

val NAMESERIAL = ref 0;

(* get a new namespace serial number *)

fun newnameserial() = (NAMESERIAL:=1+(!NAMESERIAL);(!NAMESERIAL))

(* utility for adding an item to a world at the beginning *)

fun addworld2 x ((World M)) = (World (x::M));

(* s is a declared identifier new in the current world *)

fun isnew s = Find s (hd(!CONTEXT)) <> nil;

(* s is a variable, new and not defined *)

fun isvariable s = isnew s andalso hd(Age s (hd(!CONTEXT))) <> 0;

(* utilities for taking apart an application term *)

fun appof (App(t,n,L)) = t |

appof x = "";

fun argsof (App(t,n,L)) = L |

argsof x = nil;

(* conditions for all types to be deducible in a term, 
used to test appropriate rewrite patterns *)

fun isdefvar x = isnew x andalso not (isvariable x);

(* here I exclude new defined variables from type rigid terms, which is appropriate for
patterns, but one might want to expand defined terms if one is using this for user-entered
lambda terms *)

fun typerigid0 weak (EntArg(Ent(s,0))) = 
(weak andalso not (isdefvar (EntArg(Ent(s,0)))))
 orelse not(isnew(EntArg(Ent(s,0)))) |

typerigid0 weak (EntArg(App(s,0,L))) =  
(weak andalso not (isdefvar(AbstArg(s,0))) 
andalso testall (typerigid0 false) L) orelse
(not(isnew(AbstArg(s,0))) andalso testall (typerigid0 true) L) |

typerigid0 weak (AbstArg(s,0)) = (weak 
andalso not (isdefvar (AbstArg(s,0)))) 
orelse not(isnew(AbstArg(s,0))) |

typerigid0 weak x = false;

fun typerigid x = typerigid0 false (EntArg x);

(* This is the central dependent type checking engine, not to be touched except
as absolutely necessary.  It appears to be quite stable *)

(* substitution and type checking all in one package *)

(* this is a huge block of recursively declared functions *)

(* all substitutions are of an argument for another argument,
in whatever type *)

(* substitution into general types *)

fun typesubs a A (EType t) = EType(etypesubs a A t) |

typesubs a A (AType t) = AType(atypesubs a A (renamespace t))

(* substitutions into worlds = abstraction types *)

and atypesubs a A ((World nil)) = (World nil) |

atypesubs a A ((World((n,b,t)::M))) =

addworld2 (n,argsubs a A b,typesubs a A t)  
(* argsubs on first component is used only for the defined value *)

(atypesubs a A ((World M)))

(* substitutions into entity types *)

and etypesubs a A obj = obj |

etypesubs a A prop = prop |

etypesubs a A (that P) = that (entsubs a A P)|

etypesubs a A TYPE = TYPE |

etypesubs a A (IN P) = IN (entsubs a A P) |

etypesubs a A error = error

(* substitutions into entity terms *)

(* notice that defined abstractions declared 
in the current world are expanded using defmatchcomp.  
Trivial substitutions are often made to enforce this expansion. 

The reason for this is that these substitutions
 are made for types to be recorded in the parent context, 
so abstractions defined in the current world must pass out of scope.
*)

and entsubs (EntArg x) (EntArg A) (Ent(s,n)) =

    if x = Ent(s,n) then A else Ent(s,n) |

entsubs (EntArg x) (EntArg A) (App(s,n,M)) =

    if n=0 andalso Age (AbstArg(s,n)) (hd(!CONTEXT)) = [0]

    then entsubs (EntArg x) (EntArg A) 

(defmatchcomp true (getabstype(hd(Find (AbstArg(s,n)) (hd(!
        CONTEXT)))) )
    (map(argsubs (EntArg x) (EntArg A)) M))

    else App(s,n,map(argsubs (EntArg x) (EntArg A)) M) |

entsubs (AbstArg(s,n)) (AbstArg(S,N)) (App(t,m,M)) =

   if m=0 andalso Age (AbstArg(t,m)) (hd(!CONTEXT)) = [0]
   then entsubs (AbstArg(s,n)) (AbstArg(S,N)) 
   (defmatchcomp true (getabstype(hd(Find (AbstArg(t,m)) (hd(!
        CONTEXT)))) )
    (map(argsubs (AbstArg(s,n)) (AbstArg(S,N))) M))

else if s=t andalso n=m 

then entsubs (AbstArg(s,n)) (AbstArg(S,N)) 
    (App(S,N,map(argsubs (AbstArg(s,n)) (AbstArg(S,N))) M))

else App(t,m,map(argsubs (AbstArg(s,n)) (AbstArg(S,N))) M) |

(* frank beta substitution -- I do not know if this can
actually happen *)  (* yes, it can, and it was scrambled *)

entsubs (AbstArg(s,n))(Lambda w) (App(t,m,M)) =

if m=0 andalso Age (AbstArg(t,m)) (hd(!CONTEXT)) = [0]
   then entsubs (AbstArg(s,n)) (Lambda w) 
   (defmatchcomp true (getabstype(hd(Find (AbstArg(t,m)) (hd(!
        CONTEXT)))) )
    (map(argsubs (AbstArg(s,n)) (Lambda w)) M))

else if s=t andalso n=m 

then entsubs (AbstArg(s,n)) (Lambda w)
    (defmatchcomp true w
    (map(argsubs (AbstArg(s,n)) (Lambda w)) M)) 

else App(t,m,map(argsubs (AbstArg(s,n)) (Lambda w)) M)|


entsubs x y T = T

(* substitutions into arguments -- 
notice that abstraction arguments declared in the current world 
are replaced with lambda-terms.  
Again, this is because the terms we produce with these substitution functions 
need to make sense in the parent world, 
where identifiers declared in the current world are out of scope. *)

and argsubs a A (EntArg x) = EntArg (entsubs a A x) |

argsubs x y (AbstArg(t,m)) =

   if m=0 andalso Age (AbstArg(t,m)) (hd(!CONTEXT)) = [0]

   then argsubs x y 
(Lambda(getabstype(hd(Find(AbstArg(t,m))(hd(!CONTEXT))))))

else if x=AbstArg(t,m) andalso x<>y

then argsubs x y y  (*this weird seeming maneuver forced in case
y is expandable *)

else AbstArg(t,m) |


argsubs x y (Lambda T) = 

if x=y then

Lambda (atypesubs x y T)

else  Lambda (atypesubs x y (renamespace T))

(* this is actually matching of types not strict equality:
the first type may be vaguer in having Unknown in the definition field 
at the end where the second is actually defined. *)

(* we need equalentities here, allowing definition expansion (and rewrites) to
establish equality *)

(* the exact argument if true causes us to check 
for equal lambda-terms rather than equal dependent types, 
so we pay attention to the first projection of the last item *)

and equaltypes exact (EType x)  (EType y) = x<> error 
andalso y<> error andalso (x=y orelse equalenttypes x y) |

equaltypes exact (AType(World [(n,a,t)]))  (AType(World [(m,b,u)]))
 = equaltypes false t u   andalso (not exact orelse 
equalentities (deent a) (deent b))   |

equaltypes exact (AType(World ((n,a,t)::L)))(AType(World((m,b,u)::M)))=
if not(equaltypes false t u) then false
else equaltypes exact (typesubs a b (AType(World L)))
(AType(World M)) |

equaltypes exact x y = false

(* equality of entity types *)

and equalenttypes (that P) (that Q) = equalentities P Q |

equalenttypes (IN P) (IN Q) = equalentities P Q |

equalenttypes error x = false |

equalenttypes x error = false |

equalenttypes x y = x=y

(* equality of lambda terms *)

and equivlambdas (Lambda x) (Lambda y) = if x=y then true else equaltypes true (AType x) (AType y) |

equivlambdas x y = x=y

(* equality of entity terms.  This function allows expansion of defined abstractions 
declared anywhere in the context to justify equation of entity terms *)

and equalentities  (App(s,n,M)) (App(t,u,N))

= if App(s,n,M)=App(t,u,N) then true else

let val T = expand(App(s,n,M)) and U = expand (App(t,u,N)) 
and V = rewriteonce  (App(s,n,M)) and W=rewriteonce  (App(t,u,N))in

if V <> App(s,n,M) then equalentities V (App(t,u,N)) 

else if W <> App(t,u,N) then equalentities (App(s,n,M)) W

else if T <> App(s,n,M) then equalentities T (App(t,u,N))

else if U <> App(t,u,N) then equalentities (App(s,n,M)) U

else s=t andalso n=u andalso equalentitieslist M N

end

|

equalentities (App(s,n,M)) x =

if App(s,n,M) = x then true else

let val T = expand(App(s,n,M)) and V = rewriteonce  (App(s,n,M)) in 

if V <> App(s,n,M)

then equalentities V x

else if T <> App(s,n,M)

then equalentities T x

else false end

|  equalentities x (App(s,n,M)) =

if x=App(s,n,M) then true else

let val T = expand(App(s,n,M)) in if T <> App(s,n,M)

then equalentities x T

else false end
 |

equalentities x y = x=y

(* equality of argument lists [name is deceptive] *)

and equalentitieslist nil nil = true |

equalentitieslist ((EntArg a)::L) ((EntArg b)::M) = equalentities a b andalso
  equalentitieslist L M |

equalentitieslist (a::L) (b::M) = (a=b orelse (equivlambdas (expand2 a) (expand2 b)))

andalso equalentitieslist L M |

equalentitieslist x y = false

(* work on type computation *)

(* get entity type from a general type *)

and getenttype (EType x) = x|

getenttype x = (saypause "getenttype error";error)


(* compute entity types *)

and entitytype (Ent(s,0)) = let val S = stringdef s (!CONTEXT)
in if S=nil orelse not(isenttype(pi1(hd S))) then 

(saypause ("Did not find entity "^s^" (entitytype)");Flush();
error)
else getenttype(pi1(hd S))  end |

entitytype (App (s,0,M)) = let val S = stringdef s (!CONTEXT)
in if S=nil orelse isenttype(pi1(hd S)) then 
(saypause ("Did not find abstraction "^s^" (entitytype)");Flush();error)
else typematchcomp (getabstype(pi1(hd S))) M end |

entitytype x = error

(* compute argument types.  Notice the easy case for lambda terms *)

and argtype (EntArg x) = EType(entitytype x) |

argtype (AbstArg (x,n)) = let val S = stringtype x in
   if S=nil orelse isenttype(pi1(hd S)) then 
(saypause ("Did not find abstraction "^x^" (argtype)");
EType error)
   else pi1(hd S) end |

argtype (Lambda x) = AType x

(* the type matching algorithm for dependent types; this
returns the correct type given the full type of the applied
abstraction (inputs and output) and the list of types of the input, or error if matches fail. *)

and typematchcomp (World L) M = typematchcomp0 (deworld(renamespace(World L)))
    (map (fn m=>(0,m,argtype m)) M)

and typematchcomp0 [(n,a,EType t)] nil = etypesubs (EntArg Unknown) (EntArg Unknown) t |

typematchcomp0 ((n,a,t)::L) ((m,A,T)::M) =
     if not(equaltypes false t T) then (saypause ("Type "^(display6 t)^" of "^
(display4 a)^" does not match type "^(display6 T)^" of "^(display4 A));Flush();error)
     else typematchcomp0
     (map (fn(u,v,w) => (u,argsubs a A v,typesubs a A w)) L)
     M |

typematchcomp0 x y = error

and typematchcomp1 [(n,a,t)] nil = typesubs (EntArg Unknown) (EntArg Unknown) t |

(* typematchcomp1 is used to make substitutions into initial
segments of types in FixTypeList below, to handle curried notations for abstraction arguments 10/10 mods *)

typematchcomp1 ((n,a,t)::L) ((m,A,T)::M) =
     if not(equaltypes false t T) then (saypause ("Type "^(display6 t)^" of "
^(display4 a)^" does not match type "^(display6 T)^" of "^(display4 A));Flush();EType error)
     else typematchcomp1
     (map (fn(u,v,w) => (u,argsubs a A v,typesubs a A w)) L)
     M |

typematchcomp1 x y = EType error


(* type matching with no error message reports -- used by
the implicit argument discovery feature *)

and silenttypematchcomp0 [(n,a,EType t)] nil = etypesubs (EntArg Unknown) (EntArg Unknown) t |

silenttypematchcomp0 ((n,a,t)::L) ((m,A,T)::M) =
     if not(equaltypes false t T) then 
((* saypause ("Type "^(display6 t)^" of "^(display4 a)^
" does not match type "^(display6 T)^" of "^(display4 A));Flush(); *) error)
     else silenttypematchcomp0
     (map (fn(u,v,w) => (u,argsubs a A v,typesubs a A w)) L)
     M |

silenttypematchcomp0 x y = error


(* this function computes values of defined abstractions.
It has a parameter which if set to false would do type checking,
but in fact it is only used in safe context so far 
(on things already known to have been type checked). *)

and defmatchcomp safe (World L) M = defmatchcomp0 safe (deworld(renamespace(World L))) 
    (map (if safe then (fn m=>(0,m,EType error)) else (fn m=>(0,m,argtype m))) M)

and defmatchcomp0 safe [(n,EntArg a,EType t)] nil = a |

defmatchcomp0 safe ((n,a,t)::L) ((m,A,T)::M) =
     if (not safe) andalso (not(equaltypes false t T)) then Error
     else defmatchcomp0 safe
     (map (fn(u,v,w) => (u,argsubs a A v,typesubs a A w)) L)
     M |
defmatchcomp0 safe x y = Error

(* one step of expansion of a defined abstraction in applied position *)

and expand (App(s,n,M)) = 

    if n=0 andalso stringAge s = [0]
    then entsubs (EntArg Unknown) (EntArg Unknown) 
    (defmatchcomp true(getabstype(pi1(hd(stringtype s)))) M)
    else App(s,n,M) | 

expand x = x

(* expansion of an argument to a lambda-term if it is defined *)

and expand2 (AbstArg(s,n)) =

if n=0 andalso stringAge(s) = [0]
then argsubs (EntArg Unknown) (EntArg Unknown) (Lambda(getabstype(pi1(hd(stringtype s)))))
else AbstArg(s,n) |

expand2 x=x

(* moving all bound variables in a dependent type or lambda term
to a new namespace, before a substitution into one of these is made *)

and renamespace (World L) = World (renamespace0 (newnameserial())L)

and renamespace0 N [(n,a,t)] = [(n,a,t)] |

renamespace0 N ((m,(AbstArg(s,n)),t)::L)=
(m,(AbstArg(s,N)),t)::
(deworld(atypesubs (AbstArg(s,n)) (AbstArg(s,N)) (World (renamespace0 N L)))) |

renamespace0 N ((m,(EntArg(Ent(s,n))),t)::L)=
(m,(EntArg(Ent(s,N))),t)::
(deworld(atypesubs (EntArg(Ent(s,n))) (EntArg(Ent(s,N))) (World(renamespace0 N L)))) |

renamespace0 N x = (saypause "Bad case in renamespace";nil)

(* matching function to be added at this point.   Two entity terms are matched.
A list of matches for variables is produced *)

(* it is demonstrable that if two type safe expressions, the first of which
is not a variable, match successfully, then they will in fact be of the same type *)

(* solving confluence issues by requiring that executable subterms in the body
be head-rewritten before matching:  an executable can only match anything in a context
in which it has no execution behavior (except at the top of course) *)

and ematch b (Ent(s,n)) t = 

   if n<>0

   then [[(EntArg(Ent(s,n)),EntArg t)]]

   else if (Ent(s,n)) = t 

       then [nil]
       else nil |

ematch b (App(s,n,L)) (App(t,n,M)) =

   let val T = if b then headrewrite (App(t,n,M)) else App(t,n,M) in

   if s <> appof T then nil

   else argmatch  L (argsof T) end |

ematch b x y = nil

and argmatch  nil nil = [nil] |

argmatch  x nil = nil |

argmatch  nil x = nil |

argmatch  ((EntArg x)::L) ((EntArg y)::M) =

   mergematch (ematch true x y) (argmatch  L M) |

argmatch  (AbstArg(s,n)::L) (t::M) =

  if n<>0 then

  mergematch [[(AbstArg(s,n),t)]] (argmatch  L M)

  else if AbstArg(s,n) = t then argmatch  L M else nil |

argmatch  x y = nil

and mergematch nil x = nil |

mergematch x nil = nil |

mergematch [nil] x = x |

mergematch x [nil] = x |

mergematch x y = let val M = mergematch0 (hd x) (hd y) in
   if M = nil then nil else [M] end

and mergematch0 nil L = nil |

mergematch0 L nil = nil |

mergematch0 ((s,t)::L) M = 

   let val N = abstractfind s (L @ M) in

   if N = nil then if L <> nil then (s,t)::(mergematch0 L M)
      else (s,t)::M

   else if equalarguments t (hd N) then

     if abstractdrop s L = nil andalso abstractdrop s M = nil then [(s,t)]

     else if abstractdrop s L <> nil
     then (s,t)::(abstractdrop s (mergematch0 L M))

     else (s,t)::(abstractdrop s M)

   else nil end

and equalarguments (EntArg s) (EntArg t) = equalentities s t |

equalarguments (Lambda x) (Lambda y) = equivlambdas (Lambda x) (Lambda y) |

equalarguments x y = x = y

(* tools for implementing substitutions by pattern matching over
any of the various sorts for which we have substitution functions *)

and listsubsmod subsfun nil T = T |

listsubsmod subsfun ((s,t)::L) T = subsfun s t (listsubsmod subsfun L T)

and matchsubs subsfun pattern target body = 
let val M = ematch false pattern body in if M = nil then body else
listsubsmod subsfun (hd M) target end

(* apply first applicable rewrite rule, just once *)

and rewriteoncewithalist  nil t = t |

rewriteoncewithalist  ((s,(t,u))::L) T =

   let val M = ematch false t T in

       if M = nil then rewriteoncewithalist  L T

       else matchsubs  entsubs t u T end

and rewriteoncewithalistlist  nil T = T |

rewriteoncewithalistlist  (L::M) T =

let val U =rewriteoncewithalist  L T in

if U <> T then U

else rewriteoncewithalistlist M T end

and rewriteonce T = if (!REWRITEVER) 
   then rewriteoncewithalistlist (tl(!REWRITES)) T else T

(* complete rewriting *)

and fullrewrite (App(s,0,L)) =

let val L1 = map fullrewrite2 L in

let val T1 = rewriteonce (App(s,0,L1))

in if T1 = App(s,0,L) then T1

else fullrewrite T1

end end |

fullrewrite (Ent(s,0)) = let val T1 = rewriteonce (Ent(s,0)) in

if T1 = Ent(s,0) then T1 else fullrewrite T1 end |

fullrewrite x = x

and fullrewrite2 (EntArg x) = EntArg (fullrewrite x) |

fullrewrite2 x = x

(* head rewriting (just from the top) *)

and headrewrite (App(s,0,L)) =

let val T1 = rewriteonce  (App(s,0,L))

in if T1 = App(s,0,L) then T1

else headrewrite T1

end |

headrewrite x = x;

(* extract the pattern and target from the list component 
of the type of a function justifying rewrites *)

fun getpattern nil = nil |

getpattern [x] = nil |

getpattern ((i,a,EType(that (App(y,n,[t]))))::x) = [t] |

getpattern (x::L) = getpattern L;

fun gettarget nil = nil |

gettarget ((i,a,EType(that (App(y,n,[t]))))::nil) = [t] |

gettarget [x] = nil |

gettarget (x::L) = gettarget L;


(* this computes the list of variable dependencies of terms of various kinds. *)

fun deps (Ent(s,0)) = if isnew(EntArg(Ent(s,0))) 
then (if isvariable(EntArg(Ent(s,0)))  
then [(EntArg(Ent(s,0)))] else nil)@
(typedeps(argtype(EntArg(Ent(s,0))))) else nil |

deps (App(s,0,L)) = if isnew(AbstArg(s,0))
 then (if isvariable(AbstArg(s,0))   
then [(AbstArg(s,0))]else nil)@
(typedeps(argtype(AbstArg(s,0)))@(depsarg L)) else depsarg L |

deps x = nil 

and depsarg nil = nil |

depsarg((EntArg x)::L) = (deps x)@(depsarg L) |

depsarg ((AbstArg(s,0))::L) = if isnew(AbstArg(s,0)) 
then (if isvariable(AbstArg(s,0)) then [AbstArg(s,0)] else nil)
@(typedeps(argtype(AbstArg(s,0))))@(depsarg L)
   else depsarg L |

depsarg ((Lambda x)::L) = typedeps(AType x) |

depsarg (x::L) = depsarg L

and typedeps (EType(that x)) = deps x |

typedeps (EType(IN x)) = deps x |

typedeps (EType x) = nil |

typedeps (AType(World [(i,EntArg x,T)])) = (deps x)@(typedeps T) |

typedeps (AType(World((i,A,T)::L))) = 
(typedeps T)@(typedeps (AType(World L))) |

typedeps (AType x) = nil;

(* functions for argument redundancy *)

(* a device for type casting an argument to a type, used
in the internals of the main argument reduction functions
because they are doing very general term matching disguised
as type matching, so casting is needed *)

fun arg2type (EntArg x) = EType(that x) |

arg2type (Lambda x) = AType x |

arg2type x = if expand2 x <> x then arg2type(expand2 x) else EType error;

(* moretypes discovers candidate implicit arguments in the types
of explicitly given arguments at declaration time *)

(* it is initially analyzing a type, but it looks into component application
terms, so it is constantly type casting arguments to types, weirdly *)

fun moretypes (EType(that (Ent(s,0)))) = 
if isvariable (EntArg(Ent(s,0))) then [(EntArg(Ent(s,0)),EType prop)] else nil |

moretypes (EType(IN (Ent(s,0)))) = if isvariable (EntArg(Ent(s,0)))
 then [(EntArg(Ent(s,0)),EType TYPE)] else nil |

moretypes (EType(that (App(s,n,(x::L))))) = 
(if isvariable(AbstArg(s,n)) then [(AbstArg(s,n),argtype(AbstArg(s,n)))]
@(moretypes(argtype(AbstArg((s,n))))) else nil)

@(if isvariable x then [(x,argtype x)]@(moretypes(EType(that(App(s,n,L)))))
else (moretypes((arg2type x)))@(moretypes(EType(that(App(s,n,L)))))) |

moretypes (EType(IN (App(s,n,(x::L))))) = moretypes (EType(that (App(s,n,(x::L)))))  |

moretypes (AType(World(nil))) = nil |

moretypes (AType(World([(i,a,t)]))) = (moretypes (EType(that (deent a))))@(moretypes t) |

moretypes (AType(World((i,a,t)::L))) = (moretypes t)@(moretypes(AType(World L)))|

moretypes x = nil;

(* now outline the strategy:  take the first element in the argument list.
compute the expanded alternate list of its tail.  Compute the deps of the item
and drop all dotted and undotted versions of the deps from the previous list.
Then add the item.  Then add dotted versions of its moretypes list, and replace
undotted versions with dotted versions throughout. *)

(* add or remove the initial period (.) which distinguishes
the name of an implicit argument from the name of an explicit argument *)

fun dot s = if s = "" then "" else if hd(explode s) = #"." then s else "."^s;

fun undot s = if s="" then "" 
else if hd(explode s) = #"." then implode(tl(explode s)) else s;

fun argdot (EntArg(Ent(s,n))) = (EntArg(Ent(dot s,n))) |

argdot (AbstArg(s,n)) = AbstArg(dot s,n) |

argdot x = x;

fun argundot (EntArg(Ent(s,n))) = (EntArg(Ent(undot s,n))) |

argundot (AbstArg(s,n)) = AbstArg(undot s,n) |

argundot x = x;

fun drop s nil = nil |

drop s ((i,a,t)::L) = if s=a then drop s L else ((i,a,t)::(drop s L));

fun droplist L nil = nil |

droplist nil L = L |

droplist (s::M) L = drop s(droplist M L);

fun incrementlist nil = nil |

incrementlist ((i,a,t)::L) = (i+1,a,t)::(incrementlist L);

(* make a single substitution in an argument/type list *)

fun singlesubslist s t nil = nil |

singlesubslist s T ((i,a,t)::L) = (1,argsubs s T a,typesubs s T t)::(singlesubslist s T L);

(* this function adds dotted items to an argument list. It is complicated
by the need to ensure that when dotted items are added in a bloc (output of moretypes)
to an argument list that their order is corrected if necessary to keep dependencies sound *)

fun addotlist nil L = L |

addotlist ((s,t)::M) nil = [(1,argdot s,t)] |

addotlist ((s,t)::M) ((i,a,T)::L) = addotlist (moretypes t)((i,argdot s,t)
::(singlesubslist s (argdot s) ((droplist(typedeps t)
(droplist (map argdot (typedeps t)) 
(drop s(drop (argdot s)(addotlist M (incrementlist((i,a,T)::L))))))))));

fun pi1(x,y)=x;

(* this replaces a list of items with their dotted versions throughout an argument/type list *)

fun dotsubslist ((s,t)::L)  nil = nil |

dotsubslist nil L = L |

dotsubslist ((s,t)::L) ((i,a,T)::M) = dotsubslist L 
((i, argsubs s (argdot s) a, typesubs s (argdot s) T)::(dotsubslist ((s,t)::L) M));


(* expand the list of arguments presented for an abstraction
at declaration time
with inferred implicit arguments  *)


fun expandlist nil = nil |

expandlist ((i,EntArg(Ent(s,0)),t)::L) =  addotlist(moretypes t)
(dotsubslist((moretypes t))(droplist(map pi1 (moretypes t))((i,EntArg(Ent(s,0)),t)::
(singlesubslist (EntArg(Ent(dot s,0))) (EntArg(Ent(s,0))) 
(drop (EntArg(Ent(s,0))) (drop (EntArg(Ent(dot s,0))) 
(droplist(typedeps t) (droplist(map argdot (typedeps t))(expandlist L)))))))))

|

expandlist ((i,AbstArg((s,0)),t)::L) =  addotlist(moretypes t)
(dotsubslist((moretypes t))
(droplist(map pi1 (moretypes t))((i,AbstArg((s,0)),t)::
(singlesubslist (AbstArg((dot s,0))) (AbstArg((s,0))) 
(drop (AbstArg((s,0))) (drop (AbstArg((dot s,0))) 
(droplist(typedeps t) 
(droplist(map argdot (typedeps t))(expandlist L)))))))))

|

expandlist ((i,a,t)::L) =  addotlist(moretypes t)
(dotsubslist((moretypes t))(droplist(map pi1 (moretypes t))((i,a,t)::
((* drop (EntArg(Ent(s,0))) *) (droplist(typedeps t) 
(droplist(map argdot (typedeps t))(expandlist L)))))))
;

fun guardedexpandlist L = if (!IMPLICITVER) then (expandlist L) else L;

(* functions to repair an argument list, adding values for implicit arguments in the type *)

fun firstundotted nil = EType error |

firstundotted ((i,EntArg(Ent(s,n)),t)::L) =

if s<> "" andalso hd(explode s) = #"." then firstundotted L

else t |

firstundotted ((i,AbstArg((s,n)),t)::L) =

if s<> "" andalso hd(explode s) = #"." then firstundotted L

else t |

firstundotted ((i,a,t)::L) = t;

fun deworld2 (AType(World L)) = L |

deworld2 x = nil;

(* This function finds the value of a given implicit argument
by type matching.  Since it is actually looking deep into the structure
of ordinary terms embedded in types, including ordinary terms with
variable binding, it is doing a lot of weird type casting, sometimes
quite incompatible with Lestrade's own type system :-) 

In principle, this search can fail, if an argument is implicitly typed
using a subterm of a type which can actually be eliminated by a definitional
expansion.  I have not seen this happen.  It would not represent a failure
of the logic:  the implicit argument feature actually touches the logic
not at all.  

It can fail quite easily for an implicit abstraction argument.

*)

fun initialsegment nil L = true |

initialsegment (x::L) (y::M) = if x<>y then false else initialsegment L M |

initialsegment L M = false;


fun matchsegment nil L = nil |

matchsegment (x::L) (y::M) =

   (y::(matchsegment L M)) |

matchsegment L M = M;

fun findimplicitarg Types1 Types a atype 
(EType(that(App(s,n,(x::L))))) (EType(that(App(t,m,(y::M))))) =

if a = AbstArg(s,n) andalso initialsegment (x::L) (map (fn (x1,y1,z1) => y1) Types1)
andalso silenttypematchcomp0 (deworld2 atype) (matchsegment (x::L) Types) <> error then

if (y::M) = (map (fn (x1,y1,z1) => y1) (matchsegment (x::L) Types)) then AbstArg(t,m) else

Lambda((World((matchsegment (x::L) Types)@[(0,EntArg(App(t,m,(y::M))),
EType (silenttypematchcomp0 (deworld2 atype) (matchsegment (x::L)Types)) )])))

else if a = AbstArg(s,n) andalso 
entitytype (App(t,m,y::M)) <> error andalso
equalenttypes (silenttypematchcomp0 (deworld2 atype) 
(map (fn xx => (1,xx,argtype xx)) (y::M))) 
(entitytype (App(t,m,y::M)))  then AbstArg(t,m)

else 

if s<>t orelse m<>n then 

let val T = expand (App(s,n,(x::L))) and U = expand (App(t,m,(y::M))) in

if T <> (App(s,n,(x::L))) 
   then 

   let val V = findimplicitarg Types1 Types a atype 
   (EType(that(T))) (EType(that(App(t,m,(y::M)))))

   in if V <> EntArg Error 

     then V

     else if U = (App(t,m,(y::M)))  then EntArg Error

         else findimplicitarg Types1 Types a atype 
         (EType(that(App(s,n,(x::L))))) (EType(that(U))) end 
else if U = (App(t,m,(y::M)))  then EntArg Error

         else findimplicitarg Types1 Types a atype 
         (EType(that(App(s,n,(x::L))))) (EType(that( U))) end

else if a = x then y

else let val T = findimplicitarg Types1 Types a atype (arg2type x) (arg2type y) in

if T <> EntArg Error then T

else findimplicitarg Types1 Types a atype
 (EType(that(App(s,n,L)))) (EType(that(App(t,m,M)))) end |

findimplicitarg Types1 Types a atype (EType(that(App(s,n,(L))))) (EType(that(T))) =

if a = AbstArg(s,n)

   andalso initialsegment L (map (fn (x1,y1,z1) => y1) Types1)

       andalso silenttypematchcomp0 (deworld2 atype) (matchsegment L Types) <> error

       then Lambda(World((matchsegment L Types)@
       [(0,EntArg T, EType (silenttypematchcomp0 (deworld2 atype)
       (matchsegment L Types)))]))

   else if expand (App(s,n,L)) <> App(s,n,L) then

findimplicitarg Types1 Types a atype 
(EType(that(expand(App(s,n,(L)))))) (EType(that(T))) 

 else if expand T <> T then

findimplicitarg Types1 Types a atype 
(EType(that(App(s,n,L)))) (EType(that(expand(T))))

else EntArg Error |

findimplicitarg Types1 Types a atype 
(EType(IN(App(s,n,(L))))) (EType(IN(T))) =
findimplicitarg Types1 Types a atype 
(EType(that(App(s,n,(L))))) (EType(that(T))) |

findimplicitarg Types1 Types a atype (EType(that b)) (EType(that c)) = 
if a=EntArg b then EntArg c else EntArg Error|

findimplicitarg Types1 Types a atype (EType(IN b)) (EType(IN c)) = 
if a=EntArg b then EntArg c else EntArg Error|

findimplicitarg Types1 Types a atype (AType(World([(i,b,t)]))) (AType(World([(j,c,u)]))) =

let val T = findimplicitarg Types1 Types a atype t u in

if T = EntArg Error then findimplicitarg Types1 Types a atype 
   (EType (that (deent b))) (EType (that (deent c))) 
else T end |

findimplicitarg Types1 Types a atype (AType(World((i,b,t)::L))) (AType(World((j,c,u)::M))) =

let val T = findimplicitarg Types1 Types
 a atype t u in if T = EntArg Error 
then findimplicitarg  (Types1@[(i,b,t)]) (Types@[(j,c,u)]) 
a atype (AType (World L)) (AType (World M)) else T end|

findimplicitarg Types1 Types a atype x y = EntArg Error;

(* this function repairs the argument list supplied to a function with
implicit arguments at parse time.  It is important to notice that implicit arguments
play no role in the logic at all! *)

fun fixarglist nil x = nil |

fixarglist x nil = nil |

fixarglist ((i,EntArg(Ent(s,n)),t)::L)  ((T)::M) =

if s = "" orelse hd(explode s) <> #"." then (T)::(fixarglist (singlesubslist (EntArg(Ent(s,n))) T L) M)

else 

let val MM = (findimplicitarg nil nil (EntArg(Ent(s,n))) 
(typesubs (EntArg Unknown) (EntArg Unknown) t) (firstundotted L) (argtype T)) in

(MM::(fixarglist (singlesubslist (EntArg(Ent(s,n))) MM L)((T)::M))) end |

fixarglist ((i,AbstArg((s,n)),t)::L)  ((T)::M) =

if s="" orelse hd(explode s) <> #"." then (T)::(fixarglist (singlesubslist (AbstArg(s,n)) T L) M)

else 

let val MM = findimplicitarg nil nil ((AbstArg(s,n))) 
(typesubs (EntArg Unknown) (EntArg Unknown) t)  (firstundotted L) (argtype T) in

(MM::(fixarglist (singlesubslist (AbstArg(s,n)) MM L) ((T)::M))) end |

fixarglist (x::L) (y::M) = y::(fixarglist L M);

fun guardedfixarglist L M = if (!IMPLICITVER) then fixarglist L M else M;

(* The following block of functions is used in the quite elaborate
check of the structure of inputs to the rewritec and rewrited commands *)

fun lastthree (x::y::z::w::L) = lastthree (y::z::w::L) |

lastthree L = L;

fun allbutlastthree nil = nil |

allbutlastthree [x] = nil |

allbutlastthree [x,y] = nil |

allbutlastthree [x,y,z] = nil |

allbutlastthree (x::L) = x::(allbutlastthree L);

(* variables cannot be patterns, even if old.  Constants can. *)

fun notvararg (EntArg(Ent(s,0))) = false |

notvararg (EntArg x) = true |

notvararg x = false;

fun inlist x nil = false |

inlist x (s::L) = if x=s then true else inlist x L;

fun allinlist nil L = true |

allinlist (s::M) L = inlist s L andalso allinlist M L;

(* check correctness of argument lists for rewriting commands *)

(* a good rewrite list has at least three elements in it *)

fun goodrewritelist nil = false |

goodrewritelist [x] = false |

goodrewritelist [x,y] = false |

goodrewritelist L = 

   let val [P,Q,R] = lastthree(L) in let val T = argtype Q in

   (* Q is an entity argument and not polymorphic *)

   notvararg Q andalso typerigid(deent Q) andalso

   (* the last two items have the same type *)

   equaltypes false T (argtype R) andalso

   (* the first of the last three items 
    is a predicate variable over this type *)

   isvariable P andalso

   equaltypes true (argtype P) 
   (AType(World([(~1,EntArg(Ent("???",1)),T),
   (~1,EntArg(Unknown),EType prop)])))

   andalso

   (* P does not appear in the deps of Q *)

   not (inlist P (depsarg [Q]))

   (* everything in allbutlastthree L appears in the deps of P *)

   andalso allinlist (allbutlastthree L) (depsarg [Q])

   (* everything in the deps of R appears in the deps of Q *)

   andalso allinlist (depsarg [R]) (depsarg [Q]) end end;

(* user command:  close the last world opened; 
 this sends an error message if it attempts to close world 1 *)

(* Close does not automatically save the current world:  that has to be done with Save() *)

fun Close() = (if length(!CONTEXT) > 2 then
(CONTEXT:= (tl (!CONTEXT));REWRITES:=tl(!REWRITES);
WORLDNAMES:= tl(!WORLDNAMES))  
else saypause ("Cannot undo move 1:"^(hd (!WORLDNAMES))));

fun savelist worldnames context = 

if length context <= 1 orelse length worldnames <= 1 then nil
   else (worldnames,hd context)::(savelist (tl worldnames) (tl context));

(* save all worlds on path to current world *)
(* it is not possible to save a world which has its default numerical name *)

fun Save s = if s = makestring(length(!CONTEXT)-1) 
then saypause "Cannot save a move with the default numeral name"

else

if defaultworld (tl(!WORLDNAMES)) 
then saypause "Cannot save a default move" else

(SAVEDWORLDS := abstractmerge 
(savelist (s::(tl(!WORLDNAMES))) (!CONTEXT))(!SAVEDWORLDS);
SAVEDREWRITES := abstractmerge 
(savelist (s::(tl(!WORLDNAMES))) (!REWRITES))(!SAVEDREWRITES);
WORLDNAMES:=s::(tl(!WORLDNAMES)));

fun ClearCurrent s = if defaultworld (tl(!WORLDNAMES))
 andalso not(defaultworld(s::(tl(!WORLDNAMES))))

then saypause "Named move cannot follow a default move" else

let val W = abstractfind (s::(tl(!WORLDNAMES))) (!SAVEDWORLDS) 
and R = abstractfind (s::(tl(!WORLDNAMES))) (!SAVEDREWRITES)
in

(if W = nil then

(CONTEXT:=(World nil)::(tl (!CONTEXT)); REWRITES:=nil
::(tl(!REWRITES));WORLDNAMES:=s::(tl(!WORLDNAMES)))

else (CONTEXT:=(hd W)::(tl (!CONTEXT)); 
REWRITES:=(hd R)::(tl (!REWRITES)); 
WORLDNAMES:=s::(tl(!WORLDNAMES)))) end;

(* serial number used for recording age of declarations *)

val SERIAL = ref 0

(* completely clear the Lestrade context, user command, also
issued by readfile *)

fun ClearAll() = (GREETED:=false;CONTEXT:=[World nil,World nil];
REWRITES:=[nil,nil];SERIAL:=0;NAMESERIAL:=0;WORLDNAMES:=["1","0"]);

(* load a named theory.  
This completely clears the context, supplying the world 0
declarations of the saved theory -- 
without rewrite decs, perhaps I should fix this.  *)

fun LoadTheory s = let val S = abstractfind s (!SAVEDTHEORIES) in

if s="" orelse S = nil then saypause "No such theory to load"

else let val (N1,N2,W) = hd S in

(GREETED:=false;CONTEXT:=[World nil,W];
REWRITES:=[nil,nil];SERIAL:=N1;NAMESERIAL:=N2;WORLDNAMES:=["1","0"])

end end;

(* type check an entity type.  
This checks that P in that P is a prop and T in in T is a type *)

fun typecheck obj = true |

typecheck prop = true |

typecheck TYPE = true |

typecheck (that P) =

let val ANSWER = (entitytype P = prop) in
(if not ANSWER then 
saypause ((display2 P)^" is not of type prop (typecheck)")else();Flush();ANSWER) end|

typecheck (IN P) =

let val ANSWER = (entitytype P = TYPE) in
(if not ANSWER then saypause 
((display2 P)^" is not of type 'type' (typecheck)")else();Flush();ANSWER) end|

typecheck error = false;


(* the entity declaration command will take as arguments a string
s and an EntType  T.  It needs to check that s is not already declared, 
then check that T declaration checks, 
then add (EntArg(Ent(s)),T) as an entry to the first world in the context *)

(* actually it needs to do a full type check that T is of type prop *)

(* utility adds an item to a world at the end *)

fun addtoworld (World L) x 
= (Reset();World(L @ [(fn (i,A,T)=>(i,Reindex5 A,Reindex3 T)) x]));

fun addtoworld0 (World L) x = World(L @ [x]);

(* increment the declaration age counter *)

fun newserial() = (SERIAL:=1+(!SERIAL);(!SERIAL));

(* reserved identifiers. *)

fun reserved s = s="obj" orelse s="prop" orelse s="that" 
orelse s="type" orelse s="in" orelse s = "---" orelse s = "???";

(* command for postulating an entity in the current world *)

fun Declare s T =

   if reserved s orelse stringdef s (!CONTEXT) <> nil 
       then saypause ("Identifier "^s^" is not fresh")
   else if not (typecheck T) then saypause "Type check fails"
   else (CONTEXT := (addtoworld0 (hd(!CONTEXT))(newserial(),
      EntArg(Ent (s,0)),EType T))::(tl(!CONTEXT));showdec s);

(* command for postulating a construction in the parent world *)

(* s is a name.  L could also be a list of names.  T is an entity type. *)

(* the list of names of identifiers is in order of age;
this ensures sensible dependencies without evilly recursive
checks *)

(* modified to support the notion that definitions have
age 0 *)

fun isordered nil = true |

isordered [a] = true |

isordered (a::(b::L)) = 

hd(Age a (hd(!CONTEXT))) <> 0 andalso

(hd(Age a (hd(!CONTEXT))) < hd(Age(b)(hd(!CONTEXT))) 
   andalso isordered (b::L));

(* toolkit for turning argument lists into worlds -- the
type of an abstraction is a little world *)

(* there are some functions here for interaction with the implicit arguments feature *)

fun dotfix nil t = t |

dotfix ((i,a,t)::L) t = if (!IMPLICITVER) 
andalso argundot a <> a then etypesubs (argundot a) a (dotfix L t) else dotfix L t;

fun dotfix2 nil t = t |

dotfix2 ((i,a,t)::L) t = if (!IMPLICITVER) 
andalso argundot a <> a then entsubs (argundot a) a (dotfix2 L t) else dotfix2 L t;

(* remove dotted items from a list; used by the parser *)

fun dotpurge nil = nil |

dotpurge ((i,EntArg(Ent(s,n)),t)::L) = if (!IMPLICITVER) then

   if s <> "" andalso hd(explode s) = #"." then dotpurge L
   else ((i,EntArg(Ent(s,n)),t)::(dotpurge L)) 
   else ((i,EntArg(Ent(s,n)),t)::L) |

dotpurge((i,AbstArg(s,n),t)::L) = if (!IMPLICITVER) then

   if s <> "" andalso hd(explode s) = #"." then dotpurge L
   else ((i,AbstArg(s,n),t)::(dotpurge L)) else 
  ((i,AbstArg(s,n),t)::L) |

dotpurge (x::L) = x::(dotpurge L);

fun worlditem s = (hd(Age s (hd(!CONTEXT))),s,
typesubs (EntArg Unknown) (EntArg Unknown) 
(hd(Find s (hd(!CONTEXT)))));

fun worldof L = World(guardedexpandlist(map worlditem L));

(* the current world is temporarily replaced 
during the declaration process for abstractions -- this is a place to keep it *)

val SAVECONTEXT = ref (hd(!CONTEXT));

val SAVECONTEXT2 = ref (hd(!CONTEXT));

(* in construct, the user supplies as an argument list
all notions in the current world on which the construction
depends, in order of construction.
We do a dynamic maneuver:  replace the current world with the 
part of the current world indicated by the argument list;
declaration check this world using the resulting context 
(checking that it includes all its own needed dependencies)
and type check the entity type argument in this context; then restore the context. *)

fun Construct s L T = 

if not(testall isvariable L) then saypause "Some argument is not variable"

(* testing for order is a cute way to enforce sensible dependencies;
the implicit arguments feature now does the laborious checks for this, but
if it is turned off this condition handles deps just fine *)

else if not (isordered L) then saypause "Arguments are in the wrong order"

else let val T = dotfix (deworld(worldof L)) T in

(

if reserved s orelse stringtype s <> nil 
then saypause ("Identifier "^s^" is not fresh")

else if L = nil then 
    let val TT = etypesubs (EntArg Unknown) (EntArg Unknown) T in

    (SAVECONTEXT:=hd(!CONTEXT); 
    CONTEXT:=(World nil)::tl(!CONTEXT);if not (typecheck TT) 

    then (saypause "Type check fails in declaration of constant";
    CONTEXT:=(!SAVECONTEXT)::(!CONTEXT))
   
else (CONTEXT := (!SAVECONTEXT)::
(addtoworld0 (hd(tl(!CONTEXT)))(newserial(),EntArg(Ent (s,0)),EType TT))
::(tl(tl(!CONTEXT)));showdec s)) end

else if not 

let val TT =  (etypesubs (EntArg Unknown) (EntArg Unknown) T) in

(
SAVECONTEXT:=(hd(!CONTEXT)); 

CONTEXT:= (worldof (*worldof2*) L)::(tl(!CONTEXT));

let val CHECK = deccheck4 (!CONTEXT) (hd(!CONTEXT)) 
andalso typecheck (TT) in (CONTEXT:=(!SAVECONTEXT)::(tl(!CONTEXT));CHECK) 

end
) end then saypause "Dependency or type check failure"

else (

let val newparentcontext = addtoworld0 (hd(tl(!CONTEXT)))
(newserial(),AbstArg(s,0),(Reset();
Reindex3(AType(renamespace(addtoworld0 (worldof L) 
(0,EntArg Unknown,typesubs (EntArg Unknown)  (EntArg Unknown) (EType T))))))) in

(CONTEXT:= (hd(!CONTEXT))::newparentcontext::(tl(tl(!CONTEXT)));showdec s) end)

) end;

fun Define s L T =

if not(testall isvariable L) then saypause "Some argument is not variable"

(* same remark on the argument order test as above *)

else if not (isordered L) then saypause "Arguments are in the wrong order"

else 

let val T0 = T and T = dotfix2 (deworld(worldof L)) T in

if reserved s orelse stringtype s <> nil 
then saypause ("Identifier "^s^" is not fresh")

else if not let val T2 = ( (entsubs (EntArg Unknown) (EntArg Unknown) T))
and THETYPE = dotfix (deworld(worldof L)) 
(etypesubs (EntArg Unknown) (EntArg Unknown) (entitytype T0)) in (
SAVECONTEXT:=(hd(!CONTEXT)); 

CONTEXT:= (worldof (* worldof2 *) L)::(tl(!CONTEXT));  

let val CHECK = deccheck4 (!CONTEXT) (hd(!CONTEXT)) andalso 
deccheck1 (!CONTEXT) THETYPE andalso deccheck2 (!CONTEXT) T2

in (CONTEXT:=(!SAVECONTEXT)::(tl(!CONTEXT));CHECK) 

end
) end then saypause "Type check or dependency failure"
else let val TT = fullrewrite(entsubs (EntArg Unknown) (EntArg Unknown) (T))
 and  THETYPE = dotfix(deworld(worldof L))
(etypesubs (EntArg Unknown) (EntArg Unknown) (entitytype (T0)) ) in (
let val newparentcontext = addtoworld0 (hd(tl(!CONTEXT)))
(0,AbstArg(s,0),(Reset();Reindex3(
AType(renamespace(addtoworld0 (worldof L) (0,EntArg TT,EType (THETYPE))))))) in

(CONTEXT:= (hd(!CONTEXT))::newparentcontext::(tl(tl(!CONTEXT)));showdec s) 

end)



end end;

fun deabst(AbstArg(s,0)) = s |

deabst x = "?!?!";

fun deworld2 (AType x) = deworld x;

(* construct a function witnessing validity of a rewrite rule *)

fun Rewritec s L V =

if not(!REWRITEVER) 
then saypause "Rewriting is turned off" else

if reserved s orelse stringtype s <> nil 
then saypause ("Identifier "^s^" is not fresh")

else if reserved (V) orelse stringtype (V) <> nil 
then saypause ("Identifier "^(V)^" is not fresh")

else if not (goodrewritelist L) 
then saypause "Proposed rewrite list does not sort check"

else let val [P,Q,R] = lastthree L and L1 = allbutlastthree L in

(
Declare (V) (that (App(deabst P,0,[Q])));
Construct s (L1 @ [P,EntArg(Ent(V,0))]) (that (App(deabst P,0,[R])));

let val T = stringtype s in if T = nil 
then saypause ("Construction of "^s^" failed for some reason")

else let val (L1::L2::L3) = (!REWRITES) in

REWRITES:= (L1::((s,(Negindex4((deent(hd(getpattern (deworld2(pi1(hd T))))))),
(Negindex4(deent(hd(gettarget (deworld2(pi1(hd T))))))))) :: L2)::L3)

end

end

)
end;

(* show that the function named by s
 witnesses the validity of proposed rewrite rule *)

fun Rewrited s L V =

if not(!REWRITEVER) then say "Rewriting is turned off" else

if stringtype s = nil 
then saypause ("Evidence function "^s^" is not declared")

else if reserved (V) orelse stringtype (V) <> nil 
then saypause ("Identifier "^(V)^" is not fresh")

else if not (goodrewritelist L) 
then saypause "Proposed rewrite list does not sort check"

else let val [P,Q,R] = lastthree L and L1 = allbutlastthree L in

(Declare (V) (that (App(deabst P,0,[Q])));
if equaltypes false (EType((entitytype(App(s,0,
guardedfixarglist (deworld(getabstype(pi1(hd(stringtype s)))))
(L1@[P,EntArg(Ent(V,0))])))))) 
(EType(that(App(deabst P,0,[R]))))

then (say ("Rewrite demonstration succeeded")  ;
let val T = stringtype s in if T = nil 
then saypause ("This error message should never occur")

else let val (L1::L2::L3) = (!REWRITES) in

REWRITES:= (L1::((s,((deent(hd(getpattern (deworld2(pi1(hd T)))))),
(deent(hd(gettarget (deworld2(pi1(hd T)))))))) :: L2)::L3)

end end )


else saypause ("Rewrite demonstration failed")) end;

(* parser *)

(* this was originally Polish notation with a following comma to signal that
an abstraction appears as an argument; 
it was then upgraded to suppport use of
 abstractions of arity greater than 1 
as infix or mixfix operators, and commas 
are allowed between any arguments, and mandatory before and after 
abstraction identifier arguments 
to avoid confusion with abstractions in applied or infix/mixfix position. *)

fun islower c = #"a" <= c andalso c <= #"z";

fun isupper c = #"A" <= c andalso c <= #"Z";

fun isnumeral c = #"0" <= c andalso c <= #"9"

fun isspecial  c = c= #"~" 

orelse c = #"@" orelse c = #"#" orelse c = #"$" 
orelse c = #"%" orelse c = #"^" orelse c = #"&"
 orelse c = #"*" orelse c = #"-" orelse c = #"+"
 orelse c = #"=" orelse c = #"|" orelse c = #";" orelse c = #"." orelse c = #"<" 
orelse c = #">" orelse c = #"?" orelse c = #"/"
 orelse c = #"!" orelse c = #".";

(* get first identifier from a list of characters *)

fun getident nil = nil |

getident (#"\"" :: L) = L |

getident [c] = if islower c orelse isupper c orelse isnumeral c 
orelse isspecial c orelse c = #"," orelse c= #":" 
orelse c = #"(" orelse c = #")" then [c] else nil |

(* I could fiddle with allowed shapes of identifiers here *)

getident (a::(b::L)) = 

if a = #"," orelse a = #":" orelse a = #"(" orelse a = #")" then [a] else

if a = #" " then getident (b::L)
     else if isupper a

         then if islower b orelse isnumeral b

         then a::(getident(b::L))

         else [a] 

     else if islower a

         then if islower b orelse isnumeral b

         then a::(getident(b::L))

         else [a]

     else if isnumeral a

         then if isnumeral b

         then a::(getident(b::L))

         else [a]

     else if isspecial a

         then if isspecial b

         then a::(getident(b::L))

         else [a]

     else nil;

(* the rest of the stream of characters after the first identifier is read *)

fun restident nil = nil |

restident (#"\"" :: L) = nil |

restident [c] = nil |

restident (a::(b::L)) = 

if a = #"," orelse a = #":" orelse a = #")" orelse a = #"(" then b::L else

if a = #" " then restident (b::L)
     else if isupper a

         then if islower b orelse isnumeral b

         then restident(b::L)

         else (b::L)

     else if islower a

         then if islower b orelse isnumeral b

         then restident (b::L)

         else (b::L)

     else if isnumeral a

         then if isnumeral b

         then restident(b::L)

         else b::L

     else if isspecial a

         then if isspecial b

         then restident(b::L)

         else b::L

     else nil;

(* utility for tokenization *)

fun testidentlist nil = nil |

testidentlist (#">":: #">"::L) = ">> "::[implode L] |

testidentlist L = (implode(getident(L))::(testidentlist (restident L)));

(* get a list of tokens (identifiers and punctuation) from a string *)

fun tokenize s = testidentlist(explode s);

(* repair an application term with missing arguments into a 
lambda-term 10/10 mods *)

fun pi13(x,y,z) = x;

fun pi23(x,y,z) = y;

fun pi33(x,y,z) = z;

(* convert type of a primitive construction into a lambda term *)

fun lambdaform s L = if pi23 (hd(rev L)) = EntArg Unknown
then rev((pi13(hd(rev L)),
EntArg(App(s,0,map pi23 (rev(tl(rev L))))),pi33(hd(rev L)))::(tl(rev L))) else L;

(* construct correct argument list for a curried abstraction argument *)

fun FixListType final L T1 =

if length T1 = length L + 1 then [(pi13(hd (rev T1)), 
if final then EntArg(defmatchcomp0 true T1 L) else pi23(hd(rev T1)), typematchcomp1 T1 L)]

else let val LL = FixListType false L (rev(tl(rev T1))) in

LL @ (FixListType final (L@LL) T1) end;

(* transform an argument which has explicitly closed argument
list into a curried abstraction if appropriate 10/10 mods *)

fun FixApp(EntArg(App(a,0,L))) =

let val T1 = stringtype a in

if length L = length(deworld(getabstype(pi1(hd(T1)))))-1

then EntArg(App(a,0,L))

(* else if stringAge a <> [0] then 
(saypause "Cannot curry a primitive construction";EntArg Error) *)

else Lambda (renamespace(World(FixListType true 
(map (fn m =>(0,m,argtype m)) L) (lambdaform a(deworld(getabstype(pi1(hd(T1))))))))) end |

FixApp t = t;

(* get a non-infix term from a list of tokens *)

fun getterm nil = EntArg Error |

getterm (a::L) =

if reserved a then EntArg Error else

if a = "(" then let val TERM = getterms L and REST = restterms L in

if REST<>nil andalso hd REST = ")" then TERM else EntArg Error

end 

else if a = "," then getterm L else

let val T1 = stringtype a in

if T1 = nil then EntArg Error

else if isenttype (pi1(hd T1)) then EntArg(Ent(a,0))

else if length(dotpurge(deworld(getabstype(pi1(hd(T1)))))) = 1

then EntArg(App(a,0,nil))

else if L = nil orelse hd L = "," orelse hd L = ":" 
orelse hd L = ")" orelse reserved (hd L) then AbstArg(a,0)

else if length(dotpurge(deworld(getabstype(pi1(hd T1))))) = 2

then EntArg(App(a,0,
guardedfixarglist (deworld(getabstype(pi1(hd T1))))[getterm L]))

(* adding possibility of application terms with missing
arguments representing abstractions 10/10 mods -- argument
list enclosed by parentheses can be of variable length *)

else if hd L = "(" then let val TERM = 
    guardedfixarglist (deworld(getabstype(pi1(hd T1)))) 
    (getopenarglist (tl L))
    and REST = 
    restopenarglist (tl L) 
    in  if REST <> nil andalso hd REST = ")" 
    then FixApp(EntArg(App(a,0,TERM))) else EntArg Error end


else EntArg(App(a,0,guardedfixarglist (deworld(getabstype(pi1(hd T1))))
(getarglist (length(dotpurge(deworld(getabstype(pi1(hd(T1))))))-1)  (L))))

end

(* the rest of the list of tokens after the first non-infix term is read *)

and restterm nil = nil |

restterm (a::L) = 

if reserved a then (a::L) else

if a = "(" then let val REST = restterms L in

if REST <> nil andalso hd REST = ")" then tl REST else nil

end

else if a = "," then restterm L else

let val T1 = stringtype a in

if T1 = nil then a::L

else if isenttype (pi1(hd T1)) then L

else if length(dotpurge(deworld(getabstype(pi1(hd(T1)))))) = 1

then L

else if L = nil then nil

else if hd L = ","  then L

else if reserved(hd L) then L

else if length(dotpurge(deworld(getabstype(pi1(hd T1))))) = 2

then restterm L


else if hd L = "(" then let val REST = 
    restarglist (length(dotpurge(
    deworld(getabstype(pi1(hd(T1))))))-1)  (tl L)
 in  if REST <> nil andalso hd REST = ")" then tl REST else nil end

else restarglist (length(dotpurge(deworld(getabstype(pi1(hd(T1))))))-1)  (L)

end

(* get an infix term from a stream of tokens *)

and getterms L = 

if L = nil orelse reserved(hd L) then EntArg Error else

let val TERM = getterm L and REST = restterm L in

if REST = nil orelse hd REST = "," orelse hd REST = ")" 
orelse hd REST = "(" orelse hd REST = ":" orelse reserved (hd REST) then TERM

else let val T1 = stringtype (hd REST)

in if T1=nil orelse isenttype (pi1(hd T1)) then TERM

else if length(dotpurge(deworld(getabstype(pi1(hd(T1)))))) <=2 then TERM

else if tl REST = nil orelse hd(tl REST) = "," 
orelse hd(tl REST) = ":" orelse hd(tl REST) = ")" 
orelse reserved(hd(tl(REST))) then TERM

else EntArg(App(hd REST,0,guardedfixarglist 
(deworld(getabstype(pi1(hd T1))))((TERM)::
(getarglist 
(length(dotpurge(deworld(getabstype(pi1(hd(T1))))))-2)  (tl REST)))))

end end

(* the rest of the stream of tokens after the first infix term is read *)

and restterms L =  

if L = nil orelse reserved(hd L) then L else

let val REST = restterm L in

if REST = nil orelse hd REST = "," orelse hd REST = ")" 
orelse hd REST = "(" orelse hd REST = ":" orelse reserved (hd REST) then REST

else let val T1 = stringtype (hd REST)

in if T1=nil orelse isenttype (pi1(hd T1)) then REST

else if length(dotpurge(deworld(getabstype((pi1(hd(T1))))))) <=2 then REST

else if tl REST = nil orelse hd(tl REST) = "," 
orelse hd(tl REST) = ":" orelse hd(tl REST) = ")" orelse reserved(hd(tl REST)) then REST


else restarglist (length(dotpurge(deworld(getabstype(pi1(hd(T1))))))-2)  (tl REST)

end end

(* a list of arguments of known length,
 without enclosing parentheses (these are handled by getterm *)

and getarglist 0 L = nil |

getarglist n nil = [EntArg Error] |

getarglist n L = (getterms L)::(getarglist (n-1) (restterms L))

(* what is left after reading a list of arguments
 of known length without enclosing parentheses *)

and restarglist 0 L = L |

restarglist n nil = nil |

restarglist n L = restarglist (n-1) (restterms L)

(* get a list of arguments of unknown length 
(as in an abstraction declaration or (10/10) a prefix term with
explicit argument list) *)

and getopenarglist nil = nil |

getopenarglist ((":")::L) = nil |

getopenarglist ((")")::L) = nil |


getopenarglist L = 

if restterms L = L then nil else

(getterms L)::(getopenarglist(restterms L))

(* what is left after reading a list of arguments of unknown length *)

and restopenarglist nil = nil |

restopenarglist ((":")::L)  = L |

restopenarglist ((")")::L)  = (")"::L) |


restopenarglist L = 

if restterm L = L then L else
restopenarglist(restterms L);


fun readenttype nil = error |

readenttype (a::L) = if a = "obj" then obj

else if a="prop" then prop

else if a="that" then let val P = deent(getterms L )in 

if P = Unknown then error else

that P end

else if a="type" then TYPE

else if a="in" then let val P = deent(getterms L )in 

if P = Unknown then error else

IN P end else error;

fun restenttype nil = nil |

restenttype (a::L) = if a = "obj" then L

else if a="prop" then L

else if a="that" then let val P = deent(getterms L )in 

if P = Unknown then (a::L) else

restterms L end

else if a="type" then L

else if a="in" then let val P = deent(getterms L )in 

if P = Unknown then (a::L) else

restterms L end else (a::L)


(* the command line just read *)

val THELINE = ref "";
val THELINE2 = ref "";

(* the file from which commands are being read, used by readfile
in indented and unindented versions *)

val READFILE = ref (TextIO.openIn("default"));

val LOGNAME = ref "";

val LOGNAME2 = ref "";

fun Hd nil = "" |

Hd x = hd x;

fun Tl nil = nil |

Tl x = tl x;

(* test functions *)

fun sarg s = getterms(tokenize s);

fun sent s = deent(sarg s);

fun stype s = readenttype(tokenize s);

fun sType s = EType(stype s);

fun slist s = getopenarglist (tokenize s);

fun slist2 n s = getarglist n(tokenize s);

fun Sent s = (say (display2(sent s)); say (display1(entitytype(sent s))));

fun Stype s = (say (display1(stype s)));

fun Moretypes s = say (display5(World(map (fn (x,y) => (1,x,y))
(moretypes (argtype(getterm (tokenize s)))))));

fun Expandlist s = say(display5(World(expandlist (map (fn x => (1,x,argtype x))
(getopenarglist(tokenize s))))));

fun Fixarglist s t = say(display5(World(map (fn x => (1,x,argtype x))(fixarglist

(deworld(getabstype(pi1(hd(stringtype s)))))

((getopenarglist (tokenize t)))

)))); 

val READFILEDEPTH = ref 0;

(* read a command line from a stream of tokens *)

val THEORYNAME = ref "bogus";

fun readline nil = () |

readline (a::L) =

if a = "readfile" andalso length(L)>=2 
then (Flush();READFILEDEPTH:=1+(!READFILEDEPTH);readfile (hd L) (hd(tl L)))

else if a = "parsetest" andalso L<>nil then Sent (hd L)

else if a = "parsetest2" andalso L<> nil then Stype(hd L)

else if a = "declare" andalso L<> nil andalso tl L<>nil 
then let val s = (hd L) and t = readenttype (tl L) in 

(if restenttype (tl L)<> nil andalso restenttype (tl L)<> [""] 
then saypause ("Declaration line not completely read:  "^(hd L)) else ();
TextIO.output(!LOGFILE,(!THELINE));
say0("Your command: "^(!THELINE2)); Declare s t) end

else if a = "construct" andalso L<>nil andalso tl L <> nil 
then let val s = hd L and L1 = getopenarglist (tl L) 
and T = readenttype (restopenarglist(tl L))

in (if restenttype (restopenarglist(tl L))<> nil 
andalso restenttype (restopenarglist(tl L))<> [""] 
then saypause ("Construction line not completely read:  "
^(hd (restenttype (restopenarglist(tl L))))) else ();

TextIO.output(!LOGFILE,!THELINE);
say0("Your command: "^(!THELINE2));Construct s L1 T) end


else if a = "define" andalso L<>nil andalso tl L<>nil 
then let val s = hd L and L1 = getopenarglist (tl L) 
and T = deent(getterms(restopenarglist(tl L)))

in if T = Unknown then 
saypause "Sorry, cannot define something as an abstraction" else

(if restterms(restopenarglist(tl L)) <> nil 
andalso restterms(restopenarglist(tl L)) <> [""] 
then saypause ("Definition line not completely read:  "^(hd(restopenarglist(tl L))))else ();
TextIO.output(!LOGFILE,!THELINE);say0("Your command: "^(!THELINE));Define s L1 T) end

else if a = "rewritec" andalso L <> nil andalso tl L <> nil 
then let val s = hd L and L1 = getopenarglist (tl L) 
and V = Hd(restopenarglist(tl L)) in 

(if (Tl(restopenarglist(tl L))) <> nil 
andalso (restopenarglist(tl L)) <> [""] 
then saypause ("Rewrite construction line not completely read:  "
^(hd(restterms(restopenarglist(tl L)))))else ();TextIO.output(!LOGFILE,!THELINE);
say0("Your command: "^(!THELINE));Rewritec s L1 V) end

else if a = "rewrited" andalso L <> nil andalso tl L <> nil 
then let val s = hd L and L1 = getopenarglist (tl L) 
and V = Hd(restopenarglist(tl L)) in 

(if (Tl(restopenarglist(tl L))) <> nil 
andalso (restopenarglist(tl L)) <> [""] 
then saypause ("Rewrite demonstration line not completely read:  "
 ^(hd(restopenarglist(tl L))))else ();TextIO.output(!LOGFILE,!THELINE);
say0("Your command: "^(!THELINE));Rewrited s L1 V) end

else if a = "open" then (TextIO.output(!LOGFILE,!THELINE^"\n");
say0 ("Your command:  "^(!THELINE)^"\n"); 
Open(if L=nil orelse hd L = "" then makestring(length(!CONTEXT)) else hd L))

else if a = "close" then (TextIO.output(!LOGFILE,!THELINE^"\n")
;say0("Your command: "^(!THELINE)^"\n");Close())

else if a = "save" then (TextIO.output(!LOGFILE,!THELINE^"\n")
;say0 ("Your command:  "^(!THELINE)^"\n"); Save(if L=nil orelse hd L = "" then
hd(!WORLDNAMES) else hd L))

else if a = "load" then (TextIO.output(!LOGFILE,!THELINE^"\n")
;say0 ("Your command:  "^(!THELINE)^"\n"); LoadTheory(if L=nil orelse hd L = "" then
"" else hd L))

else if a = "showall" then showall()

else if a = "showimplicit" then showimplicit()

else if a = "hideimplicit" then hideimplicit()

else if a = "displayrewrites" then displayrewrites()

else if a = "showrecent" then showrecent()

else if a = "showdec" andalso L<>nil then showdec (hd L)

else if a = "showdecs" then showdecs()

else if a = "foropen"  then say ("\n\n"^(savedforopen()))

else if a = "forclearcurrent"  then say ("\n\n"^(savedforclearcurrent()))

else if a="comment" orelse a="%" 
then (TextIO.output(!LOGFILE,!THELINE^"\n");say0((!THELINE^"\n")))

else if a="comment1" orelse a="%%" 
then (TextIO.output(!LOGFILE,!THELINE);say0((!THELINE)))

else if a =">> " then ()

else if a = "clearcurrent" then 
(ClearCurrent(if L = nil orelse hd L = "" then makestring(length(!CONTEXT)-1) else hd L);
TextIO.output(!LOGFILE,!THELINE^"\n"))

else if a = "clearall" then (ClearAll();TextIO.output(!LOGFILE,!THELINE);showall())

else if a = "basic" then (basic();TextIO.output(!LOGFILE,!THELINE))

else if a = "explicit" then (explicit();TextIO.output(!LOGFILE,!THELINE))

else if a = "fullversion" then (fullversion();TextIO.output(!LOGFILE,!THELINE))

else if a = "pause" then 
(say ("Pausing in "^(!LOGNAME)^":\n>>  type lines or type quit to resume");
TextIO.output(!LOGFILE,!THELINE);interface " ")

else if a = "" then () else saypause "Line is not a Lestrade command"

(* purge indentation from command lines *)

(* and  despace0 (#" "::L) = despace0 L |

despace0 L = L

and despace s = implode(despace0(explode s)) *)

(* read a command line from a string *)

and Readline s = (THELINE2:=(despace s)^"\n";
THELINE:=(INDENTS())^(despace s)^"\n";readline(tokenize s))

(* read command lines from standard output and receive feedback;
 output is logged to a file, end with quit *)

and interface filename =
(if filename = "" orelse filename = " " 
then () else LOGFILE:=TextIO.openOut((filename^".lti"));
 (if not(!GREETED) then (versiondate();GREETED:=true) else ();
let val LINE = TextIO.inputLine(TextIO.stdIn) in

if LINE = "quit\n" then (if filename <> " " 
then TextIO.output(!LOGFILE,"quit") else ();
TextIO.flushOut(!LOGFILE);
if filename <> " " 
then (TextIO.flushOut(!LOGFILE);closelog())
 else TextIO.flushOut(!LOGFILE);say "Bye!")

else (Readline (implode(rev(tl(rev(explode LINE)))));interface "") end))

(* read commands from a first log file after clearing the Lestrade context,
 logging to a second log file, and ending in the interface 
where you can continue to enter commands logged to the second file.  End with quit *)


and readfile filename1 filename2 =
(if filename1 <> "" then BREAKOUT:=false else ()
;if filename1 <> "" then THEORYNAME := filename1 else (); 
if filename1 = "" then () 
else (ClearAll();READFILE:=TextIO.openIn((filename1^".lti"))); 
if filename2 = "" then () 
else (ClearAll();LOGNAME:=filename1;LOGNAME2:=filename2;
LOGFILE:=TextIO.openOut((filename2^".lti"))); 
(if not(!GREETED) then (versiondate();GREETED:=true) else ();
let val LINE = TextIO.inputLine(!READFILE) in

if LINE = "quit\n" orelse (!BREAKOUT) then (BREAKOUT:=false;
TextIO.closeIn(!READFILE);
say("Done reading "^(!LOGNAME)^" to "^(!LOGNAME2)^":\n>>"
^"  type lines or type quit to exit interface\n\nquit\n"); 
SAVEDTHEORIES:= (!THEORYNAME,(!SERIAL,!NAMESERIAL,hd(rev(!CONTEXT))))
::(abstractdrop filename2 (!SAVEDTHEORIES)); if (!READFILEDEPTH)=0 then 
interface "" else READFILEDEPTH :=(!READFILEDEPTH)-1)

else (Readline (implode(rev(tl(rev(explode LINE)))));readfile "" "") end));


fun fixargtest s t n = (deworld(getabstype (pi1(hd(stringtype s)))),
guardedfixarglist (deworld(getabstype (pi1(hd(stringtype s))))) (getarglist n (tokenize t)));

(* disaster cleanup -- close the files if you crash out of the interface *)

fun Cleanup() = (TextIO.closeOut(!LOGFILE); TextIO.closeIn(!READFILE));

fun senttype s = entitytype(sent s);



\end{verbatim}}

\end{document}